\documentclass[11pt]{article}

\usepackage[dvips]{graphicx}
\usepackage{amssymb} 
\def\xt{X_{t_{k-1}}}

\def\yt{Y_{t_{k-1}}}
\def\ytk{Y_{t_k}}
\def\sumkn{\sum_{k=1}^n}

\def\hmc{\hat{\theta}_n}

\usepackage{amssymb}
\usepackage{amsthm}
\usepackage{latexsym}
\usepackage{amsmath}
\usepackage{color}
\usepackage{comment}
\usepackage{ascmac}
\newif\ifcolor
\colortrue

\includecomment{jp-text}
\excludecomment{en-text}
\newif\ifrs
\rstrue
\ifrs \usepackage{mathrsfs} \fi
\newif\ifcol
\coltrue

\newcommand{\cred}{\color{black}}
\newcommand{\cod}{\color{black}}
\newcommand{\cfg}{\color{black}}
\newcommand{\cdo}{\color{black}}
\newcommand{\colorred}{\color{black}}
\newcommand{\colorg}{\color{black}}
\newcommand{\colorblue}{\color{black}}
\newcommand{\colorn}{\color{black}}
\newcommand{\coloroy}{\color{black}}
\newcommand{\colorsb}{\color{black}}
\newcommand{\cRose}{\color{black}}
\newcommand{\chp}{\color{black}}
\newcommand{\colorr}{\color{black}}
\newcommand{\colorb}{\color{black}}
\newcommand{\bbA}{{\mathbb A}}

\newcommand{\bbD}{{\mathbb D}}

\newcommand{\bbH}{{\mathbb H}}

\newcommand{\bbL}{{\mathbb L}}

\newcommand{\bbN}{{\mathbb N}}

\newcommand{\bbQ}{{\mathbb Q}}
\newcommand{\bbR}{{\mathbb R}}
\newcommand{\bbS}{{\mathbb S}}

\newcommand{\bbX}{{\mathbb X}}
\newcommand{\bbY}{{\mathbb Y}}
\newcommand{\bbZ}{{\mathbb Z}}

\def\koko{{\coloroy{koko}}}
\def\D2{\bbD_{2,\infty-}}

\usepackage{mathrsfs} 
\usepackage{graphicx} 
\usepackage{layout} 

\newtheorem{lemma}{Lemma}
\newtheorem{proposition}{Proposition}
\newtheorem{theorem}{Theorem}
\newtheorem{remark}{Remark}
\newtheorem{example}{Example}

\def\D{{\bf D}}

\def\calc{{\cal C}}

\def\cale{{\cal E}}
\def\calf{{\cal F}}
\def\calg{{\cal G}}

\def\calj{{\cal J}}

\def\call{{\cal L}}

\def\caln{{\cal N}}

\def\calp{{\cal P}}

\def\calu{{\cal U}}

\def\calx{{\cal X}}

\def\sskip{\hspace{0.5cm}}

\def\ep{\epsilon}
\def\half{\frac{1}{2}}

\def\Iku{\Rightarrow}

\def\y{\vspace*{3mm}\\}
\def\nn{\nonumber}
\def\be{\begin{equation}}
\def\ee{\end{equation}}
\def\bea{\begin{eqnarray}}
\def\eea{\end{eqnarray}}
\def\beas{\begin{eqnarray*}}
\def\eeas{\end{eqnarray*}}
\def\l{\left}
\def\r{\right}

\begin{document}
\title{
{\colorb{Nondegeneracy of Random Field and Estimation of Diffusion
}}
}
\author{
 ${}^{1,2}${Masayuki Uchida}
  	   and
 ${}^{3,4}${Nakahiro Yoshida}
 \\
        ${}^1${\small {Graduate School of Engineering Science, 
                       Osaka University}}\\
        {\small {Toyonaka, Osaka 560-8531, Japan}}
        {\small and} ${}^2${\small Japan Science and Technology Agency, CREST}\\
        ${}^3${\small {Graduate School of Mathematical Sciences, 
                       University of Tokyo}}\\
        {\small {3-8-1 Komaba, Meguro-ku, Tokyo 153-8914, Japan}} 
        {\small and} ${}^4${\small The Institute of Statistical Mathematics}
}
\date{December 9, 2012}
\maketitle
\noindent
{\bf Abstract.} We construct a quasi likelihood analysis for diffusions 
under the high-frequency sampling over a finite time interval. 
For this, we prove a polynomial type large deviation inequality for the quasi likelihood random field. 
Then it becomes crucial to prove nondegeneracy of 
a key index $\chi_0$. 
By nature of the sampling setting, 
$\chi_0$ is random. This makes it difficult to apply 
a na\"ive sufficient condition, 
and requires a new machinery. 
In order to establish a quasi likelihood analysis, we need quantitative estimate of 
the nondegeneracy of $\chi_0$.
{\cRose{
The existence of a nondegenerate local section of a certain tensor bundle associated with the statistical random field solves this problem.
}}

\vspace{5pt} 


\vspace{5pt} 

\vspace{5pt}

{\colorb{\noindent
}}


\section{Introduction}

\begin{en-text}
In order to study asymptotic properties of statistics in likelihood analysis,
Ibragimov and Has'minskii (1972, 1973, 1981) established a new paradigm of likelihood analysis.
Let ${\cal E}^\epsilon = \{ {\cal X}^\epsilon, {\cal A}^\epsilon, (P_{\theta}^\epsilon)_{\theta \in \Theta} \}$
be a sequence of statistical experiments, where $\epsilon \in (0,1]$ and 
a parameter space $\Theta \subset {\bbR}^{{\sf p}}$.
For a $\theta^* \in \Theta$, let a statistical random field $Z_\epsilon$ based on the likelihood ratio
denote
$$
Z_\epsilon(u) = \frac{d P_{\theta^* + \phi(\epsilon)u}^\epsilon}{dP_{\theta^*}^\epsilon}(X^\epsilon)
$$
for $u \in \bbR^{{\sf p}}$, where $\phi(\epsilon)$ is a positive normalizing factor tending to zero as $\epsilon \rightarrow 0$.
Let $\hat{C}(\bbR^{{\sf p}})$ be the space of continuous functions on $\bbR^{{\sf p}}$ that
tends to zero at the infinity.
Let $C_\uparrow(\bbR^{{\sf p}})$ be the space of continuous functions on $\bbR^{{\sf p}}$ of
at most polynomial growth.
A simplified version of their result is as follows.

\begin{theorem}[Ibragimov and Has'minskii (1972, 1973, 1981)]
Suppose that $Z_\epsilon$ meets the following conditions.

(i) There exist $\alpha >p$ and $k \geq \alpha$ such that for some constant $C>0$,
$$
E_{\theta^*}^\epsilon \left[ \left| Z_\epsilon (u_2)^{1/k} - Z_\epsilon (u_1)^{1/k} \right|^k \right] \leq C | u_2 - u_1|^\alpha 
\quad \mbox{for all } u_1, u_2, \epsilon.
$$

(ii) For some $\gamma > 0$ and $c >0$,
\begin{equation} \label{I-H-2}
P_{\theta^*}^\epsilon \left[ Z_\epsilon (u) \geq e^{-c |u|^\gamma} \right] \leq e^{-c |u|^\gamma}.
\end{equation} 

(iii) Finite-dimensional convergence: $Z_\epsilon \rightarrow^{d_f} Z$, where
$Z$ is a $\hat{C}(\bbR^{{\sf p}})$-valued random variable.

Then, $(P_{\theta^*}^\epsilon)^{Z_\epsilon} \rightarrow {\cal L}\{ Z \} $. Moreover,
$$
P_{\theta^*}^\epsilon \left[ \sup_{u : |u| \geq r} Z_\epsilon (u) \geq e^{-c_1 r^\gamma} \right] \leq e^{-c_1 r^\gamma}.
$$

If $\hat{u}$ uniquely attains the maximum of $Z(u)$, then for any sequence of the maximum likelihood estimator $\hat{\theta}_\epsilon$
for $\theta$, $\hat{u}_\epsilon := \phi(\epsilon)^{-1}(\hat{\theta}_\epsilon - \theta^*) \rightarrow^d \hat{u}$
and
$$
E_{\theta^*}^\epsilon \left[ f(\hat{u}_\epsilon) \right] \rightarrow E[f(\hat{u})]
$$
for all $f \in C_\uparrow(\bbR^{{\sf p}})$.
\end{theorem}

Since it follows from the convergence $Z_\epsilon \rightarrow^d Z$ in $\hat{C}(\bbR^{{\sf p}})$
that $sup_{u \in A} Z_\epsilon \rightarrow^d sup_{u \in A} Z$ for every measurable set $A$,
one has that $\hat{u}_\epsilon = \arg \max_u Z_\epsilon \rightarrow^d \hat{u} = \arg \max_u Z$.
In the local asymptotic mixed normality (LAMN) case, 
$Z(u) = \exp \{ \Gamma(\theta^*)^{1/2} \zeta[u] - 2^{-1} \Gamma(\theta^*)[u,u] \}$,
where $\Gamma(\theta^*)$ is the random Fisher information matrix, $I_p$ is the $p \times p$ identity matrix 
and $\zeta \sim N_p(0,I_p)$ independent of $\Gamma(\theta^*)$.
In this case, the asymptotic mixed normality of the maximum likelihood estimator is obvious from the convergence
$\hat{u}_\epsilon \rightarrow^d \hat{u}$ since $\hat{u}= \Gamma(\theta^*)^{-1/2}\zeta$. 
This is the first advantage of the Ibragimov and Has'minskii method.
In order to show the convergence $Z_\epsilon \rightarrow^d Z$ in $\hat{C}(\bbR^{{\sf p}})$, however,
it is not necessary to assume such a strong condition of separation as (\ref{I-H-2}), see Yoshida (1990).
In this sense, 
the most essential part of their method is the use of 
the large deviation estimate (\ref{I-H-2})
and strong properties such as the convergence of moments of the estimator
are derived from it.

Kutoyants (1984, 1994, 1998, 2004) applied the Ibragimov-Has'minskii approach to stochastic processes
including diffusion type processes and point processes.
However, there is a serious point that the large deviation inequality (\ref{I-H-2}) is not easy to obtain for most of stochastic processes.
On the other hand, Yoshida (1990) proved the convergence of statistical random fields without such a strong assumption.
As an important observation, 
the exponential type large deviation inequality like (\ref{I-H-2}) is much stronger than our use.
In order to develop a theory, 
it is sufficient to obtain the following polynomial type large deviation inequality.
$$
P_{\theta^*}^\epsilon \left[ \sup_{u : |u| \geq r} Z_\epsilon (u) \geq r^{-N} \right] \leq \frac{C_N}{r^N}.
$$
Here we note that some exponential in place of $r^{-N}$ on the left-hand side is very often possible.
In particular, Kutoyants (2004) presented a polynomial type large deviation inequality for
one-dimensional diffusion process by means of the local time.
Another serious point is that we do not generally have an explicit form of likelihood function 
for stochastic differential equations from discrete observations.
Yoshida (2005) obtained the polynomial type large deviation inequality
for an abstract statistical random field based on the partially local asymptotically quadratic (PLAQ)
sequence of experiments.
This approach enables us to connect the Ibragimov-Has'minskii-Kutoyants program in the asymptotic decision theory for
stochastic processes to the quasi-likelihood analysis for discretely sampled continuous-time stochastic processes. 
As an example, for multi-dimensional mixing diffusion processes with unknown parameters
in both drift and diffusion coefficient,
the convergence of moments of quasi-maximum likelihood estimator, asymptotic normality of Bayes estimator
and the convergence of moments of it were shown in Yoshida (2005).
\end{en-text}

{\colorb{
In this paper, we consider estimation for a stochastic regression model specified by the 
stochastic integral equation
\begin{eqnarray} \label{sde1}
  Y_t &=& Y_0+\int_0^t b_s ds+ \int_0^t \sigma (X_s, \theta )dw_s, \quad t \in [0,T], 
\end{eqnarray}
where 
$w$ is an ${\sf r}$-dimensional standard Wiener process on a stochastic basis $(\Omega, {\cal F}, ({\cal F}_t)_{t \in [0,T]}, P)$,
$b$ and $X$ are progressively measurable processes with values in $\bbR^{{\sf m}}$ and ${\bbR}^{\sf d}$, respectively,
$\sigma$ is an $\bbR^{{\sf m}} \otimes \bbR^r$-valued function defined on $\bbR^{\sf d} \times \Theta$,
and
$\Theta$ is a bounded domain in $\bbR^{{\sf p}}$. 
As a special case, if an argument of $X_t$ is $t$, then the volatility in the model (\ref{sde1}) is time dependent. 
Furthermore, if we set $b_t = b(Y_t,t)$ and $X_t =(Y_t,t)$, then $Y$ can be a time-inhomogeneous diffusion process.
Of course, the stochastic volatility model like (\ref{sde1}) is quite commonly used in finance and econometrics. 
The data set consists of discrete observations ${\bf Z}_n = (X_{t_k},Y_{t_k})_{0 \leq k \leq n}$ 
with $t_k = k h$ for $h=h_n = T/n$. 
The process $b$ is completely unobservable and unknown. 
The asymptotics will be considered for $n\to\infty$, that is,
${\bf Z}_n$ forms high frequency data. 
 
Asymptotic theory of parametric estimation for the unknown parameter $\theta$ in the volatility 
of the stochastic differential equation 
based on high frequency data has been developed. 
Among many studies in a long history, we refer the reader to 
Plakasa Rao (1983,1988), 
Yoshida (1992,2005), 
Kessler (1997) under ergodicity, 
Shimizu and Yoshida (2006), Shimizu (2006), 
Ogihara and Yoshida (2009) for jump diffusion processes, 
Sorensen and Uchida (2003), Uchida (2003, 2004, 2008) for perturbed diffusions, 
Dohnal (1987), Genon-Catalot and Jacod (1993, 1994), Gobet (2001) 
for the fixed interval case. 
The limit distribution of the score function becomes a mixture of normal distributions 
over a finite time interval (LAMN), and a normal distribution 
over the infinite time interval (LAN) by the averaging effect. 
In this article, we will consider the LAMN (i.e., locally asymptotically mixed normal) 
quasi likelihood experiment associated with the sampling scheme over a finite time interval. 

A highlight of asymptotic decision theory is 
{\colorred{the likelihood analysis, the basic frame and functions of which} }
were established 
by Le Cam, H\'ajek, Ibra{\colorred{g}}imov and Has'minskii and others. The theory of Ibragimov and Has'minskii provides 
convergence of 
likelihood ratio random field on a function space with certain estimates for the tail probability and consequently  
convergence of moments of the estimator appearing in the likelihood analysis. 
It was Yury Kutoyants who found this methodology was effective for semimartingales, 
proving the wide applicability to various stochastic models. See Kutoyants (1984, 1994, 1998, 2004) for more information. 

Limiting distribution of the estimator is indispensable, however, it is far from sufficient to develop the elementary  
statistical theory. It is clear if we consider a problem of model selection, for example. 
The basic correction term by Akaike 
was introduced to make the Kullback-Leibler divergence between the predictive distribution and the true distribution 
asymptotically unbiased. Obviously, it is necessary to validate the existence of moments of the standardized estimator 
because the bias is described with it. The asymptotic distribution cannot provide sufficient information there. 
It is also the case in the prediction theory. Furthermore, 
the same kind of questions inevitably arise in the theory of higher-order statistical inference. 
Large deviation type estimates enable valid treatments of the higher-oder terms 
in the stochastic expansion of a statistic, and such estimates can be obtained 
by precise probabilistic estimate of the decay of the accompanying statistical random field. 

The {\bf quasi likelihood analysis} has been developing for stochastic processes. 
Here the quasi likelihood analysis means a system that gives asymptotic behavior of the quasi likelihood random field, 
its (polynomial type) large deviation estimate, 
limit theorems for the quasi maximum likelihood estimator and the quasi Bayesian estimator, and 
convergence of moments of these estimators. 
Yoshida (2005, 2011) gave a polynomial type large deviation inequality in the locally asymptotically quadratic 
(LAQ) setting to carry out the Ibragimov-Has'minskii-Kutoyants scheme for stochastic processes. 
As a corollary, the quasi-likelihood analysis for ergodic diffusion processes under sampling was presented. 
The simultaneous and adaptive Bayesian estimators were defined there. 
See Le Cam (1986), Le Cam and Yang (1990) for the fundamental notions of statistical experiments and approximation. 

The polynomial type large deviation inequality works in various settings.  
Uchida (2010) 
considered a model selection problem for discretely observed ergodic multi-dimensional diffusion processes 
and proposed a contrast-based information criterion. 
The difficulties are in existence of moments, and besides, in handling 
the exact likelihood function, that has no explicit expression. 
The polynomial type large deviation inequality and the Malliavin calculus were effectively used. 
The asymptotic results can be fairly complicated if jumps with heavy tail are involved; 
even convergence rate of the estimator can 
differ from the standard one. 
Masuda (2010) 
obtained a polynomial type large deviation estimate for the random field associated with 
a general self-weighted least absolute deviation (SLAD) in the parameter estimation of sampled 
Ornstein-Ulenbeck process driven by a heavy-tailed symmetric L\'evy process with positive 
activity index, and clarified asymptotic behavior of the estimator including convergence of moments.  
A quasi likelihood analysis was constructed by Ogihara and Yoshida (2009) for a nonlinear 
sampled diffusion process with jumps with the aid of the polynomial type large deviation inequality. 

Against these backgrounds, the first aim of this article is to construct a quasi likelihood analysis for diffusions 
under the high-frequency sampling over a finite time interval. 
For this, we will prove a polynomial type large deviation inequality for the quasi likelihood random field. 
Then we meet a question of nondegeneracy of 
a key index $\chi_0$ given in (\ref{230806-1}). 
By nature of the sampling setting, $\chi_0$ is random and this makes it difficult to apply 
a na\"ive sufficient condition often used so far because our model can easily break it. 
In order to establish a quasi likelihood analysis, we need quantitative estimate of 
the nondegeneracy of $\chi_0$.  
{\cRose{This problem is solved by the existence of a nondegenerate local section of a certain tensor bundle related to 
the statistical random field.}}
This is the second aim of this paper. 
Since such nondegeneracy argument is universal, the authors hope this part has its own interest 
even apart from statistical results presented here. 

}}

\begin{en-text}
We consider the estimation problem for 
an $m$-dimensional It{\^o} process  
satisfying the stochastic differential equation
\begin{eqnarray} \label{sde1bis}
  dY_t &=& b_t dt+ \sigma (X_t, \theta )dw_t, \quad t \in [0,T], 
\end{eqnarray}
where 
$w$ is an ${\sf r}$-dimensional standard Wiener process on some stochastic basis $(\Omega, {\cal F}, ({\cal F}_t)_{t \in [0,T]}, P)$,
$b$ and $X$ are progressively measurable processes with values in $\bbR^{{\sf m}}$ and ${\bbR}^{\sf d}$, respectively,
$\sigma$ is an $\bbR^{{\sf m}} \otimes \bbR^{\sf r}$-valued function defined on $\bbR^{\sf d} \times \Theta$,
and
$\Theta$ is 
a bounded domain in $\bbR^{{\sf p}}$
with a locally Lipschitz boundary,
which means that $\Theta$ has the strong local Lipschitz condition,
see Adams (1975) and Adams and Fournier (2003). 
$\theta^*$ denotes the true value of $\theta$.
The data are discrete observations ${\bf Z}_n = (X_{t_k},Y_{t_k})_{0 \leq k \leq n}$ 
with $t_k = k h$ for $h=h_n = T/n$, that is,
${\bf Z}_n$ forms high frequency data.
Note that if an argument of $X_t$ is $t$, then the volatility in the model (\ref{sde1}) is time dependent.
Furthermore, if we set that $b_t = b(Y_t,t)$ and $X_t =(Y_t,t)$, then under some regularity conditions,
$Y$ is the time-inhomogeneous diffusion process.
 
Genon-Catalot and Jacod (1993, 1994) proposed contrast functions for diffusion processes
and they proved the asymptotic mixed normality of the minimum contrast estimator.

In this paper, asymptotic mixed normality and convergence of moments of both the maximum contrast estimator and the Bayes estimator 
for the unknown parameter $\theta$ in the volatility of the stochastic differential equation (\ref{sde1})
are shown by means of the Ibragimov-Has'minskii-Kutoyants scheme.
The key point is to obtain the polynomial type large deviation inequality for the statistical random field.
\end{en-text}



{\colorb{
\section{Quasi likelihood analysis for diffusion and the limit theorems}

In this section, we will present the main results in statistical context. 

Suppose that $\Theta$ is 
a bounded domain in $\bbR^{{\sf p}}$
with a locally Lipschitz boundary, 
which means that $\Theta$ has the strong local Lipschitz condition,
see Adams (1975) and Adams and Fournier (2003). 
$\theta^*$ denotes the true value of $\theta$.

Let $C_\uparrow^{k,l}(\bbR^{\sf d} \times \Theta; \bbR^{{\sf m}})$ denote
the space of all functions $f$ satisfying the following conditions:
(i) $f(x,\theta)$ is an $\bbR^{{\sf m}}$-valued function on $\bbR^{\sf d} \times \Theta$,
(ii) $f(x,\theta)$ is continuously differentiable with respect to $x$ up to order $k$
for all $\theta$, and their derivatives up to order $k$ are
of polynomial growth in $x$ uniformly in $\theta$.
(iii) for $|{\bf n}|=0,1, \ldots, k$, $\partial_x^{{\bf n}} f(x,\theta)$ is continuously differentiable
with respect to $\theta$ up to order $l$ for all $x$.
Moreover, for $|{\bf \nu}|=0,1, \ldots, l$ and $|{\bf n}|=0,1,\ldots, k$, 
$\partial_\theta^{\bf \nu} \partial_x^{\bf n} f(x,\theta)$ is of polynomial growth in $x$ uniformly in $\theta$.
Here ${\bf n} =({ n}_1, \ldots, { n}_{\sf d})$ and ${\bf \nu} =({\bf \nu}_1, \ldots, {\bf \nu}_{\cRose{\sf p}})$ are multi-indices,
${\cRose{\sf p}} = \mbox{dim}(\Theta)$, $|{\bf n}| = { n}_1+ \ldots + { n}_{\sf d}$, 
$|{\bf \nu}| = {\bf \nu}_1+ \ldots + {\bf \nu}_{\cRose{\sf p}}$,
$\partial_x^{{\bf n}} = \partial_{x_1}^{{ n}_1} \cdots \partial_{x_{\sf d}}^{{ n}_{\sf d}}$, $\partial_{x_i} = \partial/\partial x_i$,
and $\partial_\theta^{\bf \nu} = \partial_{\theta_1}^{{\bf \nu}_1} \cdots \partial_{\theta_{\cRose{\sf p}}}^{{\bf \nu}_{\cRose{\sf p}}}$, 
$\partial_{\theta_i} = \partial/\partial \theta_i$.
We denote by $\rightarrow^p$ and $\rightarrow^{d_s({\cal F})}$ the convergence in probability 
and the ${\cal F}$-stable convergence in distribution,
respectively.
For matrices $A$ and $B$ of the same size,
we write $A^{\otimes 2} = A A^\star$
and $A[B] = \mbox{Tr}(AB^\star)$, where $\star$ means the transpose. 
Set $S(x, \theta) = \sigma(x,\theta)^{\otimes 2}$ and $\Delta_k Y= \ytk-\yt$.
%
We assume that 
{\cred the function $\sigma$ admits a continuous extension over $\bbR^{\sf d}\times\bar{\Theta}$}, 
and 
denote it by $\sigma$. 
{\cod Let}
\bea\label{230807-1}
Q(x,\theta,\theta^*)
&=&
{\rm Tr}\bigg(S(x,\theta)^{-1}
S(x,\theta^*)-I_d\bigg)
-\log\det\bigg(S(x,\theta)^{-1}S(x,\theta^*)\bigg) .
\eea

We consider the following conditions. 
\begin{description}
\item[[A1\!\!]] {\bf (i)} 
$\sup_{0 \leq t \leq T} \|b_t\|_p < \infty$ for all $p>1$.
\begin{description}
\item[(ii)] 
$\sigma \in C_\uparrow^{2,4}(\bbR^{\sf d} \times \Theta; \bbR^{{\sf m}} \otimes \bbR^{\sf r})$
and $\inf_{x, \theta} \det S(x,\theta) > 0$. 
\end{description}
\end{description}

\begin{description}
\item[[A2\!\!]] 
The process $X$ admits a representation 
\begin{eqnarray*}
X_t &=& X_0 + \int_0^t \tilde{b}_s ds + \int_0^t a_{s} dw_s + \int_0^t \tilde{a}_{s} d\tilde{w}_s, 
\end{eqnarray*}
where 
\begin{description}
\item[(i)] 
$\tilde{b}$, $a$ and $\tilde{a}$ are 
progressively measurable processes taking values in $\bbR^{\sf d}$, 
${\bbR}^{\sf d} \otimes {\bbR}^{r}$ and ${\bbR}^{\sf d} \otimes {\bbR}^{\sf r_1}$, respectively, 
satisfying 
\beas 
{\cred \|X_0\|_p+}
\sup_{t\in[0,T]}(\|\tilde{b}_t\|_p +\|a_t\|_p+\|\tilde{a}_t\|_p)<\infty
\eeas
for every $p>1$, 
\begin{en-text}
\beas 
\sup_{t\in[0,T]}(|\tilde{b}_t| +|a_t|+|\tilde{a}_t|)\in\cap_{p\geq1}L^p,
\eeas
\end{en-text}
and 
$\tilde{w}$ is an ${\sf r_1}$-dimensional Wiener process independent of $w$,
\item[(ii)] 
there is a stopping time $\tau$ such that $\mbox{ess.sup}_{\omega\in\Omega}\tau<T$ , 
$a^{\otimes2}_\tau+\tilde{a}^{\otimes2}_\tau$ is {\colorred bounded}, 
nondegenerate uniformly in $\omega\in\Omega$ and 
that $a$ and $\tilde{a}$ are 
right-continuous at $t=\tau$. 
\end{description} 

\begin{en-text}
$\Sigma_\tau$ is nondegenerate uniformly in $\omega\in\Omega$ and 
$\Sigma$ is right-continuous at $t=\tau$, 
where 
$\Sigma_t $ is an $\bbR^{{\sf d}}\otimes \bbR^{{\sf d}}$-valued progressively measurable process satisfying 
$\Sigma_t\Sigma_t^\star= a_t a_t^\star+\tilde{a}_t\tilde{a}^\star_t$. 
\end{description}
\end{en-text}

\end{description}


{\cfg 
We say that a function $f$ admits a $C^J$-supporting function at $(x_0,\theta_0)$ 
if there exist a function $g$ on a neighborhood $V(x_0,\theta_0){\cred \subset\bbR^{\sf d}\times\bar{\Theta}}$ 
of $(x_0,\theta_0)$ and $\xi_0\in\bbR^{\sf d}$, $|\xi_0|=1$,  
such that 
{\cred the partial derivatives $\partial_x^jg$ $(j=0,...,J)$ 
exists for each $\theta$ near $\theta_0$ and continuous in $(x,\theta)$} 
and that 
$
|f(x,\theta)| \geq |g(P_{\xi_0} x,\theta)|
$
for $(x,\theta)\in V(x_0,\theta_0)$, where $P_{\xi_0}$ is {\cred the} projection on $\bbR\xi_0$. 
Let $c_j(x,\theta)={\cod (j!)^{-1}(\partial_x^jg)(P_{\xi_0}x,\theta)[\xi_0^{\otimes j}]}$. [$g$ and $c_j$ depend on $(x_0,\theta_0)$.] 
}

\begin{description}
\item[[A3\!\!]] 
$\mbox{supp}\call\{X_\tau\}$ is compact,  and for some open neighborhood $U$ of $\mbox{supp}\call\{X_\tau\}$, 
there exist a function {\colorb{$f:U\times\bar{\Theta}\to\bbR$}} 
and a constant 
$\varrho\in(0,\infty)$ satisfying the following conditions. 
\begin{description}
\item[(i)] 
$
Q(x,\theta,\theta^*)|\theta-\theta^*|^{-2}
\geq |f(x,\theta)|^\varrho
$ 
for all $(x,\theta)\in U\times{\cred (}\bar{\Theta}{\cred \setminus\{\theta^*\})}$. 
%
%
%
\item[(ii)] 
 {\cfg $f$ admits a $C^J$-supporting function for each 
$(x_0,\theta)\in{\cred U\times\bar{\Theta}}$ with 
\beas \max_{j=0,...,J-1}\big|c_j(x_0,\theta) \big|>0.\eeas} 
%
\end{description}
\end{description}

Condition $[A1]$ is for regularity, $[A2]$ is for the nondegeneracy of the process $X$. 
The stopping time $\tau$ is often taken as $\tau=0$. 
The compactness of the support of $\call\{X_\tau\}$ can be relaxed if we assume 
stronger global nondegeneracy; a stronger condition will be inevitable in general because 
degeneracy can occur unless we assume the compactness of the support. 
Condition $[A3]$ is for the nondegeneracy of the quasi likelihood random field to which 
the nondegeneracy of $X$ can be conveyed thanks to the condition.



Since the exact transition density is not available, the inference is carried out by 
a quasi likelihood function. Let 
\beas 
{\mathbb H}_n(\theta) = -\frac{n{\cRose{\sf m}}}{2} \log ( 2 \pi h ) -\frac{1}{2} \sumkn \left\{ \log \det S(\xt, \theta) 
+ h^{-1}  S^{-1}(\xt,\theta) [ (\Delta_k Y)^{\otimes 2}] \right\}.
\eeas
Then the maximum likelihood type estimator $\hmc$ is any estimator that satisfies 
\begin{eqnarray} \label{M-est}
{\mathbb H}_n (\hmc) = \sup_{\theta \in \Theta} {\mathbb H}_n(\theta).
\end{eqnarray}
The Bayes type estimator $\tilde{\theta}_n$ for a prior density $\pi : \Theta \rightarrow {\bbR}_+$ 
with respect to the quadratic loss is defined by 
\begin{eqnarray} \label{Bayes1}
\tilde{\theta}_n = \left( \int_{\Theta} \exp ( {\mathbb H}_n(\theta) ) \pi(\theta) d\theta \right)^{-1}
\int_{\Theta} \theta \exp ( {\mathbb H}_n(\theta) ) \pi(\theta) d\theta. 
\end{eqnarray}
We assume that $\pi$ is continuous and $0 < \inf_{\theta \in \Theta} \pi(\theta) \leq \sup_{\theta \in \Theta} \pi(\theta) <\infty$.

Let 
$\Gamma(\theta^*) = (\Gamma^{ij}(\theta^*))_{i,j=1,\ldots, {\cRose{\sf p}}}$ with
\begin{eqnarray*}
\Gamma^{ij}(\theta^*) &=& \frac{1}{2 T} \int_0^T \mbox{Tr} \left( (\partial_{\theta_i} S) S^{-1}(\partial_{\theta_j} S) S^{-1}(X_t,\theta^*) \right) dt
\end{eqnarray*}
and let 
$\zeta$ be a ${\cRose{\sf p}}$-dimensional standard normal random variable
independent of $\Gamma(\theta^*)$. 
Here are restricted versions of the main results in this article. 
We will give proof of these results in Section \ref{230806-3}. 
\begin{theorem}\label{230614-1}
Suppose that Conditions $[A1]$-$[A3]$ are fulfilled. Then 
$\sqrt{n} ( \hmc -\theta^*) \rightarrow^{d_s({\cal F})} \Gamma(\theta^*)^{-1/2} \zeta$ and
$$
E \left[ f(\sqrt{n} ( \hmc -\theta^*) ) \right] \rightarrow {\mathbb E} \left[ f(\Gamma(\theta^*)^{-1/2} \zeta) \right]
$$
as $n \rightarrow \infty$ for all continuous functions $f$ of at most polynomial growth.
\end{theorem}

\begin{theorem}\label{230614-2}
Suppose that Conditions $[A1]$-$[A3]$ are fulfilled. Then 
$\sqrt{n} ( \tilde{\theta}_n -\theta^*) \rightarrow^{d_s({\cal F})} \Gamma(\theta^*)^{-1/2} \zeta$ and
$$
E \left[ f(\sqrt{n} (  \tilde{\theta}_n -\theta^*) ) \right] \rightarrow {\mathbb E} \left[ f(\Gamma(\theta^*)^{-1/2} \zeta) \right]
$$
as $n \rightarrow \infty$ for all continuous functions $f$ of at most polynomial growth.
\end{theorem}


\vspace{3mm}
We shall consider Condition $[A3]$ {\cfg (ii)}. 
The following simple example suggests degeneracy of the statistical model can easily occurs. 
\begin{example}\rm Let $X_t$ satisfy the stochastic differential equation
\bea\label{230807-10}
dX_t &=& (1+X_t^2)^\theta dw_t,\sskip X_0=0. 
\eea
Assume $\Theta\subset(0,1/2]$. 
This model is completely unidentifiable when $t=0$. 
Now $S(x,\theta)=(1+x^2)^{2\theta}$ and 
\beas 
Q(x,\theta,\theta^*) {\colorred{|\theta -\theta^*|^{-2}}} &=& 
\left\{ \exp\big(2(\theta^*-\theta)\log(1+x^2)\big)-1-2(\theta^*-\theta)\log(1+x^2) \right\}{\colorred{|\theta -\theta^*|^{-2}}}
\\&\geq&
{\colorred{\big(\log(1+x^2)\big)^2}}. 
\eeas
{\colorred{ for $(x,\theta) \in U \times \bar{\Theta}$ and $U=(-1,1)$}}. 
{\colorred{
Set $f(x,\theta)=\log(1+x^2)$ }} and $\varrho=2$.
Then 
{\colorred{
\beas 
f(x_0,\theta) &=& \log(1+x_0^2),\hspace{5mm}
\partial_x f(x_0,\theta) =\frac{2 x_0}{1+x_0^2},
\\
\partial_x^2 f(x_0,\theta)&=&
\frac{2}{1+x_0^2}-\frac{4 x_0^2}{(1+x_0^2)^2}. 
\eeas
}}
Therefore, 
$\max_{j=0,1,2}|\partial_x^jf(x_0,\theta)|>0$ for each $(x_0,\theta)\in U\times\bar{\Theta}$. 
This is a rather simple case because we found $f$ independent of $\theta$. 
Thus 
$[A3]$ {\cfg (ii)} holds. Condition $[A2]$ (ii) is also obvious if we choose $\tau=0$. 

\end{example}

\begin{example}\rm Consider a trivial model
\beas 
X_t &=& \theta \big(w_{t\vee\tau_0}- w_{\tau_0}\big),\hspace{5mm}t\in[0,1].
\eeas
for a stopping time $\tau_0$. 
Then we should take $\tau=\tau_0$. This simple example suggests the necessity of introducing 
stopping time $\tau$. Naturally, $\mbox{ess.sup}_\omega\tau_0$ should be less than one for consistent estimation. 
\end{example}

\begin{en-text}
{\colorb{Suppose that $\mbox{supp}\call\{X_0\}$ is compact, as in Section \ref{230608-3}. 
Let $\calx_0=\mbox{supp}\{X_0\}$ and let $\hat{\calx}$ be a compact set in $\bbR^{\sf d}$ 
satisfying $\calx_0\subset\hat{\calx}^o$. 
Let $X=(X_t)$ be a $\sf d$-dimensional diffusion process satisfying 
the stochastic integral equation
\begin{eqnarray}
X_t &=& X_0+\int_0^t V_0(X_s)ds +\int_0^t V(X_s)dW_s, \quad t \in [0,T], \label{u-4}
\end{eqnarray}
where $W$ is an ${\sf r}$-dimensional standard Wiener process,
$V_0$ is an ${\bbR}^{\sf d}$-valued function defined on ${\bbR}^{\sf d}$
and
$V$ is an ${\bbR}^{\sf d} \otimes {\bbR}^{\sf r}$-valued function defined on ${\bbR}^{\sf d}$.
Assume that there exists $K>0$ such that for all $x,y$,
$$
|V_0(x)-V(y)| + |V(x)-V(y)| \leq K|x-y|
$$
and that $\inf_x \det([VV^\star](x)) >0$. 
\koko
}}
\end{en-text}

\begin{example}\rm 
Consider the stochastic differential equation 
\beas 
dX_t &=& S(X_t,\theta)^\half dw_t, \hspace{5mm}t\in[0,1]\\
X_0 &=& 0,
\eeas
where  
\beas 
S(x,\theta)
&=&
\exp\big(
\sin\theta\>\sin x-\theta^2\sin^2 x\big), 
\eeas
and 
$\theta\in\Theta=(-\pi,\pi)$. 
This model is completely degenerate at $t=0$. 
Let $\theta^*=0$. 
For small neighborhood $U$ of $0$, we have 
\beas 
&&
|S(x,\theta)^{-1}S(x,\theta^*)-1
-\log\{S(x,\theta)^{-1}S(x,\theta^*)\}||\theta-\theta^*|^{-2}
\\&\geq&
c|f(x,\theta)|^2
\eeas
for $(x,\theta)\in U\times\bar{\Theta}$, 
if we take some $c>0$ independent of $(x,\theta)$, and 
\beas 
f(x,\theta)
&=&
\theta^{-1}\big(
\sin\theta\>\sin x-\theta^2\sin^2 x\big). 
\eeas
With the function $f$, 
\beas 
f(x_0,\theta)
&=& 
\frac{\sin\theta}{\theta} \sin x_0-\theta\sin^2 x_0,
\\
\partial_xf(x_0,\theta)
&=& 
\frac{\sin\theta}{\theta}\cos x_0 
-2\theta\sin x_0\cos x_0,
\\
\half\partial_x^2f(x_0,\theta)
&=&
-\frac{\sin\theta}{2\theta}\sin x_0
-\theta\cos^2 x_0 +\theta\sin^2 x_0,
\eeas
the continuous extension being applied at $\theta=0$.  
%
%
Set $\Theta_1=\Theta\setminus\Theta_2$ with 
$\Theta_2=(B(\pi,\rho)\cup B(-\pi,\rho))\cap\Theta$. 
Then it is not difficult to fix $\rho>0$ and small $\ep>0$ 
so that $\partial_xf$ is nondegenerate on $\Theta_1$ and 
so is $\partial_x^2f$ on $\Theta_2$ 
in the same time. 

A na\"ive, simple-looking sufficient condition 
is that $\inf_x\inf_\theta |f(x,\theta)|>0$. 
However, in this example, $\inf_x |f(x,\theta)|=0$, and 
the na\"ive condition does not work. 
\end{example}

\vspace*{5mm}
The later sections will be devoted to the proof of Theorems  \ref{230614-1} and \ref{230614-2}.  
Some generalization will be done on the way. 
The ingredients of the proof of these results are the polynomial type large deviation inequality 
for the quasi likelihood random field, as well as limit theorems for semimartingales. 
In Section \ref{PLDP}, we recall the polynomial type large deviation inequality for the statistical random field. 
The aim of the section is to introduce a key random index $\chi_0$ associated with $\bbH_n$ and 
to clarify its role for derivation of the large deviation estimate and as a result 
for establishing the quasi likelihood analysis. 

Thus it is necessary to prove the nondegeneracy of a random index $\chi_0$. 
To answer this question, Section \ref{230806-5} is devoted to making a new machinery to induce 
the nondegeneracy of the statistical random field in a general manner 
by connecting nondegeneracy of the associated tensor fields over the statistical manifold 
and the nondegeneracy of the underlying stochastic process. 

After laying these foundations, we will return 
to the proof of Theorems \ref{230614-1} and \ref{230614-2}, in Section \ref{230806-3}.

\begin{en-text}
$[H']$ like condition plus $L^q$-boundedness of $\sigma_t$ etc.; finally we admit polynomial type estimates as well as exponential type. 
\end{en-text}


\section{Polynomial type large deviation inequality and a generalized quasi-likelihood analysis for diffusion} \label{PLDP}

Let ${\mathbb U}_n = \{ u \in \bbR^{{\sf p}} \ ; \ \theta^* +(1/\sqrt{n})u \in \Theta \}$
and
$V_n(r) =  \{ u \in {\mathbb U}_n \ ; \ r \leq |u| \}$.


We make the following assumption.

\begin{description}
\item[[H1\!\!]]
{\bf (i)} $E [ |X_0|^q ] < \infty$ for all $q>0$.
For every $q>0$, there exists $C>0$ such that 
\beas 
E [ |X_t -X_s|^q ] \leq C |t-s|^{q/2}
\eeas
for all $t,s  \in [0,T]$. \vspace{-2mm}
\begin{description}
\item[(ii)]
$\sup_{0 \leq t \leq T} E [ |b_t|^q ] < \infty$ for all $q>0$.
\item[(iii)] 
$\sigma \in C_\uparrow^{2,4}({\bbR}^{\sf d} \times \Theta; {\bbR}^{\sf m} \otimes {\bbR}^{\sf r})$
and $\inf_{x, \theta} \det S(x,\theta) > 0$.
\vspace{0.5cm}
\end{description}
\end{description}


We define the random field ${\mathbb Z}_n$ on ${\mathbb U}_n $ by 
\bea\label{230808-1} 
{\mathbb Z}_n(u) = \exp \left\{ {\mathbb H}_n \left( \theta^* + \frac{1}{\sqrt{n}} u \right) - {\mathbb H}_n (\theta^*) \right\}
\eea
for $u \in {\mathbb U}_n $. 
%
Let ${\mathbb Y}_n(\theta) = \frac{1}{n} \left\{ {\mathbb H}_n (\theta) - {\mathbb H}_n (\theta^*) \right\}$ and 
\bea\label{230806-1}
\chi_0 = \inf_{\theta \ne \theta^*} \frac{- {\mathbb Y}(\theta) }{|\theta -\theta^*|^2}, 
\eea
where
\begin{eqnarray*}
{\mathbb Y}(\theta) &=& -\frac{1}{2 T}\int_0^T \left\{ \log \left( \frac{\det S(X_t,\theta)}{\det S(X_t,\theta^*)} \right)
+\mbox{Tr} \left( S^{-1}(X_t,\theta) S(X_t,\theta^*) -I_d \right) \right\} dt.
\end{eqnarray*}


\vspace{0.5cm}
The following condition is concerning nondegeneracy of the index $\chi_0$. 
\begin{description}
\item[[H2\!\!]] For every $L>0$, there exists $c_L >0$ such that
\beas
P\left[ \chi_0 \leq r^{-1} \right] \leq \frac{c_L}{r^{L}}
\eeas
for all $r>0$. 
\end{description}


%

\begin{theorem} \label{thm1}
Assume $[H1]$ and $[H2]$. Then, for every $L>0$, there exists a positive constant $C_L$ such that
\beas 
P \left[ \sup_{u \in V_n(r)} {\mathbb Z}_n(u) \geq e^{-r} \right] \leq \frac{C_L}{r^L}
\eeas
for all $r>0$ and $n \in {\bf N}$.
\end{theorem}

Proof of Theorem \ref{thm1} is given in Section \ref{230605-1}. 
The above theorem clarifies the essential role of the random index $\chi_0$ 
because it gives the polynomial type large deviation estimate for $\bbZ_n$, from which 
all tail properties of the estimators are deduced as the theorems below.

In order to obtain the weak convergence of the statistical random field on compact sets,
we make the following assumption.

\vspace{0.5cm}
\begin{en-text}
\noindent
$[H1']$ Assumption $[H1]$ holds and
\begin{eqnarray*}
X_t &=& X_0 + \int_0^t \tilde{b}_s ds + \int_0^t a_{s} dw_s + \int_0^t \tilde{a}_{s} d\tilde{w}_s, 
\end{eqnarray*}
where 
$\tilde{b}$, $a$ and $\tilde{a}$ are 
locally bounded predictable with values in ${\bbR}^{\sf d}$,
${\bbR}^{\sf d} \otimes {\bbR}^{\sf r}$ and ${\bbR}^{\sf d} \otimes {\bbR}^{\sf r_1}$, respectively,
$b$ is locally bounded and
$\tilde{w}$ is an ${\sf r_1}$-dimensional Wiener process independent of $w$.
\end{en-text}
\begin{description}
{\cred \item[[H1$^\sharp$\!\!]] Conditions [A1] and [A2] (i) hold. }
\end{description} 




\begin{en-text}

\begin{lemma} \label{lem3-2} 
Assume $[H1']$. Then,
for every $R>0$, 
$ {\mathbb Z}_n (u) \rightarrow^{d_s({\cal F})} {\mathbb Z}(u)$
in $C(B(R))$ as $n \rightarrow \infty$.
\end{lemma}

\end{en-text}

The following theorems generalize Theorems \ref{230614-1} and \ref{230614-2}. 

\begin{theorem} \label{thm3}
Assume ${\cred [H1^\sharp]}$ and $[H2]$. Then,
$\sqrt{n} ( \hmc -\theta^*) \rightarrow^{d_s({\cal F})} \Gamma(\theta^*)^{-1/2} \zeta$ and
$$
E \left[ f(\sqrt{n} ( \hmc -\theta^*) ) \right] \rightarrow {\mathbb E} \left[ f(\Gamma(\theta^*)^{-1/2} \zeta) \right]
$$
as $n \rightarrow \infty$ for all continuous functions $f$ of at most polynomial growth.
\end{theorem}

\begin{theorem} \label{thm2}
Assume ${\cred [H1^\sharp]}$ and $[H2]$. Then,
$\sqrt{n} ( \tilde{\theta}_n -\theta^*) \rightarrow^{d_s({\cal F})} \Gamma(\theta^*)^{-1/2} \zeta$ and
$$
E \left[ f(\sqrt{n} ( \tilde{\theta}_n -\theta^*) ) \right] 
\rightarrow {\mathbb E} \left[ f(\Gamma(\theta^*)^{-1/2} \zeta) \right]
$$
as $n \rightarrow \infty$ for all continuous functions $f$ of at most polynomial growth.
\end{theorem}

Proof of Theorems \ref{thm3} and \ref{thm2} is given in Section \ref{230605-1}. 
}}

{\colorb{
A sufficient condition for [H2] is that 
\beas 
\inf_{\omega\in\Omega,\theta\in\Theta\setminus\{{\cRose{\theta^*}}\}\atop t\in[0,T]}
\left\{ \log \left( \frac{\det S(X_t,\theta)}{\det S(X_t,\theta^*)} \right)
+\mbox{Tr} \left( S^{-1}(X_t,\theta) S(X_t,\theta^*) -I_d \right) \right\}/|\theta-\theta^*|^2>0\sskip a.s.
\eeas 
Though this kind of condition 
seems easy to handle and at hand, 
it is too na\"ive as it breaks, for example, in 
a simple model such as (\ref{230807-10}).  

The nondegeneracy condition [H2] of the statistical random field 
is a key to construction of quasi-likelihood analysis for diffusion. As we saw above, 
once the nondegeneracy of the index $\chi_0$ is established, we can obtain limit theorems 
for the quasi maximum likelihood estimator and the Bayesian type estimator, 
and moreover convergence of moments of them.

It should be remarked that the limit theorem for the Bayesian type estimator and convergence of moments of these estimators are new, and that the latter is indispensable to practical applications such as model selection, prediction and theory of asymptotic expansion. 
We will pursuit this nondegeneracy problem for statistical random fields. 
The question is when Condition [H2] holds. 
We discuss this problem in Section \ref{230806-5}. 
It involves a new technical aspect.  

}}


\section{Nondegeneracy of the statistical random field}\label{230806-5}

\subsection{{\colorb{Preliminary estimates}}}
Let $J\in\bbN$. 
For ${\bf c}=(c_0,c_1,...,c_p){\colorb{\in\bbR^{p+1}}}$, 
set 
\beas 
p({\bf c},x)
&=&
c_0+c_1x+\cdots+c_px^p.
\eeas
For $\delta>0$, let 
$
\calc_\delta
=
\{{\bf c};\>
|c_0|+|c_1|+\cdots+|c_p{|\colorblue{\geq}}\delta\}$. 
{\colorblue{
For $\ep\in(0,1)$, let 
$\calu_\ep=\{{\bf u}=(u_j)_{j=0}^p;\>\ep\leq\inf_j{\chp |u_j|}\leq\sup_j{\chp |u_j|}\leq\ep^{-1}\}$. 
Let ${\bf c}*{\bf u}=(c_ju_j)_j$. 
}}

\begin{lemma}\label{210103-1}
For any distinct positive numbers 
{\colorblue{$\{\alpha_i \}_{i=0}^p$ 
and $\delta,\ep>0$,}} 
there exist numbers $L>0$ and $n_0\in\bbN$ such that 
\beas 
\inf_{{\bf c}\in\calc_\delta}
\max_{i=0,...,p}
{\colorblue{\inf_{{\bf u}\in\calu_\ep}}}
|p({\bf c}{\colorblue{*{\bf u}}},n^{-\alpha_i})|\geq n^{-L}
\sskip\mbox{for all integers } n\geq n_0. 
\eeas
\end{lemma}
{\colorblue{
\proof 
Without loss of generality, we may assume $\alpha_0>\alpha_1>\cdots>\alpha_p$. 
Let ${\bf u}_i=(u_{ij})_{j=0}^p$ for $j=0,1,...,p$. 
Define a $(p+1)\times(p+1)$-matrix $A_n=[a^n_{ij}]_{i,j=0}^p$ by 
\beas 
a^n_{ij}
&=& 
u_{ij}n^{-\alpha_ij}. 
\eeas
Let ${\bf w}=A_n{\bf c}$, and write 
${\bf w}={^{\sf t}}(w_0,...,w_p)$ and 
${\bf c}={^{\sf t}}(c_0,...,c_p)$. 
For $i=0,...,p$, the $w_i$ is a function of $n$, ${\bf c}$ and ${\bf u}_i$: 
$w_i=w_i(n,{\bf c},{\bf u}_i)$. 
Suppose that ${\bf u}_i\in\calu_\ep$. 
By the fact that 
$\det A_n\sim \prod_{i=0}^p u_{ii}n^{-\alpha_ii}$ as $n\to\infty$, 
there exists a number $n_0\in\bbN$ such that 
${\chp |\det A_n|\geq 2^{-1}|\prod_{i=0}^p u_{ii}n^{-\alpha_ii}|}$ for all $n\geq n_0$ and all ${\bf u}_i\in\calu_\ep$ 
($i=0,...,p$). 
Since ${\bf c}=A_n^{-1}{\bf w}$, there exists a constant $K$ depending only on 
$p$ and $\ep$ and it holds that 
\beas 
|c_i| &\leq& Kn^{L'} \sum_{j=0}^p |w_j(n,{\bf c},{\bf u}_j)|
\sskip\mbox{\rm for all}\ n\geq n_0,\>
{\bf c}\in\bbR^{p+1},\>{\bf u}_j\in\calu_\ep\sskip (j=0,...,p),
\eeas
where $L'=\sum_{i=0}^p\alpha_ii$. 
For each ${\bf c}$,  the function $\calu_\ep\ni{\bf u}_j\mapsto |w_j|=|w_j(n,{\bf c},{\bf u}_j)|$ 
is obviously continuous, 
and hence 
$
\inf_{{\bf u}\in\calu_\ep} |w_j(n,{\bf c},{\bf u})|
=
|w_j(n,{\bf c},{\bf u}_j^*)|
$ 
for some ${\bf u}_j^*={\bf u}_j^*(n,{\bf c})\in\calu_\ep$. 
Applying the above inequality to ${\bf c}\in\calc_\delta$ , ${\bf u}_j^*$ ($j=0,...,p$) and $n\geq n_0$, we have  
\beas 
\delta &\leq& \sum_{i=0}^p |c_i|
\leq 
Kn^{L'} \sum_{i=0}^p \sum_{j=0}^p |w_j(n,{\bf c},{\bf u}_j^*)|
=
(p+1)Kn^{L'}  \sum_{j=0}^p\inf_{{\bf u}\in\calu_\ep} |w_j(n,{\bf c},{\bf u})|.
\eeas
The relation 
$p({\bf c}*{\bf u},n^{-\alpha_i})=w_i(n,{\bf c},{\bf u})$ 
completes the proof. 
\qed \y
}}
\begin{en-text}
\proof 
Write $\alpha=(\alpha_0,...,\alpha_p)$. 
Let $n\geq 2$. 
Define a $(p+1)\times(p+1)$-matrix $A_n=[a^n_{ij}]_{i,j=0}^p$ by 
\beas 
a^n_{ij}
&=& 
n^{-\alpha_ij}. 
\eeas
Let ${\bf w}=A_n{\bf z}$, and write 
${\bf w}={^{\sf t}}(w_0,...,w_p)$ and 
${\bf z}={^{\sf t}}(z_0,...,z_p)$. 
By the representation of 
Vandermonde's determinant, there exist 
positive numbers $L_0=L_0(p,\alpha)$ and 
$C=C(p,\alpha)$ independent of $n$ 
such that 
$|w_0|,...,|w_p|\leq n^{-L_0}$ implies 
$|z_0|,...,|z_p|\leq Cn^{-1}$. 
In other words, 
\bea\label{210102-3} 
\max_{j=0,...,p} |z_j|>Cn^{-1} 
&\Iku& 
\max_{j=0,...,p} |w_j|>n^{-L_0}. 
\eea
Therefore, if ${\bf c}\in\calc_\delta$ and 
$n\geq  n_0:=C(p+1)\delta^{-1}$, 
then 
\beas 
\max_{j=0,...,p}|p({\bf c},n^{-\alpha_j})|
&>& 
n^{-L_0}. 
\eeas
This completes the proof. \qed \y
\end{en-text}
%
Set
{\colorblue $\mathbb{S}=\{ \eta \in\bbR^{\sf d},\> | \eta|=1\}$. }
\begin{en-text}
$S^{{\sf d}-1}=\{ \eta \in\bbR^{\sf d},\> | \eta|=1\}$
and
$\bbS=S^{{\sf d}-1}\cap\bbL$. 
\end{en-text}
%
{\cod Let $\Theta$ be a subset of $\bbR^{{\colorred{{\sf p}}}}$.}
\footnote{{\cred This section is presented independently of the previous sections. 
In particular, if we assume the compactness of $\Theta$, it corresponds to $\bar{\Theta}$ of other sections. }}
Let $J\in\bbN$. 
%
{\cfg Let $f(x,\theta)$ is a function defined on a neighborhood of $\calx_0\times\Theta$. }

{\cod{Let 
$\bbL=T_x\bbR^{\sf d}\simeq\bbR^{\sf d}$.}}
A cone 
in the direction of $\xi\in\bbL\setminus\{0\}$ is 
defined by $C(\xi,a)=\{\eta\in\bbL;\>
\eta\cdot\xi\geq (1-a) |\eta||\xi|\}$ for $a\in(0,1)$. 
Let {\colorr{
$\tilde{\bbS}(\xi,\ep) =\bbS(\xi,\ep) \cup \bbS(-\xi,\ep)$
and $\tilde{S}(\xi,\ep) = S(\xi,\ep) \cup S(-\xi,\ep)$, 
where}}
$\bbS(\xi,\ep)=C(\xi,\ep)\cap\bbS$ and 
$S(\xi,\ep)=C(\xi,\ep)\cap\{ {\colorr{\eta}} \in \bbR^{{\sf d}}; | {\colorr{\eta}} |\leq\ep\}$. 
%
{\colorb{Let $\calx_\ell$ ($\ell=1,...,\bar{\ell}$) be subsets of {\cod$\calx_0$}}
{\cod and 
$\calx_0=\cup_\ell\calx_\ell$.}
{\chp Let $\tilde{S}_{x_0}(\xi,\ep)=x_0+\tilde{S}(\xi,\ep)$. }

\begin{description}
\item[[$N_0^\flat$\!\!]] 
{\colorb{There exist}} 
$\xi_{\ell,k}\in\bbS$, $\ep_{\ell,k}>0$,  and  
$\Theta_{\ell,k}\subset\Theta$ 
for $\ell=1,...,\bar{\ell}$ and $k=1,...,\bar{k}_\ell$ with $\Theta=\cup_k\Theta_{\ell,k}$ 
for each $\ell$, 
and 
{\colorb{bounded}} functions $b_{j,\ell,k}:{\chp \calx_\ell\times\Theta_{\ell, k}}\to\bbR$ 
for $j=0,...,J-1$ 
and {\colorb{a bounded function}} 
$b_{J,\ell,k}:\calx_\ell\times{\chp \tilde{S}_{x_0}(\xi_{\ell,k},\ep_{\ell,k})\times\Theta_{\ell, k}\times\bbS}\to\bbR$ 
such that the following conditions are fulfilled. 
\begin{itemize}
\item[(i)] 
For each $\ell=1,...,\bar{\ell}$ and $k=1,...,\bar{k}_\ell$, 
\beas 
{\cfg |f(x,\theta)|}
&\geq& 
|G_{\ell,k}(x_0,x,\theta,{\chp x-x_0})|
\eeas
for all $x_0\in\calx_\ell$, $x\in{\chp \tilde{S}_{x_0}(\xi_{\ell,k},\ep_{\ell,k})}$ 
and $\theta\in\Theta_{\ell,k}$, 
where
\beas 
G_{\ell,k}(x_0,x,\theta,\xi)
&=&
\sum_{j=0}^{J-1} b_{j,\ell,k}{\chp (x_0,\theta)}(\xi_{\ell,k}\cdot\xi)^j
+b_{J,\ell,k}{\chp (x_0,x,\theta,|\xi|^{-1}\xi)}
(\xi_{\ell,k}\cdot\xi)^J.
\eeas
\begin{en-text}
\beas 
|F(x_0,x,\theta,x-x_0)| 
&\geq& 
\bigg|
\sum_{j=0}^{J-1} b_j(x_0,\theta)\{\xi_k\cdot(x-x_0)\}^j
\\&&
+b_J(x_0,x,\theta)\{\xi_k\cdot(x-x_0)\}^J
\bigg|
\eeas
for all $x_0,x\in\hat{\calx}$ satisfying $x-x_0\in C_{\xi_k}$ 
and $\theta\in\Theta_k$. 
\end{en-text}
\item[(ii)] 
For each $\ell=1,...,\bar{\ell}$ and $k=1,...,\bar{k}_\ell$, 
\beas 
\inf_{(x_0,\theta)\in\calx_\ell\times \Theta_{\ell,k}}
{\chp \max_{j=0,...,J-1}
|b_{j,\ell,k}(x_0,\theta)|}>0. 
\eeas
\end{itemize}
\end{description}

\begin{en-text}
We consider {\chp symmetric} tensor fields 
{\colorb{$c_j:\calx_0\times\Theta\to\bbL^{\otimes j}$}} 
for $j=0,...,J-1$, 
and {\colorb{$c_J:\calx_0\times\hat{\calx}\times\Theta\to\bbL^{\otimes J}$.}} 
%
{\chp Suppose that $F$ admits the following expansion: }
\bea\label{210104-5}
F(x_0,x,\theta,\xi) 
&=& 
\sum_{j=0}^{J-1} c_j(x_0,\theta)[\xi^{\otimes j}]
+c_J(x_0,x,\theta)[\xi^{\otimes J}]
\eea
for 
{\colorb{$x_0 \in\calx_0$,}} 
$x\in\hat{\calx}$ {\cdo near $x_0$}, $\theta\in\Theta$ and $\xi\in\bbL$. 
{\colorb{Here the value of $c_j$ in ${\colorred{\bbL^{\otimes j}}}$ is identified with that of 
the dual space $({\colorred{\bbL^{\otimes j}}})'$. }}
{\chp Moreover {\cdo assume} the compatibility condition: 
\beas 
F(x_0,x,\theta,x-x_0) &=& F(x_0',x,\theta,x-x_0')
\eeas
for ${\cdo x_0'\in\calx_0}$ and $x\in\hat{\calx}${\cdo , both near $x_0\in\calx_0$. }
\end{en-text}

It {\colorb{will be shown}} that [N$_0^\flat$] follows 
from 
\begin{description}
\item[[$N_1^\flat$\!\!]] $\>${\bf (i)}  
{\chp $\calx_0$ and $\Theta$ are compact.} \vspace{-2mm}
\begin{description}
\item[(ii)] 
{\cfg $f$ admits a $C^J$-supporting function for each $(x_0,\theta)\in\calx_0\times\Theta$.}

\item[(iii)] 
For {\cfg the $C^J$-supporting function in (ii) for each $(x_0,\theta)\in\calx_0\times\Theta$, 
$\max_{j=0,...,J-1}
\big|c_j(x_0,\theta) \big|
>0$.} 
\end{description}
\end{description}
%
%

\begin{lemma}\label{210112-1}
$[N_1^\flat]$ implies $[N_0^\flat]$ {\cod for some $\{\calx_\ell\}_\ell$}. 
\end{lemma}
\proof 
\begin{en-text}
Let $\bbA=\{\xi\in\bbL; {\colorr{|\xi| \geq 1}} \}$. 
If 
$c_j(x_0,\theta)[\xi^{\otimes j}]=0$ for all $\xi\in\bbS$ hence 
for all $\xi\in\bbA$, then 
$c_j(x_0,\theta)=0$ as {\colorblue a symmetric} tensor. 
By induction, if necessary, express the elements of $\bbL$ 
by smaller number of coordinates, 
together with the property of the one-dimensional 
holomorphic function; or analytic function of multi-dimensional 
variable directly. 
\end{en-text}

For each $(x_0,\theta)\in\calx_0\times\Theta$, 
there exist $j(x_0,\theta)\in\{0,...,J-1\}$ and $\xi(x_0,\theta)\in\bbS$ 
such that 
\beas 
{\cod c_{j(x_0,\theta)}(x_0,\theta)
\not=0. }
\eeas 
\begin{en-text}
Let 
\beas 
\bar{F}(x_0,\theta,\xi)
&=&
\sum_{j=0}^{J-1} c_j(x_0,\theta)[\xi^{\otimes j}].
\eeas 
\end{en-text}
\begin{en-text}
Then there exists $\ep'(x_0,\theta)>0$ such that 
\beas 
\inf_{\xi\in \bbS(\xi(x_0,\theta),\ep'(x_0,\theta))}
|c_{j(x_0,\theta)}(x_0,\theta)
[\xi^{\otimes j(x_0,\theta)}]|
>0
\eeas
for all $x\in\hat{\calx}$ and all $\xi\in S(\xi(x_0,\theta),\ep'(x_0,\theta))$. %
Denote by $\Pi_{\xi(x_0,\theta)}:\bbL\to\bbL$ the orthogonal projection 
onto $\bbR\xi(x_0,\theta)$. 
Then it is easy to see that 
\beas 
|F(x_0,x,\theta,\xi)|
&\geq& 
\frac{1}{4} |c_{j(x_0,\theta)}(x_0,\theta)
[\Pi_{\xi(x_0,\theta)}\xi^{\otimes j(x_0,\theta)}]|
\eeas
for all $x\in\hat{\calx}$ and all $\xi\in S(\xi(x_0,\theta),\ep(x_0,\theta))$ 
for some $\ep(x_0,\theta)>0$. 
Indeed, the functional on the right-hand side is homogeneous in $\xi$, so 
we can compare the two terms appearing on the right-hand sides 
by assuming $\xi\in\bbS$. 
\beas 
|F(x_0',x,\theta',\xi)|
&=&
\big| \sum_{j=0}^{J-1}c_j(x_0',\theta')[\xi^{\otimes j}]
+c_J(x_0',x,\theta')[\xi^{\otimes J}]\big|
\\&\geq&
!!!!!!
\eeas
\end{en-text}
{\cod By continuity, there exist an open neighborhood $V(x_0,\theta)$  
of $(x_0,\theta)$} 
such that 
{\cfg 
\beas &&
\inf\bigg\{
\big|c_{j(x_0,\theta)}(x_0',\theta')
\big| 
;\>
(x_0',\theta')\in V(x_0,\theta)
\bigg\}
>0 
\eeas 
}
%
\begin{en-text}
and 
\bea\label{210224-1} 
|c_{j(x_0,\theta)}(x_0',\theta')[\xi(x_0,\theta)^{\otimes j(x_0,\theta)}]
|>\half |c_{j(x_0,\theta)}(x_0,\theta)
[\xi(x_0,\theta)^{\otimes j(x_0,\theta)}]|(>0)
\eea 
for all $x\in\hat{\calx}$, 
$(x_0',\theta')\in V(x_0,\theta)$ and 
$\xi\in S(\xi(x_0,\theta),\ep(x_0,\theta))$. 
\end{en-text}
\begin{en-text}
Let 
\beas 
M(x_0,\theta,\xi)=\max_{j=0,...,J-1}
|c_j(x_0,\theta)[\xi^{\otimes j}]|. 
\eeas
For $(\xi,m)\in\bbA\times\bbN$, let 
$V(\xi,m)=\{(x_0,\theta)
\in\calx_0\times\Theta;M(x_0,\theta,\xi)>m^{-1}\}$. 
Since the mapping $(x_0,\theta)\mapsto M(x_0,\theta,\xi)$ is continuous, 
$V(\xi,m)$ is open. 
By Assumption [N$_1^\flat$], 
$\Theta=
\big\{(x_0,\theta)\in\calx_0\times\Theta;\>
M(x_0,\theta,\xi)>0\mbox{ for some }\xi\in\bbA
\big\}$. 
\end{en-text}
The family $\{V(x_0,\theta)\}_{(x_0,\theta)\in\calx_0\times\Theta}$ is 
an open covering of the compact set $\calx_0\times\Theta$, 
consequently there are $(x_p,\theta_p)\in\calx_0\times\Theta$ 
$(p=1,...,\bar{p})$ 
such that $\cup_{p=1}^{\bar{p}}V(x_p,\theta_p)=\calx_0\times\Theta$. 
For each $n\in\bbN$, consider a family $\calu=\{U(n,m)\}_{m\in\bbN}$ of 
sets each of which is of the form 
\beas 
\bigg\{\big(2^{-n}i_1,2^{-n}(i_1+1)\big)\times\cdots\times
\big(2^{-n}i_{{\sf d}+{\colorred{{\sf p}}}},2^{-n}(i_{{\sf d}+{\colorred{{\sf p}}}}+1)\big)\bigg\}
\bigcap
(\calx_0\times\Theta)
\eeas 
for some $i_1,...,i_{{\sf d}+{\colorred{{\sf p}}}}\in\bbZ$. 
For each $(x_0,\theta)\in\calx_0\times\Theta$, 
there exists $U(n(x_0,\theta),m(x_0,\theta))\in\calu$ 
such that the closure $\overline{U(n(x_0,\theta),m(x_0,\theta))}$ is 
included in some $V(x_p,\theta_p)$. 
The family 
$\{U(n(x_0,\theta),m(x_0,\theta))\}_{(x_0,\theta)\in\calx_0\times\Theta}$ 
forms an open covering of the compact set $\calx_0\times\Theta$, 
therefore $\calu$ is already covered by a finite family 
$\{U(n_i,m_i)\}_{i=1,...,\bar{i}}$ 
of $\calu$. 
Set $\bar{n}=\max_{i=1,...,\bar{i}}\> n_i$ and 
divide each $\overline{U(n_i,m_i)}$ 
into some of elements of $\{\overline{U(\bar{n},m)}\}_{m\in\bbN}$. 
Thus we obtained a partition 
$\{\calx_\ell\times\Theta_k\}_{\ell=1,...,\bar{\ell},\>
k=1,...,\bar{k}}$ 
of $\calx_0\times\Theta$ such that 
$\calx_0=\cup_{\ell=1}^{\bar{\ell}}\calx_\ell$, 
$\Theta=\cup_{k=1,...,\bar{k}}\Theta_k$, and 
each $\calx_\ell\times\Theta_k$ is included 
in some $V(x_{p(\ell,k)},\theta_{p(\ell,k)})$. 
{\chp Set $\xi_{\ell,k}=\xi(x_{p(\ell,k)},\theta_{p(\ell,k)})$.} 

{\cod Based on the above construction of the covering, 
the function $G_{\ell,k}(x_0,x,\theta',\xi)$ is defined in a neighborhood of  $(x_0,\theta)$ 
through the expansion of the support function $g(P_{\xi_0}x,\theta')$ in $x$ 
around $x_0'$ of $(x_0',\theta')$ 
near $(x_0,\theta)$ as follows. 
Recall $c_j(x,\theta')=(j!)^{-1}\partial_x^j\{g(P_{\xi_0}x,\theta')\}[\xi_0^{\otimes j}]$. 
[$c_j(x,\theta')$ depends on $(x_0,\theta)$. ] 
Let 
$b_{j,\ell,k}(x_0',\theta')=c_j(x_0',\theta')$ for $j=0,...,J-1$ and 
\beas 
b_{J,\ell,k}(x_0',x,\theta',|\xi|^{-1}\xi)
&=&
J\int_0^1(1-s)^{J-1}c_J(x_0'+s(x-x_0'),\theta')ds.
\eeas
[In this case, $b_J$ does not depend on $\xi$. ] 
Then obviously $G_{\ell,k}(x_0',x,\theta',\xi)$ defined by the expansion of  [N$_0^\flat$] (i) satisfies 
$G_{\ell,k}(x_0',x,\theta',x-x_0')=g(P_{\xi_0}x,\theta')$. 
We choose $\ep(x_0,\theta)$ sufficiently small so that the inequality of  [N$_0^\flat$] (i) is valid. }
\qed \ \\

\begin{en-text}
as the following argument would be valid. 
Fix $(x_0,\theta)$. 

We may assume that $j_0=j(x_0,\theta)$ is the minimum $j\in\{0,1,...,J-1\}$ for which 
$c_j(x_0,\theta)\not\equiv0$. 
Let $\xi_0:=\xi(x_0,\theta)$. 
Let $P_0$ denote the orthogonal projection from $\bbL$ to the line $\bbR \xi_0$, and 
$P_0^\perp$ is the orthogonal projection to the orthogonal space of $\bbR \xi_0$. 
For $(x_0,\theta)$, define function $g_0$ of $(x',\theta')$ by 
\beas 
g_0(x',\theta') &=& \frac{1}{4} c_{j_0}(x_0,\theta')\big[\big(P_0(x'-x_0)\big)^{\otimes j_0}\big]
\sskip((x',\theta')\in V(x_0,\theta)). 
\eeas
If we take small $V(x_0,\theta)$, then 
\beas 
|F(x_0,x',\theta',x'-x_0)|
&\geq&
|g_0(x',\theta') |
\sskip((x',\theta')\in V(x_0,\theta)) 
\eeas
and that $\inf_{(x',\theta')\in V(x_0,\theta)}|g_0(x',\theta')|>0$. 
[Take $V(x_0,\theta)$ on which $|P_0^\perp(x'-x_0)|$ is small. 
Consider first $|2^{-1}c_{j_0}(x_0,\theta)[(P_0(x'-x_0)^{\otimes j_0}]|$ and make it less than 
$|F(x_0,x',\theta,x'-x_0)|$. Next use continuity of 
$c_{j_0}(x_0,\theta')$ and $F(x_0,x',\theta',x'-x_0)$ in $\theta'$. ] 
Thus $g_0$ is a supporting function for $F$ around $(x_0,\theta)$. 

Condition [N$_0^\flat$] requires a function $G(x_0',x,\theta',\xi)$ 
having an expansion around $x_0'$ 
for each $(x_0',\theta')\in V(x_0,\theta)$. 
Expand $g_0({\cdo x'},\theta')$ in $x'$ around $x_0'$ by using $x'-x_0=(x'-x_0')+(x_0'-x_0)$, 
to obtain 
\beas 
g_0(x',\theta')&=& \sum_{j'=0}^{j_0}b_{j,P_0}(x'_0,\theta')[(P_0(x'-x_0'))^{\otimes j'}] . 
\eeas
Thus to make $G_{\ell,k}$, we can use the function  
\beas 
G(x_0',x',\theta',\xi)&=& \sum_{j'=0}^{j_0}b_{j,P_0}(x'_0,\theta')[(P_0\xi)^{\otimes j'}] .
\eeas
[In this construction,$G(x_0',x',\theta',\xi)$ is independent of $x'$. ]
By the compatibility, 
$|F(x_0',x',\theta',x'-x_0')|=|F(x_0,x',\theta',x'-x_0)|\geq |g_0(x',\theta')|
=|G(x_0',x',\theta',x'-x_0')|$, which gives the inequality in [N$_0^\flat$] (i). 
The inequality of [N$_0^\flat$] (ii) is met if we make $V(x_0,\theta)$ sufficiently small. 
\qed\ \\
\end{en-text}
\begin{en-text}
Now we set 
\beas 
b_{j,\ell,k}(x_0,\theta,\xi)
&=&
{\colorblue c_{j}}(x_0,\theta)\big[
\big\{\big(\xi\cdot \xi(x_{p(\ell,k)},\theta_{p(\ell,k)})\big)^{-1}\xi\big\}
^{\otimes {\colorblue j}}\big] 
\eeas
for $(x_0,\theta,\xi)\in\calx_{p(\ell,k)}\times\Theta_{p(\ell,k)}
\times {\colorr{\tilde{\bbS} \big(\xi(x_{p(\ell,k)},\theta_{p(\ell,k)}),
\ep(x_{p(\ell,k)},\theta_{p(\ell,k)})\big)}}$, 
$\ell=1,...,\bar{\ell}$, $k=1,...,\bar{k}$ 
and $j=0,...,J-1$, and 
set 
\beas 
b_{J,\ell,k}(x_0,x,\theta,\xi)
&=&
c_J(x_0,x,\theta)\big[\big\{\big(\xi\cdot\xi(x_{p(\ell,k)},\theta_{p(\ell,k)})
\big)^{-1}\xi\big\}^{\otimes J}\big]
\eeas
for $(x_0,x,\theta,\xi)
\in\calx_{p(\ell,k)}\times\hat{\calx}\times\Theta_{p(\ell,k)}
\times {\colorr{\tilde{\bbS}(\xi(x_{p(\ell,k)},\theta_{p(\ell,k)}),
\ep(x_{p(\ell,k)},\theta_{p(\ell,k)}))}}$, 
$\ell=1,...,\bar{\ell}$, $k=1,...,\bar{k}$. 
{\colorb{Then $[N_0^\flat]$ holds for the so defined random fields $G_{\ell,k}$,  
each of which coincides with $F$ on each domain. }}
\qed \y
\end{en-text}

\begin{en-text}
{\bf koko} 
We know 
\beas 
\inf_{
(x_0,\theta)
\in\calx_\ell\times\Theta_k}
M(x_0,\theta,\xi_{p(\ell,k)})>m_{p(\ell,k)}^{-1}. 
\eeas
We write $\xi_\ell=\xi_{p(\ell,k)}$ and we may assume 
that $\xi_\ell\in\bbS$ by replacing $m_{p(\ell,k)}$ by 
some positive constant. 
By continuity, there exist $a_\ell>0$ 
such that 
\beas 
\inf_{
(x_0,\theta,\xi)
\in\calx_\ell\times\Theta_k\times (C_\ell\cap\bbA)}
M(x_0,\theta,\xi)>a_\ell
\eeas
for the cone $C_\ell=\{\eta\in\bbL;\>
\eta\cdot\xi_\ell\geq (1-a_\ell) |\eta|\}$. 
Therefore, there exists $j(\ell,k)\in\{0,1,...,J-1\}$ such that 
\beas 
\inf_{
(x_0,\theta,\xi)
\in\calx_\ell\times\Theta_k\times (C_\ell\cap\bbA)}
|c_{j(\ell,k)}(x_0,\theta)[\xi^{\otimes j}]|
>a_\ell,
\eeas

{\bf 連続性によって,
$\inf$を取る集合上絶対値の中身が符号を変えることはない．}

and hence 
\beas 
\inf_{
(x_0,\theta)
\in\calx_\ell\times\Theta_k}
|c_{j(\ell,k)}(x_0,\theta)[\xi^{\otimes j(\ell,k)}]|
\geq a_\ell|\xi|^{j(\ell,k)} 
\geq a_\ell(\xi\cdot\xi_\ell)^{j(\ell,k)}
\eeas
for all $\xi\in C_\ell$. 

Thus if we take 
\beas 
b_{j,\ell,k}(x_0,\theta) 
&=& 
\bigg\{
\begin{array}{ll}
a_\ell& \mbox{ if } j=j(x_0,\theta), \y
0& \mbox{otherwise}
\end{array}
\eeas
for 
$(x_0,\theta)\in\calx_\ell\times\Theta_k$,

{\bf 各cone上各項が定符号になるように$\xi_\ell$は
すでに選ばれている.
$F$の各項を符号に応じて上下から評価する．}

\vspace{1cm}
{\bf koko}

初期値$x_0$も摂動しておくか．

\bea\label{210104-5}
F(x_0,x,\theta,\xi) 
&=& 
\sum_{j=0}^{J-1} c_j(x_0,\theta)[\xi^{\otimes j}]
+c_J(x_0,x,\theta)[\xi^{\otimes J}]
\eea

there exists $\xi_\theta\in\bbL$ 
such that $\max_{j\in0,...,J-1}|c_j(x^*,\theta)[\xi^{\otimes j}]|>0$. 
In fact, if this value were $0$, then 

Divide $\Theta$ into very small subsets $\Theta_k$. 
Take $\theta_k\in\Theta_k$. 
If $|c_j(x^*,\theta_k)[\xi^{\otimes j}]|=0$ 
for all $\xi\in\bbL\cap\{\half\leq |\xi|\leq 1\}$, 
then $c_j(x^*,\theta_k)=0$ as the tensor, 
by induction together with the property of the one-dimensional 
holomorphic function. This contradicts [N$_1^\flat$]. 
Thus we can find a $\xi_k\in\bbS$ for which 
$|F(x^*,x^*,\theta_k,\xi_k)|>0$. 
Then we can make a $G$-function by replacing $\xi$ 
by $(\xi\cdot\xi_k)\xi_k$ with small cone $C_{\xi_k}$, 
and also 
$(x^*,x^*)$ by $(x_0,x)$ in a small set $\hat{\calx}^2$. 
Each $b_{j,k}$ $(j=0,...,J-1,J)$ is roughly to be set as 
either $(1+\ep')$ or $(1-\ep')$ times of ``$c_j$'' 
according to the sign of its value at $(x^*,x^*,\theta_k,\xi_k)$. 
\y
\end{en-text}

{\colorr{
Let
$E(\xi_0)=\{\xi\in\bbL;\>\xi\cdot\xi_0=1\}$ 
for $\xi_0\in\bbS$.
Set
$D(\xi_0,\ep) = E(\xi_0) \cap C(\xi_0,\ep)$
and   
$\tilde{D}(\xi_0,\ep) = D(\xi_0,\ep) \cup D(-\xi_0,\ep)$.
}}

\begin{lemma}\label{210102-2}
Suppose that 
$[N_0^\flat]$ is fulfilled. 
Then for 
any distinct positive numbers $\{\alpha_j\}_{j=0}^J$, 
there exist $L>0$ and $n_0\in\bbN$ such that 
\beas 
\min_{\ell=1,...,\bar{\ell},\>k=1,...,\bar{k}_\ell}
{\cdo \inf_{(x_0,\theta)\in\calx_\ell\times \Theta_{\ell,k}}}
\>\max_{i=0,...,J}
{\cfg \>\inf_{
\xi\in \tilde{D}(\xi_{\ell,k},\ep_{\ell,k})}
|f(x_0+n^{-\alpha_i}\xi,\theta)|
}
&\geq& 
n^{-L}
\eeas
for all $n\geq n_0$. 
{\cred This estimate is also valid for $|G_{\ell,k}(x_0,\cdot,\cdot,\cdot-x_0)|$ in place of $|f|$. 
}
\end{lemma}
\begin{en-text}
\begin{remark}\rm （書き直す！）
If $\hat{\calx}$ and $\Theta$ are compact and 
$c_j$ $(j=0,...,J-1)$ are continuous, then 
Condition (\ref{210102-1}) is equivalent to that 
\bea\label{210102-3}
\>\max_{j=0,...,J-1}\chi_j(\theta,\xi)>0 
\eea
for every $(\theta,\xi)\in\Theta\times\bbS$. 
\end{remark}
\end{en-text}
%
\proof 
It follows from [N$_0^\flat$](ii) that 
for some $\delta>0$, 
\beas 
\min_{(\ell,k)\in\{1,...,\bar{\ell}\}\times\{1,...,\bar{k}\}}
\inf_{(x_0,\theta)\in\calx_\ell\times{\colorb{\Theta_{\ell,k}}}}
{\chp 
\>\max_{j=0,...,J-1}
\big|b_{j,\ell,k}(x_0,\theta)
}
\big|
&>&
\delta. 
\eeas
{\chp 
For $(x_0,\theta)\in\calx_\ell\times\Theta_{\ell,k}$, 
for some $j_1=j_1(x_0,\theta)\in\{0,...,J-1\}$, $|b_{j_1,\ell,k}(x_0,\theta)|>\delta$. 
Then there exists $n_1\in\bbN$ such that 
\beas 
&&\sum_{j= j_1}^{J-1} b_{j,\ell,k}(x_0,\theta)n^{-\alpha_ij}
+b_{J,\ell,k}(x_0,x,\theta,|\xi|^{-1}\xi)n^{-\alpha_iJ}
\\&=&
b_{j_1,\ell,k}(x_0,\theta)n^{-\alpha_ij_1}\bigg(1+\ep\big(n,x_0,x,\theta,\xi,i,j_1\big)\bigg)
\eeas
for $n\geq n_1$, where 
$|\ep\big(n,x_0,x,\theta,\xi,i,j_1\big)|\leq1/2$ and $n_1$ depends only on 
$\delta$, $J$, $\alpha_i$ and $\|b_{j,\ell,k}\|_\infty$. 
Now we can apply Lemma \ref{210103-1} to $p=j_1$ with 
$u_{ij}=1$ for $i=0,...,j_1-1$ and $u_{ij_1}=1+\ep\big(n,x_0,x,\theta,\xi,i,j_1)$. 
}
\qed 
\begin{en-text}
{\colorb{\noindent 
Since functions $b_{j,\ell,k}$ are bounded and $\xi_{\ell,k}\cdot \xi$ is away from zero uniformly in 
$\xi\in\tilde{\bbS}(\xi_{\ell,k},\ep_{\ell,k})$, the definitive order for large $n$  between the sum and the last term by the nondegeneracy condition 
implies 
\beas 
|G_{\ell,k}(x_0,x,\theta,n^{{\colorblue -\alpha_i}}\xi)|
&\geq&
\min\bigg\{ 
\bigg|\sum_{j=0}^{J-1} b_{j,\ell,k}(x_0,\theta,|\xi|^{-1}\xi)(\xi_{\ell,k}\cdot\xi)^jn^{-j{\colorblue \alpha_i}}
\\&&
+\ep \>\| b_{J,\ell,k}\|_\infty
(\xi_{\ell,k}\cdot\xi)^Jn^{-J{\colorblue \alpha_i}}\bigg|;
\ep=\pm 1
\bigg\}
\eeas
for large $n$ and $(x_0,x,\theta)\in\calx_\ell\times\hat{\calx}\times \Theta_k$ and $\xi\in \tilde{D}(\xi_{\ell,k},\ep_{\ell,k})$. 
Apply Lemma \ref{210103-1} to $p=J$ with 
\beas 
c_j 
&=&
{\colorblue \>\inf_{\xi\in \tilde{D}(\xi_{\ell,k},\ep_{\ell,k})}}
b_{j,\ell,k}(x_0,\theta,|\xi|^{-1}\xi)
(\xi_{\ell,k}\cdot\xi)^j
\eeas
for $j=0,...,J-1$ and 
\beas 
c_J
&=&
{\colorblue \>\inf_{\xi\in \tilde{D}(\xi_{\ell,k},\ep_{\ell,k})}}
\ep\>
\|b_{J,\ell,k}\|_\infty
(\xi_{\ell,k}\cdot\xi)^J, 
\eeas 
$x\in\hat{\calx}$, 
and twice for $\ep=1$ and $\ep=-1$, 
to obtain the result. 
}}
\qed \\
\end{en-text}
}

\begin{en-text}
According to the assumption (\ref{210102-1}), 
there exists a number $n_0$ independent of 
$(\theta,\xi_0,...,\xi_{J-1})\in\Theta\times\bbS^J$
 as well as $x$, and 
\beas 
\max_{j=0,...,J-1} |\chi_j(\theta,\xi_j)|>Cn^{-1} 
\sskip(\forall (\theta,\xi_0,...,\xi_{J-1})\in\Theta\times\bbS^J)
\eeas
for all $n\geq n_0$. 
We apply the relation (\ref{210102-3}) to 
$z_j=\chi_j(\theta,\xi_j)$ $(j=0,...,J-1)$ and 
$z_J=\chi_J(x,\theta,\xi)$ 
to obtain 
\beas
\max_{k=0,...,J-1}|F(x,\theta,n^{-\alpha_k}\xi_j)|>n^{-L}
\eeas
for all $(x,\theta,\xi_0,...,\xi_J)\in\hat{\calx}\times\Theta\times\bbS^J$ 
and $n\geq n_0$. 
\qed\y
\end{en-text}

\subsection{{\colorb{Nondegeneracy of the index $\chi_0$}}}\label{230608-3}

{\cdo 
Suppose that $X=(X_t)_{t\in[0,T]}$ is a ${\sf d}$-dimensional separable process.
\footnote{This section gives a way to the estimate of the key index $\chi_0$. 
Since the method is general, we write it for a general stochastic process $X$, 
apart from the It\^o process $X$ in Section \ref{PLDP}. 
%
}
$\Theta$ is a set in $\bbR^{\sf p}$. 
$\calx_0$ and $\hat{\calx}$ are subsets of $\bbR^{\sf d}$ with $\calx_0\subset\hat{\calx}^o$, the interior of $\hat{\calx}$ in $\bbR^{\sf d}$. 
Let $\theta^*\in\Theta$ and let
\beas 
\chi_0 &=& \inf_{\theta\in\Theta:\theta\not=\theta^*}
\frac{\int_0^T\bbQ(X_t,\theta)dt}{|\theta-\theta^*|^2}
\eeas
for a function $\bbQ:\hat{\calx}\times\Theta\to\bbR$.
\footnote{We use the same symbol $\chi_0$ for the key index as Section \ref{PLDP} since we will apply 
the nondegeneracy results here to $\bbQ=Q(\cdot,\cdot,\theta^*)/2T$ in Section \ref{230806-3}}
}
\begin{en-text}
By (\ref{230807-1}), 
\beas 
-2\bbY(\theta)
&=&
\frac{1}{T} \int_0^T Q(X_t,\theta,\theta^*)dt. 
\eeas
\end{en-text}

{\colorb{\noindent 
Furthermore, 
suppose 
{\cfg{
\begin{description}
\item[[R\!\!]] 
There exist a function {\colorb{$f:\hat{\calx}\times\Theta\to\bbR$}} 
and a constant 
$\varrho\in(0,\infty)$ satisfying 
\beas 
{\cdo \bbQ(x,\theta)}|\theta-\theta^*|^{-2}
\geq |f(x,\theta)|^\varrho
\eeas
for all $(x,\theta)\in\hat{\calx}\times\Theta$, 
{\cred and the function $f(\cdot,\theta)$ is Lipschitz continuous on $\hat{\calx}$ uniformly in $\theta\in\Theta$}. 
%
\end{description}
}}
\noindent Here 
$f$ and $\varrho$ possibly depend on $\theta^*$. 
An example of $f$ is ${\cdo \bbQ(x,\theta)|\theta-\theta^*|^{-2}}$ itself 
for $\varrho=1$, however we have much more 
freedom of choice of $f$ and $\varrho$. 
{\colorb{Introducing the subfield $f$ facilitates application of the result. }}

\begin{en-text}
Define $F(x_0,x,\theta,\xi)$ by 
(\ref{210104-5}) 
\bea\label{210323-1}
F(x_0,x,\theta,\xi) 
&=& 
\sum_{j=0}^{J-1} c_j(x_0,\theta)[\xi^{\otimes j}]
+c_J(x_0,x,\theta)[\xi^{\otimes J}]
\eea
for 
{\colorb{$x_0 \in\calx_0$}}, 
$x\in\hat{\calx}$ near $x_0$, $\theta\in\Theta$ and $\xi\in\bbL$ 
with 
$
c_j(x_0,\theta) 
= 
(j!)^{-1}
\partial_x^jf(x_0,\theta)
$ 
and 
\beas 
c_J(x_0,x,\theta)
&=&
((J-1)!)^{-1}
\int_0^1 (1-s)^{J-1}
\partial_x^Jf(x_s,\theta)ds, 
\eeas
where $x_s=(1-s)x_0+sx$. 
Then 
$f(x,\theta)=F(x_0,x,\theta,x-x_0)$ 
{\cdo for $x\in\hat{\calx}$ near $x\in\calx_0$. 
}
\end{en-text}

%
\begin{en-text}
\textcolor{red}{
Furthermore, we assume that
for some $C>0$,
$$
\max_{j=0,\ldots, J-1}
\sup_{ (x_0, \theta) \in \hat{\calx} \times \Theta}
|c_j(x_0,\theta)|
+
\max_{|{\bf n}|=0,1}
\sup_{ {(x_0, x) \in \hat{\calx} \times \hat{\calx}}
\atop \theta \in \Theta}
|\partial_x^{\bf n} c_J(x_0,x,\theta)| <C.
$$
}
\end{en-text}

{\colorr{We denote by $\calp$ the set of sequences $({\mathfrak a}_n)_{n\in\bbN}$ of nonnegative numbers 
satisfying the condition that 
for every $L>0$, there exists a number $C_L$ such that 
${\mathfrak a}_n\leq C_L/n^L$ for all $n\in\bbN$. 
$\cale$ denotes the set of sequences $({\mathfrak a}_n)_{n\in\bbN}$ of nonnegative numbers such that 
for some $c>0$, ${\mathfrak a}_n\leq c^{-1}e^{-cn^c}$ for all $n\in\bbN$. 
}}
{\colorb{
Now we assume the following conditions for 
the nondegeneracy of the deterministic field and the variation of the underlying stochastic process. 
Condition [C] is for estimate of a modulus of continuity of $X$. 
\begin{description}
\item[[$N_0$\!\!]] 
There exist $T_0\in(0,T)$, 
subsets $\calx_\ell{\colorsb{\subset\hat{\calx}}}$ ($\ell=1,...,\bar{\ell}$) 
with $\calx_0{\colorsb{\subset}}\cup_\ell\calx_\ell$ and $\overline{\calx}_\ell\subset\hat{\calx}^o$
\begin{en-text}
, 
$\xi_{\ell,k}\in\bbS$, $\ep_{\ell,k}>0$,   
sets $\Theta_{\ell,k}\subset\Theta$ 
$(\ell=1,...,\bar{\ell}$ and $k=1,...,\bar{k}_\ell)$ with $\Theta=\cup_k\Theta_{\ell,k}$ 
for each $\ell$, 
\end{en-text}
{\cdo 
for which the following conditions hold: }
\begin{description}
\item[(i)] 
{\cdo 
[N$_0^\flat$] holds and $\cup_{\ell,k}\cup_{x_0\in\calx_\ell}\overline{B(x_0,\ep_{\ell,k})}\subset\hat{\calx}$.
\footnote{$B(x,\ep)$ is the open ball centered $x$ with radius $\ep$. }
}
\begin{en-text}
For each $\ell=1,...,\bar{\ell}$ and $k=1,...,\bar{k}_\ell$, 
\beas 
{\colorblue \inf_{(x_0,\theta)\in\calx_\ell\times \Theta_{\ell,k}}}
\max_{j=0,...,J-1}
{\colorblue \inf_{\xi \in \tilde{\bbS}(\xi_{\ell,k},\ep_{\ell,k})}}
|\partial_x^jf(x_0,\theta)[\xi^{\otimes j}]|>0. 
\eeas
\end{en-text}

\item[(ii)] 
For each $(\ell,k)$, 
there exist 
a positive constants $c_0$ 
and 
distinct positive numbers $\{\alpha_j:=\alpha_j(\ell,k)\}_0^J$ such that 
the sequence 
\bea
\bigg(1-
 P\bigg[ &&
\hspace{-0.8cm}
{\colorsb{\bigcup_{\ell=1,...,\bar{\ell}}}}
{\colorred{\bigcap_{k=1,...,\bar{k}_\ell}}}
{\colorsb{\bigcup_{s\in[0,T_0]}}}
\bigcap_{j=0,...,J} \bigcup_{t\in[0,n^{-c_0}]}
\bigg\{ X_s\in\calx_\ell,\>
X_{s+t}-X_s\in \tilde{S}(\xi_{\ell,k},\ep_{\ell,k}) 
\nn\\&&
\mbox{ and }
\big|\xi_{\ell,k}\cdot (X_{s+t}-X_s)\big|=n^{-\alpha_j} 
\bigg\}
{\colorred{\bigg]}}
\bigg)_{n\in\bbN}
\eea
is in $\calp$. 
\end{description}
\end{description}
}}
\begin{en-text}
\begin{description}
\item[[$N_0$\!\!]] 
Condition [$N_0^\flat$] is satisfied, and for each $(\ell,k)$, there exist 
positive constants 
$c_0,c_1, \textcolor{red}{c_2}$ and 
distinct positive numbers $\{\alpha_j\}_0^J$ such that 
\bea
\min_{j=0,...,J}\inf_{x_0\in\calx_\ell}
& P\bigg[ &
\bigg\{ X_t-x_0\in S(\xi_{\ell,k},\ep_{\ell,k}) \mbox{ and }
\xi_{\ell,k}\cdot (X_t-x_0)=n^{-\alpha_j} \bigg\} \mbox{ or } 
\nonumber 
\\
& &
\textcolor{red}{
\bigg\{ X_t-x_0\in S(-\xi_{\ell,k},\ep_{\ell,k}) \mbox{ and } 
-\xi_{\ell,k}\cdot (X_t-x_0)=n^{-\alpha_j} \bigg\}
}
\nonumber 
\\
& & 
\mbox{ for some }
t\in[0,n^{-c_0}]
\ \bigg| \  
X_0=x_0 \ \bigg]
\geq
1-c_1^{-1}e^{-c_1 n^{\textcolor{red}{c_2}}} \label{u-n0}
\eea
for all $n\in\bbN$. 
\end{description}
\end{en-text}

{\colorb{
\begin{description}
\item[[C\!\!]] There exist positive constants $\beta_0$ 
such that 
the sequence 
\beas 
{\colorsb{\bigg(
P\bigg[}}
\sup_{s,t\in[0,T]: 
X_s\in {\colorg{\hat{\calx}}},\>t\in[s,s+n^{-1}]}
|X_t-X_s|\geq n^{-\beta_0}
\bigg] 
\bigg)_{n\in\bbN}
\eeas
is in $\calp$. 
\end{description}
}}
\begin{en-text}
\begin{description}
\item[[C\!\!]] There exist positive constants $\beta_0$ (\textcolor{red}{$\beta_0<1$}), $\beta_1$ , $\textcolor{red}{d_0}$ and $d_1$ 
such that 
\beas 
\sup_{x_0 \in\calx_0}
P\bigg[
\sup_{t \in[0,n^{-1}] \atop 
\textcolor{red}{s \in [0,n^{-{d_0}/{l}}]} } 
|X_{t+s}-X_s|\geq n^{-\beta_0}\>\big| \> X_0=x_0 
\bigg] 
&\leq&
d_1^{-1}\exp(-d_1n^{\beta_1/{\textcolor{red}{l}}})
\eeas
for all $n\in\bbN$ \textcolor{red}{and all $l\geq1$}. 
\end{description}
\end{en-text}
%

{\colorb{
\begin{proposition}\label{230608-1}
Under $[N_0]$, $[R]$ and $[C]$,  
Condition $[H2]$ holds true.\footnote{{\cred Of course, ``$\chi_0$'' is the one in this section. }} 
\end{proposition}
}}

{\colorb{
\proof 
Let 
\beas 
\Omega_{n,\ell,k,j,s}
&=&
\big\{
X_s\in\calx_\ell\big\}
\bigcap\bigg\{
X_{s+t} -X_s \in \tilde{S}(\xi_{\ell,k}, \ep_{\ell,k})
\\&&
 \mbox{ and } 
\big|\xi_{\ell,k}\cdot (X_{s+t}-X_s)\big| = n^{-\alpha_j}
\mbox{ for some } t \in [0,n^{-c_0}] \bigg\}.
\eeas
In what follows, we consider sufficiently large $n$. 
For $\omega\in\bigcap_{j=0,...,J}\Omega_{n,\ell,k,j,s}$ and $\theta\in\Theta_{\ell,k}$, 
there are random times $\tau_j={\colorred{\tau_j(\omega,n,\ell,k,s)}}\in(s,s+n^{-c_0}]$ such that 
\beas 
X_{\tau_j}-X_s\in\tilde{S}(\xi_{\ell,k},\ep_{\ell,k})\sskip\mbox{and}\sskip
|\xi_{\ell,k}\cdot (X_{\tau_j}-X_s)|=n^{-\alpha_j}.
\eeas
Then Lemma \ref{210102-2} yields 
\beas 
{\colorblue \max_{j=0,...,J}} |f(X_{\tau_j},\theta)|
&\geq&
\min_{\ell=1,...,\bar{\ell},\>k=1,...,\bar{k}_\ell}
{\cdo \inf_{(x_0,\theta')\in\calx_\ell\times\Theta_{\ell,k}}}
{\colorblue \max_{j=0,...,J}}
\inf_{\cfg 
\xi\in\tilde{D}(\xi_{\ell,k},\ep_{\ell,k})} 
{\cfg |f(x_0+n^{-\alpha_j}\xi{\cod ,\theta'})|}
\geq n^{-L}
\eeas
for $n\geq n_0$, whre $L>0$ and $n_0$ are depending only on $\{\alpha_j\}_{j=0}^J$ and independent of 
$\omega\in\bigcap_{j=0,...,J}\Omega_{n,\ell,k,j,s}$ and $\theta\in\Theta_{\ell,k}$. 
%
%

Take $\kappa\in\bbN$ such that $\kappa>L\varrho [\beta_0(\varrho\wedge1)]^{-1}$, and let $L'>\kappa+L\varrho$. 
We have 
\beas 
P\big[{\cdo \chi_0}\leq n^{-L'}\big]
&=&
{\colorred{
P\bigg[\inf_{\theta\not=\theta^*}\int_0^T\frac{{\cdo \bbQ(X_t,\theta)}}{|\theta-\theta^*|^2}dt\leq  n^{-L'}
\bigg]}}
\\&\leq&
{\colorred{
\sum_{\ell=1}^{\bar{l}}
P\bigg[ \bigg\{\bigcup_{k=1}^{\bar{k}_\ell}\bigg(
\inf_{\theta\in\Theta_{\ell,k}}\int_0^T |f(X_t,\theta)|^\varrho dt\leq  n^{-L'}\bigg)\bigg\}
\bigcap\bigg(
{\colorsb{\bigcap_{k=1}^{{\colorg{\bar{k}_\ell}}}\bigcup_{s\in[0,T_0]}}}
\bigcap_{j=0,...,J}\Omega_{n,\ell,k,j,s}\bigg) \bigg]+{\mathfrak a}_n
}}
\\&\leq&
{\colorred{
\sum_{\ell=1}^{\bar{l}}
P\bigg[ \bigcup_{k=1}^{\bar{k}_\ell}\bigg\{\bigg(
\inf_{\theta\in\Theta_{\ell,k}}\int_0^T |f(X_t,\theta)|^\varrho dt\leq  n^{-L'}\bigg)
\bigcap\bigg(
\bigcup_{s\in[0,T_0]}\bigcap_{j=0,...,J}\Omega_{n,\ell,k,j,s}\bigg)\bigg\} \bigg]+{\mathfrak a}_n
}}
%
%
%
\\&\leq&
\sum_{\ell=1}^{\bar{\ell}}\sum_{k=1}^{\bar{k}_\ell}
{\colorsb{P}}\bigg[\bigg\{\inf_{\theta\in\Theta_{\ell,k}}\int_0^T |f(X_t,\theta)|^\varrho dt\leq  n^{-L'}\bigg\}
\bigcap 
\bigg(
{\colorred{ \bigcup_{s\in[0,T_0]} }}
\bigcap_{j=0,...,J}\Omega_{n,\ell,k,j,s}\bigg)
\bigg]
+{\mathfrak a}_n
\eeas
with $({\mathfrak a}_n)_{n\in\bbN}\in\calp$. 

Thanks to [C], we have 
\beas 
&&
P\bigg[\bigg\{\inf_{\theta\in\Theta_{\ell,k}}\int_0^T |f(X_t,\theta)|^\varrho dt\leq  n^{-L'}\bigg\}
\bigcap \bigg(\bigcup_{s\in[0,T_0]}\bigcap_{j=0,...,J}\Omega_{n,\ell,k,j,s}\bigg)
\bigg]
\\&\leq&
P\bigg[\inf_{\theta\in\Theta_{\ell,k}}{\colorblue \max_{j=0,...,J}}|f(X_{\tau_j},\theta)|^\varrho n^{-\kappa}
-\sup_{\theta\in\Theta_{\ell,k}}{\colorblue \max_{j=0,...,J}}
\int_{\tau_j}^{\tau_j+n^{-\kappa}}\big| |f(X_t,\theta)|^\varrho - |f(X_{\tau_j},\theta)|^\varrho\big| dt
\leq n^{-L'}
\bigg]
\\&\leq&
{\mathfrak b}_n
\eeas
{\colorsb{for}} some $({\mathfrak b}_n)_{n\in\bbN}\in\calp$. 
Indeed, 
$
n^{-L\varrho-\kappa}-A n^{-\kappa-\kappa\beta_0(\varrho\wedge1)}> n^{-L'}
$ 
as $n$ becomes large for every $A>0$,  
and 
\beas 
\sup_{\theta\in\Theta_{\ell,k}}\int_{\tau_j}^{\tau_j+n^{-\kappa}}
\big| |f(X_t,\theta)|^\varrho - |f(X_{\tau_j},\theta)|^\varrho\big| dt
&\leq&
n^{-\kappa-\kappa\beta_0(\varrho\wedge1)}
\eeas
on the event 
{\colorg{
$ \big\{ \sup_{s,t\in[0,T]: 
X_s \in \hat{\calx}',\>t\in[s,s+n^{-\kappa}]}
|X_t-X_s| < n^{-\kappa \beta_0}
\big\},
$
{\cdo where $\hat{\calx}'$ is a suitable neighborhood of $\overline{\cup_\ell\calx_\ell}$. }
}}
%
Therefore, we obtain 
$
P\big[\chi_0\leq  n^{-L'}\big]
= 
{\mathfrak c}_n
$ 
for some $({\mathfrak c}_n)_{n\in\bbN}\in\calp$. 
This gives [H2]. \qed

\begin{remark}\rm 
In Proposition \ref{230608-1}, 
if the sequences in $[N_0]$(ii) and [C] are in $\cale$, 
then we obtain 
$P[\chi_0\leq r^{-1}]\leq c^{-1}e^{-cr^c}$ ($r>0$) for some $c>0$. 
The same remark is also for Proposition \ref{2306-3}. 
\end{remark}

{\cred 
\begin{remark}\label{241205-1}\rm 
If we strengthen $[N_0]$ (ii) by replacing $\tilde{S}(\xi_{\ell,k},\ep_{\ell,k})$ by $S(\xi_{\ell,k},\ep_{\ell,k})$, 
then the inequality in Lemma \ref{210102-2} but with $D(\xi_{\ell,k},\ep_{\ell,k})$ for 
$\tilde{D}(\xi_{\ell,k},\ep_{\ell,k})$ 
is still sufficient to prove the same result as Proposition \ref{230608-1}. 
This formulation will work for nondegenerate diffusions. 
However the original one is worth stating because it is easy to give an example such that  
the process $X$ moves toward $\xi_{\ell,k}$ or $-\xi_{\ell,k}$ with probability $1/2$. 
\end{remark}
}

\begin{en-text}
We assume $[N_0]$ and $[C]$. 
Let 
\begin{eqnarray*}
\Omega_{n,\ell,k,j}^{(1)} 
&=&
\left\{ X_t -X_0 \in S(\xi_{\ell,k}, \ep_{\ell,k}) \mbox{ and } \xi_{\ell,k}\cdot(X_t-X_0) = n^{-\alpha_j}
\mbox{ for some } t \in [0,n^{-c_0}] \right\}, 
\\
\Omega_{n,\ell,k,j}^{(2)} 
&=&
\left\{ X_t -X_0 \in S(-\xi_{\ell,k}, \ep_{\ell,k}) \mbox{ and } -\xi_{\ell,k}\cdot(X_t-X_0) = n^{-\alpha_j}
\mbox{ for some } t \in [0,n^{-c_0}] \right\},
\\
\Omega_{n,\ell,k,j} &=& \Omega_{n,\ell,k,j}^{(1)}\cup \Omega_{n,\ell,k,j}^{(2)}, 
\\
\Omega_{n} &=& \bigcup_{\ell=1}^{\bar{\ell}} \bigcup_{k=1}^{\bar{k}_\ell} 
\left[ \left( \bigcap_{j=0}^J \Omega_{n,\ell,k,j} \right) \cap \left\{ X_0 \in \calx_\ell \right\} \right].
\end{eqnarray*}
Since $\calx_0 =\cup_\ell \calx_\ell$, $\calx_\ell \cap \calx_m =\phi$ ($\ell \ne m$),
\begin{eqnarray*}
P[ \Omega_n ]
&=& \sum_{\ell=1}^{\bar{\ell}}
P \left[ \left( \bigcup_{k=1}^{\bar{k}_\ell} \bigcap_{j=0}^J \Omega_{n,\ell,k,j} \right) \cap \left\{ X_0 \in \calx_\ell \right\} \right]
\\
&\geq & \sum_{\ell=1}^{\bar{\ell}} 
\left( \min_{k=1, \ldots, \bar{k}_\ell} \ \inf_{x_0 \in \calx_\ell} 
P \left[ \left. \bigcap_{j=0}^J \Omega_{n,\ell,k,j} \right|  X_0 =x_0 \ \right] \right) P[ X_0 \in \calx_\ell ] 
\\
&\geq & 
\min_{\ell=1, \ldots, \bar{\ell}} \ \min_{k=1, \ldots, \bar{k}_\ell} \ \inf_{x_0 \in \calx_\ell} 
P \left[ \left. \bigcap_{j=0}^J \Omega_{n,\ell,k,j} \right|  X_0 =x_0 \ \right].
\end{eqnarray*}
By $[N_0]$, for every $\ell$ and $k$, there exist positive constants $c_0, c_1, c_2$ and
distinct positive numbers $\{\alpha_j\}_{j=0,\ldots,J}$ such that
$$
\inf_{x_0 \in \calx_\ell} P \left[ \left. \bigcap_{j=0}^J \Omega_{n,\ell,k,j} \right|  X_0 =x_0 \ \right]  
\geq 1 - c_1^{-1} \exp (-c_1 n^{c_2})
$$
since
\begin{eqnarray*}
\inf_{x_0 \in \calx_\ell} P \left[ \left. \bigcap_{j=0}^J \Omega_{n,\ell,k,j} \right|  X_0 =x_0 \ \right] 
&\geq& 
\inf_{x_0 \in \calx_\ell} P \left[ \left. \bigcap_{j=0}^{J-1} \Omega_{n,\ell,k,j} \right|  X_0 =x_0 \ \right]
+ \inf_{x_0 \in \calx_\ell} P \left[ \left. \Omega_{n,\ell,k,J} \right|  X_0 =x_0 \ \right]
\\
& &
- \sup_{x_0 \in \calx_\ell} P \left[ \left. \left( \bigcap_{j=0}^{J-1} \Omega_{n,\ell,k,j} \right) \cup \Omega_{n,\ell,k,J} \right|  X_0 =x_0 \ \right]
\\
&\geq& 
\inf_{x_0 \in \calx_\ell} P \left[ \left. \bigcap_{j=0}^{J-1} \Omega_{n,\ell,k,j} \right|  X_0 =x_0 \ \right]
+1 - c_{1,J}^{-1} \exp (-c_{1,J} n^{c_{2,J}}) -1
\\
&\geq& 
\inf_{x_0 \in \calx_\ell} P \left[ \left. \Omega_{n,\ell,k,0} \right|  X_0 =x_0 \ \right]
- \sum_{j=1}^J c_{1,j}^{-1} \exp (-c_{1,j} n^{c_{2,j}}) 
\\
&\geq& 
1 - \sum_{j=0}^J c_{1,j}^{-1} \exp (-c_{1,j} n^{c_{2,j}}).   
\end{eqnarray*}
Thus, there exist positive constants $c_0, c_1$ and $c_2$ such that
$P[\Omega_n] \geq 1 - c_1^{-1} \exp (-c_1 n^{c_2}).$

Fix $\ell$ and $k$. 
We consider $J+1$ stopping times $\tau_j =\min \{ \tau_j^{(1)}, \tau_j^{(2)} \}$, where
\beas 
\tau_j^{(1)} &=& \inf\{t>0; 
X_t-X_0\in S(\xi_{\ell,k},\ep_{\ell,k}) 
\mbox{ and }
\xi_{\ell,k}\cdot (X_t-X_0)= n^{-\alpha_j}\},
\\
\tau_j^{(2)} &=& \inf\{t>0; 
X_t-X_0\in S(-\xi_{\ell,k},\ep_{\ell,k}) 
\mbox{ and }
-\xi_{\ell,k}\cdot (X_t-X_0)= n^{-\alpha_j}\}.
\eeas
Let
$$
A_n = \left\{ \inf_{\theta \in \Theta} \max_{j=0,\ldots, J} |F(X_0, X_{\tau_j}, \theta, X_{\tau_j}-X_0)| \geq n^{-L} \right\}
$$
for all $n \geq n_0$, where $L$ and $n_0$ are given Lemma \ref{210102-2}.
Let 
$$
B_n = \left\{ \sup_{t \in [0,n^{-k}] \atop s \in \left[ 0,n^{-k d_1/{l}} \right]} |X_{t+s}-X_s| \geq n^{-k \beta_0} \right\}
$$
for all $n \in \mathbb{N}$, $k>0$ and $l \geq 1$, where $d_1$ and $\beta_0$ are given in $[C]$.
By Lemma \ref{210102-2}, one has that for all $n \geq n_0$,
\beas
& & 1-c_1^{-1} \exp (-c_1 n^{c_2})
\\
&\leq& P[\Omega_n] = P[\Omega_n \cap A_n] \leq P[\Omega_n \cap A_n \cap B_n] + P[\Omega_n \cap A_n \cap B_n^c]
\\
&\leq& P[B_n] + P[\Omega_n \cap A_n \cap B_n^c].
\eeas
It follows from $[C]$ that $P[B_n] \leq \sup_{x_0 \in \calx_0} P[B_n |X_0=x_0] P[X_0 \in \calx_0] \leq d_1^{-1} \exp (-d_1 n^{k \beta_1/l})$
for some positive constants $d_1$ and $\beta_1$.

Next we estimate $P[\Omega_n \cap A_n \cap B_n^c]$. On $\Omega_n \cap A_n \cap B_n^c$,
there exists $j \in \{0,1,\ldots,J \}$ such that $\inf_{\theta \in \Theta} |F(X_0,X_{\tau_j},\theta,X_{\tau_j}-X_0)| \geq n^{-L}$.
Let $f(x,\theta)=F(X_0,x,\theta,x-X_0)$.
We then have that for any $k>0$,
\beas
\chi_0 &=& \inf_{\theta \ne \theta^*} \int_0^T \frac{Q(X_t,\theta)}{|\theta - \theta^*|^2}dt
\\
&\geq& \inf_{\theta \ne \theta^*} \int_{\tau_j}^{\tau_j+n^{-k}}  |F(X_0,X_t,\theta,X_t-X_0)|^\varrho dt
\\
&\geq& \inf_{\theta \ne \theta^*} \int_{\tau_j}^{\tau_j+n^{-k}}  |F(X_0,X_{\tau_j},\theta,X_{\tau_j}-X_0)|^\varrho dt
\\
& & - \sup_{\theta \ne \theta^*} \int_{\tau_j}^{\tau_j+n^{-k}} | |F(X_0,X_t,\theta,X_t-X_0)|^\varrho-|F(X_0,X_{\tau_j},\theta,X_{\tau_j}-X_0)|^\varrho | dt.
\\
&\geq& n^{-k} \times n^{-\varrho L} 
- \sup_{\theta \ne \theta^*} \int_{\tau_j}^{\tau_j+n^{-k}} | |f(X_t,\theta)|^\varrho-|f(X_{\tau_j},\theta)|^\varrho | dt.
\eeas
Let $Y_{t,s} =X_{\tau_j} +s(X_t-X_{\tau_j})$. It can be shown that on $\Omega_n \cap A_n \cap B_n^c$,
\begin{eqnarray}
& & \sup_{\theta \in \Theta} \ \sup_{t \in [\tau_j, \tau_j+n^{-k}] \atop s \in [0,1] } \ |f(Y_{t,s},\theta)| <C, \label{u-1}
\\
& & \sup_{\theta \in \Theta} \ \sup_{t \in [\tau_j, \tau_j+n^{-k}] \atop s \in [0,1] } \ |\partial_x f(Y_{t,s},\theta)| <C, \label{u-2}
\\
& &
\inf_{\theta \in \Theta} \ \inf_{t \in [\tau_j, \tau_j+n^{-k}] \atop s \in [0,1] } \ |f(Y_{t,s},\theta)| \geq \frac{1}{2}n^{-L} \label{u-3}
\end{eqnarray}
for sufficiently large $n$.
Indeed, noting that $\tau_j \leq n^{-c_0}$ on $\Omega_n \cap A_n \cap B_n^c$,
for $l$ satisfying that $l \geq k d_0/{c_0}$,
\beas
\sup_{t \in [\tau_j, \tau_j+n^{-k}] \atop s \in [0,1] } \ |Y_{t,s} - X_{\tau_j} | 
&\leq& \sup_{t \in [0, n^{-k}] } \ |X_{t+\tau_j} - X_{\tau_j} | 
\leq \sup_{t \in [0, n^{-k}] \atop s \in [0, n^{-k d_0/{l}}] } \ |X_{t+s} - X_{s} | 
< n^{-k \beta_0},
\\
\sup_{t \in [\tau_j, \tau_j+n^{-k}] \atop s \in [0,1] } \ |Y_{t,s} - X_0 | 
&\leq& |X_{\tau_j} -X_0| + \sup_{t \in [0, n^{-k}] } \ |X_{t+\tau_j} - X_{\tau_j} | 
< \ep_{\ell.k} + n^{-k \beta_0},
\eeas
which complete the proofs of both (\ref{u-1}) and (\ref{u-2}). 
Furthermore, 
\beas
& & 
\inf_{\theta \in \Theta} \inf_{t \in [\tau_j, \tau_j+n^{-k}] \atop s \in [0,1] }
|f(Y_{t,s},\theta)| 
\\
&\geq& \inf_{\theta \in \Theta} |f(X_{\tau_j},\theta)| 
- \sup_{\theta \in \Theta} \sup_{t \in [\tau_j, \tau_j+n^{-k}] \atop s \in [0,1] }
|f(Y_{t,s},\theta) - f(X_{\tau_j},\theta)|. 
\\
&\geq&
n^{-L} 
- \sup_{\theta \in \Theta} \sup_{t \in [\tau_j, \tau_j+n^{-k}] \atop s \in [0,1] }
\left| \int_0^1 \partial_x f(X_{\tau_j} +u(Y_{t,s}-X_{\tau_j}),\theta)du [Y_{s,t}-X_{\tau_j}] \right|
\\
&\geq&
n^{-L} 
- \sup_{\theta \in \Theta} \sup_{t \in [\tau_j, \tau_j+n^{-k}] \atop s, u \in [0,1] }
| \partial_x f(X_{\tau_j} +u(Y_{t,s}-X_{\tau_j}),\theta) | \sup_{t \in [0, n^{-k}]} \ |X_{t+\tau_j} - X_{\tau_j} | 
\\
&\geq&  n^{-L} -C n^{-k \beta_0} \geq \frac{1}{2}n^{-L}
\eeas
for sufficiently large $n$ and $k$, which completes the proof of (\ref{u-3}).

Since it follows from (\ref{u-1}), (\ref{u-2}) and  (\ref{u-3}) that on $\Omega_n \cap A_n \cap B_n^c$, for $\varrho \in (0,1)$,
\beas
& & \sup_{\theta \ne \theta^*} \int_{\tau_j}^{\tau_j+n^{-k}}  |  |f(X_t,\theta)|^\varrho-|f(X_{\tau_j},\theta)|^\varrho | dt
\\
&\leq&  \varrho n^{-k}
\sup_{\theta \in \Theta} \sup_{t \in [\tau_j, \tau_j+n^{-k}] \atop s \in [0,1] }
\{ | f(Y_{t,s},\theta) |^{\varrho-1} |\partial_x f(Y_{t,s},\theta) | |X_t - X_{\tau_j}| \}
\\
&\leq& C \times n^{-k} \times n^{-k \beta_0} \times (2n^L)^{1-\varrho}
\eeas
and that on $\Omega_n \cap A_n \cap B_n^c$, for $\varrho \geq 1$,
\beas
\sup_{\theta \ne \theta^*} \int_{\tau_j}^{\tau_j+n^{-k}}  |  |f(X_t,\theta)|^\varrho-|f(X_{\tau_j},\theta)|^\varrho | dt
&\leq&  C \times n^{-k} \times n^{-k \beta_0},
\eeas
one has that on $\Omega_n \cap A_n \cap B_n^c$, for $\varrho \in (0,\infty)$,  
\beas
\chi_0
&\geq& n^{-k} \times n^{-\varrho L} - C \times n^{-k} \times (n^L)^{1-\varrho} \times n^{-k \beta_0} 
\\
&=& n^{-k -\varrho L} (1-C n^{L -k \beta_0}) \geq \frac{1}{2} n^{-k -\varrho L}
\eeas
for sufficiently large $n$ and $k$.
Thus, for sufficiently large $n$ and $k$, and for $l$ satisfying that $l \geq k d_0/{c_0}$,
\beas
1 - c_1^{-1} \exp (-c_1 n^{c_2}) 
&\leq& P[\Omega_n] 
\\
&\leq& d_1^{-1} \exp (-d_1 n^{k \beta_1/l}) 
+ P \left[ 
\chi_0
\geq \frac{1}{2} n^{-k -\varrho L} \right],
\eeas
and for any $M>0$, there exists $C_M>0$ such that for all $n \in \mathbb{N}$ and sufficiently large $k$,
\beas
P \left[ 
\chi_0
< \frac{1}{2} n^{-k -\varrho L} \right]
\leq \frac{C_M}{n^M}.
\eeas
Therefore, for any $r>0$ satisfying that $2 n^{k+\varrho L} < r \leq 4 n^{k+\varrho L}$,
\beas
P\left[ \chi_0 \leq \frac{1}{r} \right]
\leq P\left[ \chi_0 < \frac{1}{2 n^{k+\varrho L}} \right]
\leq \frac{C_M}{n^M} \leq C_M \left( \frac{4}{r} \right)^{\frac{M}{k+\varrho L} }. 
\eeas
Consequently, for any $K>0$, there exists $C_K>$ such that
$$
P\left[ \chi_0 \leq \frac{1}{r} \right] \leq \frac{C_K}{r^K}
$$
for all $r>0$.
which completes the proof. \qed
\end{en-text}

{\colorb{

When the process $X$ varies in any direction, it finds a nondegenerating direction $\xi$ locally uniformly 
in $(x,\theta)$. 
\begin{description}
\item[[$N_1$\!\!]] \ {\bf (i)} 
{\cdo $\calx_0$ and $\Theta$ are compact.} \vspace{-2mm}
\begin{description}
\item[(ii)] {\cfg $f$ admits a $C^J$-supporting function for each $(x_0,\theta)\in\calx_0\times\Theta$ with 
$\max_{j=0,...,J-1}\big|c_j(x_0,\theta) \big|>0$.}


\item[(iii)] 
{\colorsb{
There exists a stopping time $\tau$ satisfying $\mbox{ess.sup}_\omega\tau<T$ and 
$\mbox{supp}\call\{X_\tau\}\subset(\calx_0)^o$ and there exist 
a positive constant $c_0$ 
and 
distinct positive numbers $\{\alpha_j\}_{j=0,...,J}$ such that $\min_j\alpha_j>c_0/2$ and that 
for every $\xi\in\bbS$ and $\ep>0$, 
the sequence 
\beas 
\bigg(
1-P\bigg[&&\hspace{-6mm}\bigcup_{s\in[\tau,\tau+n^{-c_0}]}
\bigcap_{j=0,...,J}\bigcup_{t\in[0,n^{-c_0}]}
\bigg\{ X_{s+t}-X_s\in\tilde{S}(\xi,\ep) 
\nn\\&&
\mbox{ and }
\big|\xi\cdot (X_{s+t}-X_s)\big|=n^{-\alpha_j} 
{\colorred{\bigg\}}}
\nonumber \bigg]\bigg)_{n\in\bbN}
\eeas
is in $\calp$. 
}}
\end{description}
\end{description}

\begin{en-text} 

\item[(iii)]
For every $\xi\in\bbS$ and $\ep>0$,
there exist $T_0\in(0,T)$, positive constants $c_0,c_1$ and 
distinct positive numbers $\{\alpha_j\}_0^J$ such that 
the sequence 
\bea
\bigg(1-
 P\bigg[ &&
\hspace{-0.8cm}\bigcup_{s\in[0,T_0]}\bigcap_{j=0,...,J} \bigcup_{t\in[0,n^{-c_0}]}
\bigg\{ X_s\in\calx_0,\>
X_{s+t}-X_s\in \tilde{S}(\xi,\ep) 
\nn\\&&
\mbox{ and }
\big|\xi\cdot (X_{s+t}-X_s)\big|=n^{-\alpha_j} 
\nonumber \bigg\}
\bigg]
\bigg)_{n\in\bbN}
\eea
is in $\calp$. 
\end{en-text} 

}}

{\colorb{
\begin{proposition}\label{2306-3}
Under $[N_1]$, $[R]$ and $[C]$,  
Condition $[H2]$ holds. 
\end{proposition}
\proof 
\begin{en-text}
{\cdo We may assume that $\hat{\calx}$ is compact and that for some $\ep>0$, $f$ is defined and 
can be expanded in $x\in B(x_0,\ep)\subset\hat{\calx}$ for every $x_0\in\calx_0$ .}
Since $f\in C(\hat{\calx}\times \Theta)$ for compact sets $\hat{\calx}$ and $\Theta$, 
$c_j:\calx_0\times\Theta\to {\colorred{\bbL^{\otimes J}}}$ ($j=0,...,J-1$)  and $c_J:\calx\times\hat{\calx}\times\Theta\to {\colorred{\bbL^{\otimes J}}}$ are 
bounded continuous. Thus 
\end{en-text}
{\cod Conditions $[N_1]$} (i) and (ii) imply $[N_1^\flat]$, and hence $[N_0^\flat]$ by Lemma \ref{210112-1}. 
{\colorsb{Then $\xi_{\ell,k},\ep_{\ell,k},\calx_\ell$ and $\Theta_{\ell,k}$ are given by $[N_0^\flat]$. 
{\cdo By observing the proof of Lemma \ref{210112-1}, we may take sufficiently small $\calx_\ell$ and 
$\ep_{\ell,k}$. }
{\cod Choose small $\hat{\calx}$ as a neighborhood of $\calx_0$.}

For $\xi=\xi_{\ell,k}$ and $\ep=\ep_{\ell,k}$, 
we apply $[N_1]$(iii) to ensure $[N_0]$(ii) as follows. 
\begin{en-text}
Note that $\calx_0\subset \hat{\calx}^o$ and $f(x,\theta)$ is defined and continuous 
on $\hat{\calx}\times\Theta$.  
\end{en-text}
By the representation of $b_{j,\ell,k}(x_0,\theta,\xi)$ in the proof of Lemma \ref{210112-1}, 
we can replace $\calx_\ell$ by an open set $\check{\calx}_\ell$ 
($\bar{\calx}_\ell\subset\check{\calx}_\ell\subset\overline{\check{\calx}_\ell}\subset\hat{\calx}$) 
in $[N_0^\flat]$.
}}
Denote by $({\mathfrak a}_n)_{n\in\bbN}$ a generic element of $\calp$. It changes from line to line. 
{\colorsb{
Let 
\beas 
B(n,\ell,k,j,s,t) 
&=& 
\bigg\{
X_{s+t}-X_s\in \tilde{S}(\xi_{\ell,k},\ep_{\ell,k}) 
\mbox{ and }
\big|\xi_{\ell,k}\cdot (X_{s+t}-X_s)\big|=n^{-\alpha_j} 
\nonumber \bigg\}.
\eeas
By $[N_1]$ (iii), 
\beas 
P\bigg[
{\colorg{\bigcap_{\ell=1,...,\bar{\ell}}}}
\bigcap_{k=1,...,\bar{k}_\ell}
\bigcup_{s\in[\tau,\tau+n^{-c_0}]}
\bigcap_{j=0,...,J}\bigcup_{t\in[0,n^{-c_0}]}B(n,\ell,k,j,s,t) 
\bigg]
{\colorblue{=}}
1-{\mathfrak a}_n.
\eeas
Let $A_\ell=\{X_{\tau}\in\calx_\ell\}$ for $\ell=1,...,\bar{\ell}$. 
Since 
\beas 
&&\bigg(\bigcup_{\ell}A_{\ell}\bigg)
\bigcap
\bigg(
{\colorg{\bigcap_{\ell}}}
\bigcap_{k}
\bigcup_s\bigcap_{j}\bigcup_tB(n,\ell,k,j,s,t) 
\bigg)
{\colorg{\subset}}
\bigcup_\ell\bigcap_{k}\bigcup_s\bigcap_j\bigcup_t \bigg(A_\ell\cap B(n,\ell,k,j,s,t) \bigg)
\\&&\subset
\bigg[\bigcup_\ell\bigcap_k\bigcup_s\bigcap_j\bigcup_t 
\bigg(\{X_s\in\check{\calx}_\ell\}\cap B(n,\ell,k,j,s,t) \bigg)\bigg]
\bigcup\bigg[\bigcup_\ell\bigcup_s\big\{X_\tau\in\calx_\ell,X_s\not\in\check{\calx}_\ell\big\}\bigg],
\eeas
we have 
\beas
1-{\mathfrak a}_n 
&\leq&
P\bigg[\bigcup_{\ell=1,...,\bar{\ell}}\bigcap_{k=1,...,\bar{k}}\bigcup_{s\in[0,T_0]}
\bigcap_{j=0,...,J}\bigcup_{t\in[0,n^{-c_0}]}
\bigg\{ X_s\in\check{\calx}_\ell,\>
X_{s+t}-X_s\in \tilde{S}(\xi_{\ell,k},\ep_{\ell,k}) 
\nn\\&&
\mbox{ and }
\big|\xi_{\ell,k}\cdot (X_{s+t}-X_s)\big|=n^{-\alpha_j} 
\nonumber \bigg\} 
\bigg]. 
\eeas
This inequality implies $[N_0]$(ii). 
\begin{en-text}
Once again we recall 
the construction of $b_{j,\ell,k}$ and hence $G_{\ell,k}$ given in the proof of Lemma \ref{210112-1}; 
we then see $[N_0]$(i) is satisfied by $[N_0^\flat]$(ii). 
\end{en-text}
{\cdo 
Consequently, we can apply Proposition \ref{230608-1}  to conclude $[H2]$.} \qed \\
}}

{\cred
\begin{remark}\label{241205-2}\rm
{\cred (i)} Corresponding to Remark \ref{241205-1}, if we assume 
Condition $[N_1]$ (iii) with $\tilde{S}(\xi_{\ell,k},\ep_{\ell,k})$ replaced by $S(\xi_{\ell,k},\ep_{\ell,k})$, 
then the inequality in Lemma \ref{210102-2} with $D(\xi_{\ell,k},\ep_{\ell,k})$ for 
$\tilde{D}(\xi_{\ell,k},\ep_{\ell,k})$ 
is still sufficient to prove the same result as Proposition \ref{2306-3}. 
Note that Lemma \ref{210102-2} was implicitly used in the proof of Proposition \ref{2306-3} through Proposition \ref{230608-1}. 
{\cred 
(ii) For Proposition \ref{2306-3}, we can eliminate the Lipschitz continuity condition in $[R]$ because 
we can replace $f$ by the supporting function in the last part of the proof of Proposition \ref{230608-1}. }
\end{remark}
}

\begin{en-text}
Next, we assume $[N_1]$ and $[C]$.
Note that for any $\alpha>0$, there exist positive constants $c_0$ and $\beta_0$ such that $c_0 \beta_0 <\alpha$.
Let $\ep_n = 1- n^{-(\alpha-c_0 \beta_0)}$. 
Owing to Lemma \ref{210112-1}, 
$[H2]$ holds if (\ref{u-n1}) implies (\ref{u-n0}).
Indeed, by $[N_1]$, for each $\xi \in \mathbb{S}$ and $\alpha>0$, there exist positive constants $c_0$,$c_1$ and $c_2$ such that
\beas
& & 1-c_1^{-1}\exp( -c_1 n^{-c_2} ) 
\\
&\leq& \inf_{x_0 \in \calx_0} P[ \ X_t -x_0 \in \mathbb{L} \mbox{ and } |\xi \cdot(X_t-x_0)| = n^{-\alpha}
\mbox{ for some } t \in [0,n^{-c_0}] \ | \ X_0=x_0 ]
\\
&\leq&
\inf_{x_0\in\calx_0}
 P\bigg[ 
\bigg\{ X_t-x_0 \in S(\xi,\ep_n) \mbox{ and }
\xi \cdot (X_t-x_0)=n^{-\alpha} \bigg\} \mbox{ or } 
\nonumber 
\\
& &
\hspace{1.5cm} 
\textcolor{black}{
\bigg\{ X_t-x_0\in S(-\xi, \ep_n) \mbox{ and } 
-\xi \cdot (X_t-x_0)=n^{-\alpha} \bigg\}
}
\mbox{ for some }
t\in[0,n^{-c_0}]
\ \bigg| \  
X_0=x_0 \ \bigg]
\\
& & + \sup_{x_0 \in \calx_0} P[ \ X_t -x_0 \in S(\xi,\ep_n)^c \mbox{ and } \xi \cdot(X_t-x_0) = n^{-\alpha}
\mbox{ for some } t \in [0,n^{-c_0}] \ | \ X_0=x_0 ]
\\
& & + \sup_{x_0 \in \calx_0} P[ \ X_t -x_0 \in S(-\xi,\ep_n)^c \mbox{ and } -\xi \cdot(X_t-x_0) = n^{-\alpha}
\mbox{ for some } t \in [0,n^{-c_0}] \ | \ X_0=x_0 ].
\eeas     
For all $x_0 \in \calx_0$, 
\beas
& & P[ \ X_t -x_0 \in S(\xi,\ep_n)^c \mbox{ and } \xi \cdot(X_t-x_0) = n^{-\alpha}
\mbox{ for some } t \in [0,n^{-c_0}] \ | \ X_0=x_0 ]
\\
&\leq& P[ |X_t-x_0| \geq n^{-c_0 \beta_0}
\mbox{ for some } t \in [0,n^{-c_0}] \ | \ X_0=x_0 ]
\\
& & + P[ \ \xi \cdot(X_t-x_0) < (1-\ep_n) |X_t -x_0| \mbox{ and } \xi \cdot(X_t-x_0) = n^{-\alpha}  \mbox{ and } |X_t-x_0| < n^{-c_0 \beta_0}
\\
& & \hspace{0.7cm} 
\mbox{ for some } t \in [0,n^{-c_0}] \ | \ X_0=x_0 ]
\\
& & + P[ \ |X_t -x_0| > \ep_n \mbox{ and } \xi \cdot(X_t-x_0) = n^{-\alpha}  \mbox{ and } |X_t-x_0| < n^{-c_0 \beta_0}
\\
& & \hspace{0.7cm} 
\mbox{ for some } t \in [0,n^{-c_0}] \ | \ X_0=x_0 ]
\\
&\leq& d_1^{-1} \exp \left( -d_1 n^{c_0 \beta_1} \right) 
\\
& &  + P[ \ \xi \cdot(X_t-x_0) < n^{-\alpha} \mbox{ and } \xi \cdot(X_t-x_0) = n^{-\alpha}  
\mbox{ for some } t \in [0,n^{-c_0}] \ | \ X_0=x_0 ]
\\
& & + P[ \ |X_t -x_0| > 1- n^{-(\alpha-c_0 \beta_0)}  \mbox{ and } |X_t-x_0| < n^{-c_0 \beta_0}
\mbox{ for some } t \in [0,n^{-c_0}] \ | \ X_0=x_0 ]
\\
&\leq& d_1^{-1} \exp \left( -d_1 n^{c_0 \beta_1} \right)
\eeas
for sufficiently large $n$.
In the same way, for all $x_0 \in \calx_0$, 
\beas
& &  P[ \ X_t -x_0 \in S(-\xi,\ep_n)^c \mbox{ and } -\xi \cdot(X_t-x_0) = n^{-\alpha}
\mbox{ for some } t \in [0,n^{-c_0}] \ | \ X_0=x_0 ]
\\
&\leq& d_1^{-1} \exp \left( -d_1 n^{c_0 \beta_1} \right)
\eeas
for sufficiently large $n$.
Consequently, for each $\xi \in \mathbb{S}$ and $\alpha>0$, there exist positive constants $c_0$,$\tilde{c}_1$ and $\tilde{c}_2$ such that
\beas
& & 1-\tilde{c}_1^{-1}\exp( -\tilde{c}_1 n^{-\tilde{c}_2} ) 
\\
&\leq&
\inf_{x_0\in\calx_0}
 P\bigg[ 
\bigg\{ X_t-x_0 \in S(\xi,\ep_n) \mbox{ and }
\xi \cdot (X_t-x_0)=n^{-\alpha} \bigg\} \mbox{ or } 
\nonumber 
\\
& &
\hspace{1.5cm} 
\textcolor{black}{
\bigg\{ X_t-x_0\in S(-\xi, \ep_n) \mbox{ and } 
-\xi \cdot (X_t-x_0)=n^{-\alpha} \bigg\}
}
\mbox{ for some }
t\in[0,n^{-c_0}]
\ \bigg| \  
X_0=x_0 \ \bigg],
\eeas
which completes the proof.
\qed \\
\end{en-text}

\begin{en-text}
Let
\beas
C_n^{(1)} 
&=&
\left\{ X_t -X_0 \in \mathbb{L} \mbox{ and } \xi \cdot(X_t-X_0) = n^{-\alpha}
\mbox{ for some } t \in [0,n^{-c_0}] \right\}, 
\\
C_n^{(2)} 
&=&
\left\{  X_t -X_0 \in \mathbb{L} \mbox{ and } -\xi \cdot(X_t-X_0) = n^{-\alpha}
\mbox{ for some } t \in [0,n^{-c_0}] \right\}, 
\\
D_n^{(1)} 
&=&
\left\{ |X_t -X_0| \geq n^{-\beta_0} \mbox{ and } -\xi \cdot(X_t-X_0) = n^{-\alpha}
\mbox{ for some } t \in [0,n^{-c_0}] \right\}, 
\\
D_n^{(2)} 
&=&
\left\{ |X_t -X_0| < n^{-\beta_0} \mbox{ and } -\xi \cdot(X_t-X_0) = n^{-\alpha}
\mbox{ for some } t \in [0,n^{-c_0}] \right\}, 
\\
E_n^{(1)} 
&=&
\left\{ |X_t -X_0| \geq n^{-\beta_0}  \mbox{ and } \xi \cdot(X_t-X_0) = n^{-\alpha}
\mbox{ for some } t \in [0,n^{-c_0}] \right\}, 
\\
E_n^{(2)} 
&=&
\left\{ |X_t -X_0| < n^{-\beta_0} \mbox{ and } -\xi \cdot(X_t-X_0) = n^{-\alpha}
\mbox{ for some } t \in [0,n^{-c_0}] \right\}, 
\\
F_n^{(1)} 
&=&
\left\{ X_t -X_0 \in S(\xi, \ep) \mbox{ and } \xi \cdot(X_t-X_0) = n^{-\alpha}
\mbox{ for some } t \in [0,n^{-c_0}] \right\},
\\
F_n^{(2)} 
&=&
\left\{ X_t -X_0 \in S(-\xi, \ep) \mbox{ and } -\xi \cdot(X_t-X_0) = n^{-\alpha}
\mbox{ for some } t \in [0,n^{-c_0}] \right\}.
\eeas
By $[C]$ with $k=1$, there exist positive constants $\beta_0$, $\beta_1$, $d_0$ and $d_1$ such that
\beas
& &
\sup_{x_0 \in \calx_0} P\left[  \left. \ \sup_{t \in [0,n^{-1}]} |X_t -X_0| \geq n^{-\beta_0} \ \right| \ X_0=x_0 \right]
\\
&\leq&
\sup_{x_0 \in \calx_0} P\left[  \left. \ \sup_{t \in [0,n^{-1}] \atop s\in [0, n^{-d_0}]} |X_{t+s} -X_s| \geq n^{- \beta_0} \ \right| \ X_0=x_0 \right]
\leq
d_1^{-1}\exp( -d_1 n^{\beta_1} ).
\eeas
By $[N_1]$, for each $\xi \in \mathbb{S}$ and $\alpha>0$, there exist positive constants $c_0$,$c_1$ and $c_2$ such that
\beas
& & 1-c_1^{-1}\exp( -c_1 n^{-c_2} ) 
\\
&\geq& \inf_{x_0 \in \calx_0} P[ \ X_t -X_0 \in \mathbb{L} \mbox{ and } |\xi \cdot(X_t-X_0)| = n^{-\alpha}
\mbox{ for some } t \in [0,n^{-c_0}] \ | \ X_0=x_0 ].
\eeas
We then have that for any $\alpha$ satisfying that $\alpha > \beta_0$ and for sufficiently large $n$,
\beas
& & 1-c_1^{-1}\exp( -c_1 n^{-c_2} ) 
\\
&\geq& P[ \ X_t -X_0 \in \mathbb{L} \mbox{ and } |\xi \cdot(X_t-X_0)| = n^{-\alpha}
\mbox{ for some } t \in [0,n^{-c_0}] \ | \ X_0=x_0 ]
\\
&\geq& 
P[ \ |X_t -X_0| < n^{-\beta_0} \mbox{ and } \xi \cdot(X_t-X_0) = n^{-\alpha}
\mbox{ for some } t \in [0,n^{-c_0}] \ | \ X_0=x_0 ]
\\
& & +
P[ \ |X_t -X_0| < n^{-\beta_0} \mbox{ and } -\xi \cdot(X_t-X_0) = n^{-\alpha}
\mbox{ for some } t \in [0,n^{-c_0}] \ | \ X_0=x_0 ]
\\
& & +
P[ \ |X_t -X_0| \geq n^{-\beta_0} \mbox{ for some } t \in [0,n^{-c_0}] \ | \ X_0=x_0 ]
\\
&\geq& 
P[ \ X_t -X_0 \in S(\xi, 1-n^{-\alpha+\beta_0})  \mbox{ and } \xi \cdot(X_t-X_0) = n^{-\alpha}
\mbox{ for some } t \in [0,n^{-c_0}] \ | \ X_0=x_0 ]
\\
& & +
P[ \ X_t -X_0 \in S(-\xi, 1-n^{-\alpha+\beta_0})  \mbox{ and } -\xi \cdot(X_t-X_0) = n^{-\alpha}
\mbox{ for some } t \in [0,n^{-c_0}] \ | \ X_0=x_0 ]
\\
& & +
c_1^{-1}\exp( -c_1 n^{c_0 \beta_1} ).
\eeas
\qed\\
\end{en-text}

\begin{en-text}
Then there exists a sequence of events $\Omega_n$ 
satisfying $P[\Omega_n]>1-c_1^{-1}e^{-c_1n}$ ($n\in\bbN$) for some $c_1>0$ and 
such that 
\beas 
\max_{j=0,...,J}
|F(X_0,X_{\tau_j},\theta,X_{\tau_j}-X_0)|
&\geq 
n^{-L}
\eeas
for all $n\geq n_0$, where 
$L$ and $n_0$ are given in Lemma \ref{210102-2}. 
%
%
%
Once (the original functional of) 
$F$ gained order $n^{-L}$, then the integration in time 
turns out to keeps the amount of order $n^{-L-c_2}$ for some $c_2>0$ 
except for a large deviation probability due to [C]. 
\end{en-text}
\begin{en-text}
Also, note that for $\xi_1\in\bbR^{\sf d}$ 
and for the multi-index $\alpha$, 
there exists $a>0$ such that 
$a|\xi_1\cdot\xi|^{|\alpha|}\leq |\xi^\alpha|
\leq a^{-1}|\xi_1\cdot\xi|^{|\alpha|}$ for all $\xi\in C_{\xi_1}$ 
if the cone has small angle. 
This observation helps to find $G(x_0,x,\theta,\xi)$ 
in multi-dimensional case. 
It is also possible to replace the above criterion 
by a similar but slightly simpler condition in 
case ${\sf d}=1$.
\qed\\ 
\end{en-text}

{\colorb{
\section{Proof of Theorems \ref{230614-1} and \ref{230614-2}}\label{230806-3}
\begin{en-text}
The ingredients of the proof of these results are the polynomial type large deviation inequality 
for the quasi likelihood random field, as well as limit theorems for semimartingales. 
In order to ensure the the polynomial type large deviation inequality, 
it is necessary to prove the nondegeneracy of a random index $\chi_0$. 
To answer this question, we will prepare a new machinery to induce the nondegeneracy of the statistical random field 
by connecting a nondegeneracy of a newly introduced tesor field over the statistical manifold 
and the nondegeneracy of the underlying stochastic process. 

{\it Proof of Theorems \ref{230614-1}and \ref{230614-2}. } 
\end{en-text}
It suffices to verify ${\cred [H1^\sharp]}$ and $[H2]$ due to Theorems \ref{thm3} and \ref{thm2}. 
We take $\bar{\Theta}$ for ``$\Theta$'' in $[N_1]$. 
Condition ${\cred [H1^\sharp]}$ is obviously satisfied under the assumptions. 

We shall show $[H2]$. For this, we will apply Proposition \ref{2306-3}. 
Let $\calx_0$ and $\hat{\calx}$ {\colorred be} compact sets in $\bbR^{{\sf d}}$ such that 
$\mbox{supp}\{X_\tau\}\subset \calx_0^o \subset \calx_0 \subset \hat{\calx}^o
\subset \hat{\calx}\subset U$. 
Condition $[C]$ is easily verified for $\beta_0\in(0,1/2)$. 
Conditions $[R]$, $[N_1](i),(ii)$ are obvious. 

Take numbers $c_0$, $c_1$, $\alpha_j$ ($j=0,...,J$) such that 
$\alpha_0>\cdots>\alpha_J>c_1/2>c_0/2$. 
\begin{en-text}
Since the support of $\call\{X_0\}$ is compact, 
\beas 
\sup_{x_0\in\mbox{supp}\call\{X_0\}}P\big[ \sup_{t\in[0,n^{-c_0}]}|X_t|\geq C\big| X_0=x_0\big]
\leq C_L n^{-L}
\eeas
for all $n\in\bbN$ and for every $L$ if we choose sufficiently large constant $C_L$, where 
$C$ is a number larger than the diameter of the support. 
Therefore we may assume 
\end{en-text}
Let $T(\xi_0,\ep)=C(\xi_0;\ep)\cap\{\xi\in\bbL;{\colorred 3}<\xi\cdot\xi_0<{\colorred 4}\}$ for $\ep>0$, 
{\cred $T(\xi_0,\ep,\eta)=\{\eta \xi;\xi\in T(\xi_0,\ep)\}$} 
for $\eta>0$. 
Let $s(n,k)=\tau+kn^{-c_1}$ for $k=1,...,k(n)$, $k(n)=[n^{c_1-c_0}]$. 
Let $s(n,k,-1)=s(n,k)$,
$s(n,k,0)=s(n,k)+n^{-2\alpha_0}$
and $s(n,k,j)=s(n,k,j-1)+n^{-2\alpha_j}$ for $j=1,...,J$. 
Obviously, $s(n,k)$ and $s(n,k,j)$ are stopping times. 
We may assume that $s(n,k,J)\leq s(n,k+1)$. 

Let $\ep_0>0$. 
Let 
\beas 
A^\bbX(\xi,\ep_0,n,k,j) &=& 
\bigg\{\bbX_u-\bbX_{s(n,k,j-1)}\in {\cred T}(\xi,\ep_0, n^{-\alpha_j}) \mbox{ for some } u\in (s(n,k,j-1),s(n,k,j)]\bigg\}
\\&&
\bigcap 
\bigg\{
\sup_{u\in (s(n,k,j-1),s(n,k,j)]}|P^\perp_\xi (\bbX_u-\bbX_{s(n,k,j-1)})|<\ep_0 n^{-\alpha_j}\bigg\}
\eeas
for a process $\bbX$, 
where $P^\perp_\xi:\bbL\to\bbL$ is the orthogonal projection on $\bbL$ to the subspace orthogonal to $\xi$. 
We write $X_t$ for $t\geq\tau$ as 
$
X_t = X_\tau +M_t+R_t
$, 
where 
\beas 
M_t &=& a_\tau (w_t-w_\tau)+\tilde{a}_\tau(\tilde{w}_t-\tilde{w}_\tau)
\eeas
and 
\beas 
R_t &=& 
\int_\tau^t (a_s-a_\tau) dw_s+\int_\tau^t (\tilde{a}_s-\tilde{a}_\tau)d\tilde{w}_s
+\int_\tau^t \tilde{b}_sds.
\eeas
The process
$(M_{t-\tau})_{t\in[\tau,\tau+n^{-c_0}]}$ has the same law on $C([0,n^{-c_0}];\bbR^{{\sf d}})$ 
as $((a_\tau^{\otimes2}+\tilde{a}_\tau^{\otimes2})^{1/2}B_t)_{t\in[0,n^{-c_0}]}$, where 
$B_t$ is a ${\sf d}$-dimensional standard Wiener process independent of $\calf_\tau$.  
Using the scaling property and independency between increments of the Wiener process, and also a classical result of the distribution of its absolute deviation 
or a support theorem, 
it is easy to see 
\beas 
q:=\mbox{ess.inf}_{\omega,\xi\in\bbS,n,k,j}P[A^M(\xi,\ep_0,n,k,j)|\calf_{s(n,{\colorred{k,j-1)}}}]>0.
\eeas
{\colorred{It should be noted that the uniform (in $\omega$) boundedness and the uniform (in $\omega$) nondegeneracy of the matrix 
$(a_\tau^{\otimes2}+\tilde{a}_\tau^{\otimes2})^{1/2}$ was used to control random linear transform of the Brownian motion $B_t$. }}

For any $\ep_1,\ep_2>0$,  there exists $n_0\in\bbN$ such that 
\beas 
\mbox{ess.sup}_\omega\sup_{n\geq n_0\atop {k=1,...,k(n) \atop j=0,..,,J}}
P\bigg[\sup_{t\in(s(n,k,j-1),s(n,k,j)]}\big|n^{\alpha_j}(R_t-R_{s(n,k,j-1)})\big|>\ep_1
|\calf_{s(n,{\colorred{k,j-1)}})}\bigg]<\ep_2. 
\eeas
Indeed, Lenglart's inequality gives uniform estimates for stochastic integrals with the aid of the right-continuity of 
$a$ and $\tilde{a}$ as well as the $L^p$-boundedness provided in $[A2]$. 
For the integral of $\tilde{b}$, the H\"older inequality with $L^p$-estimate for $\tilde{b}_t$ yields the estimate. 
{\colorred{
In order to check $[N_1]$(iii), we consider arbitrary $\xi\in\bbS$ and $\ep>0$. 
We choose positive constants $\ep_0$, $\ep_1$ and $\ep_2$ such that $\ep_1<<\ep_0<<\ep$ and $\ep_2<q$. Then 
\beas 
\mbox{ess.inf}_\omega\inf_{n\geq n_0\atop {k=1,...,k(n) \atop j=0,..,,J}}
P\big[A(\xi,\ep,n,k,j)\big|\calf_{s(n,k,j-1)}\big]
>q-\ep_2=:q'>0,
\eeas
where $A(\xi,\ep,n,k,j)$ is the event defined in the same way as $A^{\bbX}$ for $\bbX=X$ and $\ep_0=\ep$,  
with $T(\xi_0,\ep)$ replaced by 
$C(\xi_0;\ep)\cap\{\xi\in\bbL;2<\xi\cdot\xi_0<5\}$. 

Let $\ep_3 >0$.
Let 
\beas 
B(n,k,j-1)
&=&
\bigg\{n^{\alpha_j}|X_{s(n,k)}-X_{s(n,k,j-1)}| <\ep_3 \bigg\}. 
\eeas
We see that for $\ep_4>0$, there exists $n_1\geq n_0$ such that 
\beas
\mbox{ess.sup}_\omega\sup_{n\geq n_1\atop {k=1,...,k(n) \atop j=0,..,,J}}
P[B(n,k,j-1)^c|\calf_{s(n,k)}]
&<&
\ep_4.
\eeas
Here the ordering $\alpha_0>\cdots>\alpha_J$ was used. 

Since 
\beas 
&&
{\colorg{P}} \bigg[B(n,k,j-1)\cap A(\xi,\ep,n,k,j)\big|\calf_{s(n,k,j-1)}\bigg]
\\&\geq&
P\big[A(\xi,\ep,n,k,j)\big|\calf_{s(n,k,j-1)}\big]-1_{B(n,k,j-1)^c},
\eeas
we have 
\beas
&&
P\bigg[\bigcap_{j=0}^{j'}\bigg(B(n,k,j-1)\cap A(\xi,\ep,n,k,j)\bigg)\bigg|\calf_{s(n,k)}\bigg]
\\&{\cred \geq}&
P\bigg[\prod_{j=0}^{j'-1}1_{(B(n,k,j-1)\cap A(\xi,\ep,n,k,j)}\>
\bigg(
P\big[A(\xi,\ep_0,n,k,j')\big|\calf_{s(n,k,j'-1)}\big]-1_{B(n,k,j'-1)^c}\bigg)\bigg|\calf_{s(n,k)}\bigg]
\\&\geq&
q'P\bigg[\bigcap_{j=0}^{j'-1}\bigg(B(n,k,j-1)\cap A(\xi,\ep,n,k,j)\bigg)\bigg|\calf_{s(n,k)}\bigg]-\ep_4
\eeas
for all $n\geq n_1$, $k$, $j$ and a.s. $\omega$. 
Let 
\beas 
C(\xi,\ep,n,k)&=&
\bigcap_{j=0}^{J}\bigg(B(n,k,j-1)\cap A(\xi,\ep,n,k,j)\bigg). 
\eeas
We use the above inequality repeatedly to obtain 
\beas 
\mbox{ess.inf}_\omega\inf_{n\geq n_1\atop {k=1,...,k(n)}}
P\bigg[C(\xi,\ep,n,k)\bigg|\calf_{s(n,k)}\bigg]
&\geq&
(q')^{J+1}-\ep_4\sum_{j=0}^{J-1}(q')^j
\geq q''
\eeas
with some positive constant $q''$ 
if we take a sufficiently small $\ep_4$ for $q'$. 
Similarly by conditioning, 
\beas &&
P\bigg[\bigcap_{k=1,...,k(n)} C(\xi,\ep,n,k)^c\bigg]
\\&=&
P\bigg[\bigcap_{k=1,...,k(n)-1} C(\xi,\ep,n,k)^c\> P\big[C(\xi,\ep,n,k(n))^c|\calf_{s(n,k(n))}\big]\bigg]
\\&\leq&
(1-q'')P\bigg[\bigcap_{k=1,...,k(n)-1} C(\xi,\ep,n,k)^c\bigg]
\\&\leq& \cdots
\\&\leq&
(1-q'')^{k(n)}
\eeas
for $n\geq n_1$. 
}}

\begin{en-text}
In order to check $[N_1]$(iii), we consider arbitrary $\xi\in\bbS$ and $\ep>0$. 
We choose sufficiently small $\ep_0,\ep_1,\ep_2>0$ and use continuity of $X$ to observe 
\beas 
\inf_{\xi\in\bbS, n,k,j}
\mbox{ess.inf}_\omega && 
\hspace{-5mm} P\bigg[
X_{s(n,k,j-1)}\in\calx_0,\>
X_t-X_{s(n,k,j-1)}\in \tilde{S}(\xi,\ep) 
\nn\\&&
\mbox{ and }
\big|\xi\cdot (X_t-X_{s(n,k,j-1)})\big|=n^{-\alpha_j} 
\\&&
\mbox{ for some }t\in(s(n,k,j-1),s(n,k,j)]
\bigg| \calf_{s(n,k,j-1)}\bigg]
\\&&
>q-\ep_2=:q'>0
\eeas
for large $n$. 
\end{en-text}
{\colorred{
{\colorblue We choose a sufficiently small $\ep_3$.}
Now it is easy to see that for large $n$, 
$X_t-X_{{\colorblue s(n,k)}}\in {\cred S}(\xi,\ep) $ and 
$\big|\xi\cdot (X_t-X_{{\colorblue s(n,k)}})\big|=n^{-\alpha_j} $ 
for some {\colorblue $t=t(j)\in(s(n,k),s(n,k)+n^{-{\cred c_1}}]$ for every $j$} 
on the event $C(\xi,\ep,n,k)$.
%
\begin{en-text}
\beas 
P\bigg[ \ &&
\hspace{-0.8cm}\bigcup_{k=1,...,k(n)}\bigcap_{j=0,...,J} \bigcup_{t\in[0,n^{-c_0}]}
\bigg\{ X_s\in\calx_0,\>
X_{s+t}-X_s\in \tilde{S}(\xi,\ep) 
\nn\\&&
\mbox{ and }
\big|\xi\cdot (X_{s+t}-X_s)\big|=n^{-\alpha_j} 
\nonumber \bigg\}
\bigg]
\bigg)_{n\in\bbN}
\\&\geq&
1-C(1-{q'}^{J+1})^{k(n)}
\eeas  
\end{en-text} 
\begin{en-text}
Therefore, for any $\kappa>0$, 
\bea\label{230616-1}
P\bigg[ \ &&
\hspace{-0.8cm}\bigcup_{k=1,...,k(n)}\bigcap_{j=0,...,J} \bigcup_{t\in[0,n^{-c_1}]}
\bigg\{ X_{s(n,k)}\in\calx_0,\>
X_{s(n,k)+t}-X_{s(n,k)}\in \tilde{S}(\xi,\ep) 
\nn\\&&
\mbox{ and }
\big|\xi\cdot (X_{s(n,k)+t}-X_{s(n,k)})\big|=n^{-\alpha_j} 
\nonumber \bigg\}
\bigg]
\nn\\\geq&&
P\bigg[ \ 
\bigcup_{k=1,...,k(n)}C(\xi,\ep,n,k)\bigg]
-P\bigg[\bigcup_{t\in(\tau,\tau+n^{-2c_0}]}\big\{X_t\not\in\calx_0\big\}\bigg]
\nn\\&\geq&
1-C(1-q'')^{k(n)}-Cn^{-\kappa}
\eea
for $n\in\bbN$ and for some constant $C$. 
\end{en-text}
{\colorsb{\noindent 
Therefore, 
\bea\label{230616-1}
P\bigg[ \ &&
\hspace{-0.8cm}\bigcup_{k=1,...,k(n)}\bigcap_{j=0,...,J} \bigcup_{t\in[0,n^{-c_1}]}
\bigg\{
X_{s(n,k)+t}-X_{s(n,k)}\in {\cred S}(\xi,\ep) 
\nn\\&&
\mbox{ and }
\big|\xi\cdot (X_{s(n,k)+t}-X_{s(n,k)})\big|=n^{-\alpha_j} 
\nonumber \bigg\}
\bigg]
\nn\\\geq&&
P\bigg[ \ 
\bigcup_{k=1,...,k(n)}C(\xi,\ep,n,k)\bigg]
\nn\\&\geq&
1-C(1-q'')^{k(n)}
\eea
for $n\in\bbN$ and for some constant $C$.
%
\begin{en-text}
Here the ordering $\alpha_0>\cdots>\alpha_J$ was used to deform the initial value of the increments of $X$ 
in question 
from $X_{s(n,k,j-1)}$ to $X_{s(n,k)}$; the existence of a suitable hitting time of the increment is deduced 
even after the change of the initial values. 
\end{en-text}
\begin{en-text}
This inequality implies $[N_1]$(iii). 
Indeed, if $\omega\in\Omega$ is in the event on the left-hand side of (\ref{230616-1}), 
then set $\tau^*(\omega)=s(n,k)$ with minimum $k$ for $\omega$. If there is no such number $k$, then set 
$\tau^*(\omega)=\tau(\omega)$; note that we do not require $\tau^*$ is a stopping time. 
\end{en-text}
%
This inequality implies $[N_1]$(iii). 
}}
}}
 \qed

\begin{en-text}
{\colorb{
\section{Nondegenerate diffusion}
For the elliptic diffusion process $X$, 
it does not matter to verify 
the second half of Condition [$N_1$]. 
More precisely, we will assume that $X$ satisfies the stochastic differential 
equation ..... 

\koko}}

{\colorb{Suppose that $\mbox{supp}\call\{X_0\}$ is compact, as in Section \ref{230608-3}. 
Let $\calx_0=\mbox{supp}\{X_0\}$ and let $\hat{\calx}$ be a compact set in $\bbR^{\sf d}$ 
satisfying $\calx_0\subset\hat{\calx}^o$. 
Let $X=(X_t)$ be a $\sf d$-dimensional diffusion process satisfying 
the stochastic integral equation
\begin{eqnarray}
X_t &=& X_0+\int_0^t V_0(X_s)ds +\int_0^t V(X_s)dW_s, \quad t \in [0,T], \label{u-4}
\end{eqnarray}
where $W$ is an ${\sf r}$-dimensional standard Wiener process,
$V_0$ is an ${\bbR}^{\sf d}$-valued function defined on ${\bbR}^{\sf d}$
and
$V$ is an ${\bbR}^{\sf d} \otimes {\bbR}^{\sf r}$-valued function defined on ${\bbR}^{\sf d}$.
}}

\begin{proposition} 
Suppose that $[N_1^\flat]$ is satisfied.
For the diffusion process $X=(X_t)$  satisfying (\ref{u-4}),
assume that there exists $K>0$ such that for all $x,y$,
$$
|V_0(x)-V(y)| + |V(x)-V(y)| \leq K|x-y|
$$
and that $\inf_x \det([VV^\star](x)) >0$. 
Then $[N_1]$ holds as well as $[C]$, 
in particular, $[H2]$ holds true. 
\end{proposition}
\proof 
Let
\beas
C_n &=& \left\{ \sup_{t \in [0,n^{-c_0}]} |\xi \cdot (X_t-X_0)| \geq n^{-\alpha} \right\}, \\
D_n &=& \left\{ \sup_{t \in [0,n^{-c_0}]} |X_t-x_0| \geq n^{-\beta_0} \right\}. 
\eeas
Since
\beas
& & \inf_{x_0 \in \calx_0} P \left[ \left. C_n \ \right| X_0 =x_0 \right]
\\
&\leq& 
\inf_{x_0 \in \calx_0} P[ \ |X_t -x_0| < n^{-\beta_0} \mbox{ and } |\xi \cdot(X_t-x_0)| = n^{-\alpha}
\mbox{ for some } t \in [0,n^{-c_0}] \ | \ X_0=x_0 ]
\\
& & + \sup_{x_0 \in \calx_0} P \left[ \left. D_n \ \right| X_0 =x_0 \right]
\\
&\leq& 
\inf_{x_0 \in \calx_0} P[ \ X_t -x_0 \in \mathbb{L} \mbox{ and } |\xi \cdot(X_t-x_0)| = n^{-\alpha}
\mbox{ for some } t \in [0,n^{-c_0}] \ | \ X_0=x_0 ]
\\
& & + \sup_{x_0 \in \calx_0} P \left[ \left. D_n \ \right| X_0 =x_0 \right] 
\eeas
and
\beas
\sup_{x_0 \in \calx_0} P \left[ \left. C_n^c \ \right| X_0 =x_0 \right]
\leq
\sup_{x_0 \in \calx_0} P[ \ C_n^c \cap D_n^c 
\ | \ X_0=x_0 ]
+ \sup_{x_0 \in \calx_0} P \left[ \left. D_n \ \right| X_0 =x_0 \right], 
\eeas
one has that
\bea
& &
\inf_{x_0 \in \calx_0} P[ \ X_t -x_0 \in \mathbb{L}  \mbox{ and } |\xi \cdot(X_t-x_0)| = n^{-\alpha}
\mbox{ for some } t \in [0,n^{-c_0}] \ | \ X_0=x_0 ]
\nonumber
\\
&\geq& 1- \sup_{x_0 \in \calx_0} P[ \ C_n^c \cap D_n^c \ | \ X_0=x_0 ]
-2 \sup_{x_0 \in \calx_0} P \left[ \left. D_n \ \right| X_0 =x_0 \right].  \label{u-p2-1}
\eea

We will estimate $P \left[ \left. D_n \ \right| X_0 =x_0 \right]$.
Fix $c_0$.
Let $\tau_0 = \left\{ t>0 \ ; \ |X_t-X_0| \geq n^{-\beta_0} \right\}$.
For each $i=1, \ldots, {\sf d}$, there exists a Brown motion $B$ such that $\int_0^s V^{i \cdot} (x_u) dW_u = B(A_s^i), s \geq 0$,
where $A_s^i = \int_0^s [V V^\star]^{i i}(X_u)du$.
It follows from Grownwall's lemma that 
$$
\sup_{t \in [0,\tau_0]} |X_t-X_0| \leq C (\tau_0 (1+|X_0|) + C \sup_{t \in [0,\tau_0]} \left| \int_0^t V(X_s) dW_s \right|.
$$
We then have that
for sufficiently large $n$, and for all $\beta_0 \in (0, \frac{c_0}{2})$,
\beas
& & P \left[ \left.  \sup_{t \in [0,n^{-c_0}]} |X_t-x_0| \geq n^{-\beta_0}  \ \right| X_0 =x_0 \right]
\\
&=& P \left[ \left. \{ \tau_0 \leq n^{-c_0} \} \cap \left\{ n^{-\beta_0} -C_1 n^{-c_0} \leq 
\sum_{i=1}^{\sf d} \sup_{t \in [0,\tau_0]} \left| \int_0^t V^{i \cdot} (x_u) dW_u \right| \right\} \ \right| X_0 =x_0 \right]
\\
&\leq& {\sf d} \times P \left[ \left.  \frac{n^{-\beta_0}}{2 {\sf d}} \leq 
\sup_{u \in [0,C_2 n^{-c_0}]} |B(u)|  \ \right| X_0 =x_0 \right] 
\leq 2 {\sf d} \exp \left(- \frac{ \left( n^{-\beta_0}/{2 {\sf d}} \right)^2 }{2n^{-c_0}} \right)
\eeas
for all $x_0 \in \calx_0$.
Thus, for each $c_0>0$, there exist positive constants $\beta_0$ ($\beta_0 < \frac{c_0}{2}$) and $\gamma_1$  
such that
\bea
\sup_{x_0 \in \calx_0} P \left[ \left. D_n \ \right| X_0 =x_0 \right] \leq \gamma_1^{-1} \exp (-\gamma_1 n^{c_0 -2\beta_0}) \label{u-p2-2}
\eea
for all $n \in \mathbb{N}$.

We will estimate $P[ \ C_n^c \cap D_n^c \ | \ X_0=x_0 ]$. Fix $\alpha$. 
On $C_n^c \cap D_n^c$, for all $\beta_0$ and $c_0$
satisfying that $ \beta_0 < \frac{c_0}{2} < \alpha$,
\beas
& & \sup_{t \in [0,n^{-c_0}]} \left| \xi \cdot \left( X_t -X_0 -\int_0^t V_0(X_s) ds \right) \right|
\\
&\leq& \sup_{t \in [0,n^{-c_0}]} \left|  \xi \cdot (X_t -X_0) \right| 
+ \sup_{t \in [0,n^{-c_0}]} \left| \xi \cdot \int_0^t V_0(X_s) ds \right|
\\
&\leq& n^{-\alpha} + C(1+|X_0|)n^{-c_0} + Cn^{-c_0} n^{-\beta_0}
\\
&\leq& n^{-c_0/{2}}(n^{-(\alpha-c_0/2)}+C_1 n^{-c_0/2}) \leq C_2 n^{-c_0/2} n^{-\beta_1/2},
\eeas
where $\beta_1= \min \{ 2\alpha-c_0, c_0 \}$.
Furthermore, for some $\ep_1>0$,
\beas
\sup_{t \in [0,n^{-c_0}]} \left| \xi \cdot \int_0^t V(X_s) dW_s \right|
&=& \sup_{t \in [0,n^{-c_0}]} \left| B \left( \int_0^t \xi^\star V V^\star (x_s) \xi ds \right) \right| 
\\
&\geq&  \sup_{t \in [0,\ep_1 n^{-c_0}]} |B(u)|
\eeas
because of the uniform ellipticity. 
It follows from Lemma A of Watanebe (1984) that for all $x_0 \in \calx_0$, 
\beas
P \left[ \left.  C_n^c \cap D_n^c  \ \right| X_0 =x_0 \right]
&\leq&
P \left[ \left. \sup_{t \in [0,\ep_1 n^{-c_0}]} |B(u)| \leq C_2 n^{-c_0/2} n^{-\beta_1/2}  \ \right| X_0 =x_0 \right]
\\
&\leq&
c_3 \exp \left( -c_4 \frac{\ep_1 n^{-c_0}}{(C_2 n^{-c_0/2-\beta_1/2})^2 } \right)
\eeas
for some positive constants $c_3$ and $c_4$.
Hence, for any $\alpha >0$, there exist positive constants $\beta_0$, $c_0$ ($\beta_0 < \frac{c_0}{2} < \alpha$), $\gamma_2$ and $\beta_1$  
such that 
\bea
\sup_{x_0 \in \calx_0} P \left[ \left.  C_n^c \cap D_n^c  \ \right| X_0 =x_0 \right]
\leq
\gamma_2^{-1} \exp \left( -\gamma_2 n^{\beta_1}  \right) \label{u-p2-3}
\eea
for all $n \in \mathbb{N}$.
By (\ref{u-p2-1}), (\ref{u-p2-2}) and (\ref{u-p2-3}), for every $\alpha>0$, 
there exists positive constants $c_0$, $c_1$ and $c_2$ 
such that
\beas
& &
\inf_{x_0 \in \calx_0} P[ \ X_t -x_0 \in \mathbb{L}  \mbox{ and } |\xi \cdot(X_t-x_0)| = n^{-\alpha}
\mbox{ for some } t \in [0,n^{-c_0}] \ | \ X_0=x_0 ]
\\
&\geq& 1- \gamma_2^{-1} \exp \left( -\gamma_2 n^{-\beta_1}  \right) 
-2 \gamma_1^{-1} \exp (-\gamma_1 n^{c_0 -2\beta_0})
\\
&\geq& 1- c_1^{-1} \exp \left( -c_1 n^{-c_2}  \right), 
\eeas
which completes the proof of $(\ref{u-n1})$.

Next, we will show $[C]$. 
Let $\tau = \left\{ t>0 \ ; \ \sup_{s \in [0, n^{-d_0/l}]} |X_{t+s}-X_s| \geq n^{-\beta_0} \right\}$,
\beas
E_n = \left\{ \sup_{t \in [0,n^{-1}] \atop s \in [0, n^{-d_0/l}]} |X_{t+s}-X_s| \geq n^{-\beta_0} \right\}, \quad
F_n = \left\{ \sup_{s \in [0, n^{-d_0/l}]} |X_s-X_0| \geq n^{-\beta_0/l} \right\}.
\eeas
In the same way as the proof of (\ref{u-p2-2}), 
for each $d_0>0$, there exist positive constants $\beta_0$ ($\beta_0 < \frac{d_0}{2}$) and $\gamma_1$  
such that
\bea
\sup_{x_0 \in \calx_0} P \left[ \left. F_n \ \right| X_0 =x_0 \right] \leq \gamma_1^{-1} \exp (-\gamma_1 n^{(d_0-2\beta_0)/l}) \label{u-p2-4}
\eea
for all $n \in \mathbb{N}$ and $l \geq 1$.

Fix $d_0$.
Set
$\alpha(t,s) = C(1+|X_s|)t + \sup_{v \in [0,t]} \left| \int_s^{v+s} V(X_u) dW_u \right|$ and
$m(t,s) = |X_{t+s} -X_s|$.
Since Grownwall's lemma implies that 
$$
m(t,s) \leq \alpha(t,s) + C \int_0^t \alpha(u,s) e^{C(t-u)}du,
$$
one has that on $\{ \tau \leq n^{-1} \} \cap F_n^c$,
\beas
\sup_{t \in [0,\tau] \atop s \in [0, n^{-d_0/l}]} m(t,s)
&\leq& C_1 \sup_{t \in [0,\tau] \atop s \in [0, n^{-d_0/l}]} \alpha(t,s)
\\
&\leq& C_1(1+|X_0|+n^{-\beta_0/l})\tau 
+ \sum_{i=1}^{\sf d} \sup_{v \in [0,\tau] \atop s \in [0, n^{-d_0/l}]}  \left| \int_s^{v+s} V^{i \cdot} (X_u) dW_u \right|. 
\eeas
Noting that on $\{ \tau \leq n^{-1} \} \cap F_n^c$,
\beas
\sup_{v \in [0,\tau] \atop s \in [0, n^{-d_0/l}]} \int_s^{v+s} [V V^\star]^{i i}(X_u) du
\leq C_2(1+|X_0|^2 +n^{-2 \beta_0/l}+n^{-2 \beta_0}) n^{-1}, 
\eeas
one has that for sufficiently large $n$ and for $\beta_0 < \frac{1}{2}$,
\beas
& & P \left[ \left. E_n \cap F_n^c \ \right| X_0 =x_0 \right] 
=
P \left[ \left. \{ \tau \leq n^{-1} \} \cap F_n^c \ \right| X_0 =x_0 \right] 
\\
&\leq& 
\sum_{i=1}^{\sf d} P \left[ \left.  
\{ \tau \leq n^{-1} \} \cap F_n^c \cap \left\{
\frac{n^{-\beta_0} -C_2 n^{-1}}{\sf d} 
\leq C_3 \sup_{v \in [0,\tau] \atop s \in [0, n^{-d_0/l}]}  \left| \int_s^{v+s} V^{i \cdot} (x_u) dW_u \right| \right\}
\ \right| X_0 =x_0 \right] 
\\
&\leq&
{\sf d} \times P \left[ \left.  \frac{n^{-\beta_0}}{2} \leq \sup_{v \in [0, C_3n^{-1}]} |B(u)| \ \right| X_0 =x_0 \right] 
\leq 2 {\sf d} \exp \left( - \frac{(n^{-\beta_0}/2)^2}{2 C_3n^{-1}} \right)
\eeas
for all $x_0 \in \calx_0$.
Thus, for each $d_0>0$, there exist positive constants $\beta_0<\frac{1}{2}$ and $c_1$ such that
\bea
\sup_{x_0 \in \calx_0} P \left[ \left. E_n \cap F_n^c \ \right| X_0 =x_0 \right] 
\leq c_1^{-1} \exp (-c_1 n^{1-2\beta_0}) \label{u-p2-5}
\eea
for all $n \in \mathbb{N}$ and $l \geq 1$. It follows from (\ref{u-p2-4}) and (\ref{u-p2-5}) that
there exist positive constants $d_0$, $\beta_0$ ($\beta_0 < \frac{d_0}{2} < \frac{1}{2}$) and $d_1$ and $\beta_1$ 
such that
\beas
\sup_{x_0 \in \calx_0} P \left[ \left. E_n \ \right| X_0 =x_0 \right]
&\leq& 
\sup_{x_0 \in \calx_0} P \left[ \left. F_n \ \right| X_0 =x_0 \right] 
+ \sup_{x_0 \in \calx_0} P \left[ \left. E_n \cap F_n^c \ \right| X_0 =x_0 \right] 
\\
&\leq&
\gamma_1^{-1} \exp (-\gamma_1 n^{(d_0-2\beta_0)/l})
+ c_1^{-1} \exp (-c_1 n^{1-2\beta_0})
\\
&\leq& d_1^{-1} \exp (-d_1 n^{\beta_1/l})
\eeas
for all $n \in \mathbb{N}$ and $l \geq 1$. This completes the proof.

\qed\\
}}
\end{en-text}

\begin{en-text}
We will make a sketch of the proof. 
We need to change $x_0$, so it is not necessarily $X_0$. 
$x_0$ will become the initial value $X_{(m-1)n^{-c_0}}$ 
of the $m$-th trial by 
$(X_t)_{t \in [(m-1)n^{-c_0},mn^{-c_0}]}$. 
We take advantage of the projected increment 
$\xi_{\ell,k}\cdot(X_{mn^{-c_0}}-X_{(m-1)n^{-c_0}})$ 
only when $X$ exits from the ball 
$B(X_{(m-1)n^{-c_0}},n^{-\alpha})$ through 
$B(X_{(m-1)n^{-c_0}},n^{-\alpha})\cap C_{\xi_{\ell,k}}$. 
The probability of such an event is positive each time and 
by the nondegeneracy of $X$, 
for $c_2\in(0,c_0)$, 
it occurs at least one time in consecutive $n^{c_0-c_2}$ trials 
up to time $n^{-c_2}$ almost surely except for 
exponentially small probability. 
\qed\\
\end{en-text}


\begin{en-text}

\begin{example}\rm 
$\theta^*=0$, 
\beas 
S(x,\theta) 
&=& 
(1+x^2)^\theta+(1+2x^2)^{-\frac{\theta}{2}}.
\eeas
\end{example}

\begin{example}\rm 
$2$-dimensional $\theta^*=0$, $2$-dimensional $x^*=0$. 
\beas 
S(x,\theta) 
&=& 
(1-|\theta|^2)\sin x^{(1)}
+\theta^{(1)}\theta^{(2)}\sin x^{(2)}. 
\eeas
\end{example}

\end{en-text}

\begin{en-text}
{\bf kokomade} 
Let 
\beas 
\chi_J(x,\theta,\xi) 
&=& 
c_J(x,\theta)[\xi^{\otimes J}]. 
\eeas
Let $\bbL$ be a linear subspace of $\bbR^{\sf d}$ and 
$S^{{\sf d}-1}=\{\xi\in\bbR^{\sf d},\>|\xi|=1\}$. 
Let $\bbS=S^{{\sf d}-1}\cap\bbL$. 
\end{en-text}
\begin{en-text}
Suppose that $c_j$ $(j=0,...,J-1)$ are continuous on $\Theta$, 
and $c_{J+1}$ is continuous on $\hat{\calx}\times\Theta$. 
bounded 
in that 
\beas 
\sup_{(x,\theta)\in(\hat{\calx}\times K),
\xi\in \bbS}|c_J(x,\theta)[\xi^{\otimes J}]|<\infty. 
\eeas
Let 
\beas 
\chi_j(\theta,\xi) 
&=& 
c_j(\theta)[\xi^{\otimes j}]. 
\eeas

I use $S$ for $\Xi$ for ease of writing here. 
\end{en-text}

\begin{en-text}

\vspace{3mm}
{\bf koko}
We will discuss the expansion of the function 
that appears in Condition [H2]. 
For $(x_0,\theta^*)\in\calx\times\Theta$, let 
\beas 
U(x,\theta)
&=&
S(x,\theta^*)S(x,\theta)^{-2}S(x,\theta^*)
\eeas
and 
\beas 
D(x,\theta)
&=&
U(x_0,\theta)^{-\half}
\big(U(x,\theta)-U(x_0,\theta)\big)U(x_0,\theta)^{-\half}. 
\eeas
Then 
\beas 
\det U(x,\theta)
&=&
\det U(x_0,\theta)
\det \big( I_d + D(x,\theta) \big). 
\eeas
Therefore
\beas 
\log\det\bigg(S(x,\theta)^{-1}S(x,\theta^*)\bigg) 
&=&
\half\log\det U(x,\theta)
\\&=&
\half\log \det U(x_0,\theta)+\half\log \det \big( I_d + D(x,\theta) \big)
\eeas

We denote by $B(x_0,\rho)$ an open ball centered at $x_0$ 
with radius $\rho$ included in $\calx$. 
Let $\zeta\sim N_d(0,I_d)$. 
For sufficiently small $\rho$, 
$\sup_{(x,\theta)\in B(x_0,\rho)\times\Theta}|D(x,\theta)|<1$, 
therefore 
we have the asymptotic expansion 
\beas 
&&
-\half\log\det\big( I_d + D(x,\theta) \big)
\\&=&
\log\bigg\{
(2\pi)^{-\frac{d}{2}} \int\exp(-\half ( I_d +  D(x,\theta)) 
[z^{\otimes2}]\big)dz
\bigg\}
\\&=&
\log\bigg\{
\sum_{k=0}^\infty 
\frac{1}{k!}
\l(-\half\r)^k
E\bigg[\big(D(x,\theta)[\zeta^{\otimes2}]\big)^k\bigg]
\bigg\}
\\&=&
\log\bigg\{
\sum_{k=0}^\infty 
\frac{1}{k!}
\l(-\half\r)^k
E\bigg[\bigg(
\sum_{m=1}^\infty\frac{1}{m!}\partial_x^mD(x_0,\theta)
[(x-x_0)^{\otimes m},\zeta^{\otimes2}]\bigg)^k\bigg]
\bigg\}
\\&=&
\sum_{\lambda=1}^\infty C_\lambda(x_0,\theta)
[(x-x_0)^{\otimes\lambda}]
\eeas
as $x\to x_0$; in other words, 
the difference between the left-hand side and 
the partial sum of the terms with $(x-x_0)^{\otimes j}$ $(j\leq k)$ 
on the right-hand side 
is of $O(|x-x_0|^{k+1})$ as $x\to x_0$ 
for each $k\in\bbN$. 
In particular, 
\beas 
C_1(x_0,\theta)[\xi]
&=&
-\half E\big[\partial_x D(x_0,\theta)[\xi,\zeta^{\otimes2}]\big]
\\&=&
-\half E\big[\big(U(x_0,\theta)^{-\half}
\partial_x U(x_0,\theta)[\xi]U(x_0,\theta)^{-\half}\big)
[\zeta^{\otimes2}]\big]
\\&=&
-\half {\rm Tr}\big(U(x_0,\theta)^{-1}
\partial_x U(x_0,\theta)[\xi]\big). 
\eeas
Hence 
\beas 
C_1(x_0,\theta)[\xi]
&=&
- {\rm Tr}\bigg(S(x_0,\theta^*)^{-1}\big(\partial_xS(x_0,\theta^*)
-\partial_xS(x_0,\theta)\big)[\xi]\bigg)
\eeas
for $\theta$ satisfying $S(x_0,\theta)=S(x_0,\theta^*)$, and then 
\beas 
\partial_xQ(x_0,\theta,\theta^*)&=&0. 
\eeas
This equality can be obtained also from the nonnegativity of 
$Q(x,\theta,\theta^*)$ whenever $S(x_0,\theta)=S(x_0,\theta^*)$. 

\beas 
C_2(x_0,\theta)
[\xi^{\otimes 2}]
&=& \frac{1}{4} \mbox{Tr} \left( U^{-1}(x_0,\theta) (\partial_x U)(x_0,\theta)[\xi] U^{-1}(x_0,\theta) (\partial_x U)(x_0,\theta)[\xi]
- U^{-1}(x_0,\theta) (\partial_x^2 U)(x_0,\theta)[\xi^{\otimes 2}] \right)
\eeas

\begin{description}
\item[[B\!\!]] 
$\Xi(X_{t^*},\theta^*)$ is positive-definite a.s., 
where 
\beas 
\Xi(x,\theta)=
\bigg(\mbox{Tr}((\partial_{\theta^{(i_1)}} S)S^{-1}
(\partial_{\theta^{(i_2)}} S)S^{-1})(x,\theta)
\bigg)_{i_1,i_2=1}^{d_0}, \sskip
\theta=(\theta^{(i)}). 
\eeas 
\end{description}

\begin{remark}\rm 
Since $\Xi(x,\theta^*)$ is the Hessian matrix of 
the nonnegative function 
\beas 
Q(x,\theta^*,\theta):=
\mbox{Tr}\bigg(S(x,\theta^*)^{-1}S(x,\theta)-I_d\bigg)
-\log\det\bigg(S(x,\theta^*)^{-1}S(x,\theta)\bigg)
\eeas
of $\theta$ at $\theta^*$, 
$\Xi(x,\theta^*)$ is nonnegative-definite. 
\end{remark}

\bibliographystyle{plain}
\bibliography{naka-db}

\end{en-text}

\begin{en-text}

We denote $B(R) = \{ u \in \bbR^{{\sf p}} \ ; \ |u| \leq R \}$.

\begin{lemma} \label{lem3} 
Assume $[H1]$. Then,
for every $q>p$ and $R>0$, there exists  $C_0 >0$ such that
$$
\sup_{n \in {\bf N}} E_{\theta^*} \left[ \left| \log {\mathbb Z}_n (u)  \right|^q \right] \leq C_0 |u|^q
$$
for all $u \in B(R)$.
\end{lemma}

\begin{remark}\rm
By using Lemma 2 in Yoshida (2005), Lemma \ref{lem3} implies that for every $R>0$,
$$
\sup_{n \in {\bf N}} E_{\theta^*} \left[ \left( \int_{|u| \leq R}  {\mathbb Z}_n (u)  du  \right)^{-1} \right] <\infty.
$$
Consequently, one has that
$$
\sup_{n \in {\bf N}} E_{\theta^*} \left[ \left( \int_{{\mathbb U}_n}  {\mathbb Z}_n (u)  du  \right)^{-1} \right] <\infty.
$$
\end{remark}

\end{en-text}




\section{Examples and simulation results}

As an example, 
we consider the one-dimensional diffusion process
\begin{eqnarray}
dX_t &=&  X_t dt + \exp \{ \theta \sin^2 X_t \} dw_t, \quad t \in [0,1], \quad X_0 =0, \label{ex-sde0}
\end{eqnarray}
where $\theta \in (-\pi,\pi)$.

For the simulations, in order to get the maximum likelihood type estimator,
we used the MATHEMATICA 6.0, 
concretely,
``{\bf FindMinimum}" with an initial value.
We examine the asymptotic behaviour of the estimators,
which are  
the maximum likelihood type estimator $\hat{\theta}_{n}$
obtained by using "FindMinimum"
with the initial value
$\theta_0 =0.5$
the Bayes type estimator $\tilde{\theta}_{n}$ 
with respect to the uniform prior $\pi(\theta)=1/(2 \pi)$
and the maximum likelihood type estimator $\hat{\theta}_{n}^{(1)}$
obtained by using "FindMinimum"
with initial value being the Bayes type estimator $\tilde{\theta}_{n}$,
through the simulations, which were done 
for each $h_n=1/50, 1/250, 1/500$. 
For the true model (\ref{ex-sde0}) with $\theta^* =1$,
$10000$ independent sample paths are generated by
the Milstein scheme,
and the means and the standard deviations of the estimators are computed and shown in Table 1 below. 


In Table 1, even if $h_n =1/50$, the Bayes estimator $\tilde{\theta}_{n}$ has good performance, 
but both the maximum likelihood type estimators $\hat{\theta}_{n}$ and  $\hat{\theta}_{n}^{(1)}$ have biases.
In case that $h_n=1/250$, all three estimators are unbiased and they have good behaviors.
In this example, it is better to use the Bayes estimator than the maximum likelihood type estimators.

\begin{en-text}
\begin{table}[h]
\caption[Tab_TAN01]{
The mean and standard deviation (s.d.) 
of the  maximum likelihood type estimators with the initial value $\theta_0 =0, 0.5, 1$ for $10000$ independent simulated sample paths
with $\theta^* = 1$. 
}
  \begin{center}
    \begin{tabular}{|c|c|c|c|c|c|c|}
      \hline
          & \multicolumn{2}{c|}{$\hat{\theta}_{n}$ with $\theta_0 =0$} 
          &   \multicolumn{2}{c|}{$\hat{\theta}_{n}$ with $\theta_0 =0.5$} 
          &   \multicolumn{2}{c|}{$\hat{\theta}_{n}$ with $\theta_0 =1$} \\
      \cline{2-7}
        $h_n$ &  mean & s.d. & mean & s.d. & mean & s.d.  \\
      \hline
  1/50 & 0.89850 & 0.55160 & 0.89850 & 0.55160 & 0.89850 & 0.55160 \\
\hline
  1/250 & 0.97816 & 0.23723 & 0.97816 & 0.23723 & 0.97816 & 0.23723 \\
\hline
  1/500 & 0.99145 & 0.16041 & 0.99145 & 0.16041 & 0.99145 & 0.16041 \\
\hline
    \end{tabular}
  \end{center}
\end{table}
\end{en-text}

\begin{table}[h]
\caption[Tab_TAN01]{
The mean and standard deviation (s.d.) 
of the three kinds of estimators for $10000$ independent simulated sample paths
with $\theta^* = 1$. 
}
  \begin{center}
    \begin{tabular}{|c|c|c|c|c|c|c|}
      \hline
          & \multicolumn{2}{c|}{$\hat{\theta}_{n}$ with $\theta_0 =0.5$} 
          &   \multicolumn{2}{c|}{$\tilde{\theta}_{n}$}  
          &   \multicolumn{2}{c|}{$\hat{\theta}_{n}^{(1)}$ with $\theta_0 =\tilde{\theta}_{n}$} \\
      \cline{2-7}
        $h_n$ &  mean & s.d. & mean & s.d. & mean & s.d.  \\
      \hline
  1/50 & 0.89850 & 0.55160 & 0.96473 & 0.48914 & 0.89850 & 0.55160 \\
\hline
  1/250 & 0.97816 & 0.23723 & 0.99392 & 0.23112 & 0.97816 & 0.23723 \\
\hline
  1/500 & 0.99145 & 0.16041 & 0.99969 & 0.15874 & 0.99145 & 0.16041 \\
\hline
    \end{tabular}
  \end{center}
\end{table}

As another example, 
we consider the one-dimensional diffusion process
\begin{eqnarray}
dX_t &=&  X_t dt + \exp \{ \sin \theta \sin X_t - \theta^2 \sin^2 X_t \} dw_t, \quad t \in [0,1], \quad X_0 = 0,  \label{ex-sde1} 
\end{eqnarray}
where $\theta \in (-\pi,\pi)$.

For the true model (\ref{ex-sde1}) with $\theta^* =0$, simulations were done in the same way as the previous example.
The means and the standard deviations of 
the maximum likelihood type estimator $\hat{\theta}_{n}$ with the initial value 
$\theta_0=0.5$,
the Bayes type estimator $\tilde{\theta}_{n}$  
with respect to the uniform prior $\pi(\theta)=1/(2 \pi)$
and the maximum likelihood type estimator $\hat{\theta}_{n}^{(1)}$ with the initial value $\theta_0= \tilde{\theta}_{n}$
are computed and shown in Table 2 below. 



In Table 2, the maximum likelihood type estimator $\hat{\theta}_{n}$ with $\theta_0 =0.5$ has a bias in all cases,
while the Bayes type estimator $\tilde{\theta}_{n}$ and 
the maximum likelihood estimator $\hat{\theta}_{n}^{(1)}$ with $\theta_0 =\tilde{\theta}_{n}$ have good behaviors in all cases.
Furthermore, we see that the standard deviation of the Bayes estimator $\tilde{\theta}_{n}$
is smaller than the one of $\hat{\theta}_{n}^{(1)}$ in all cases.

\begin{en-text}
>From the theoretical point of view,
both the maximum likelihood type estimator and
the Bayes type estimator with any prior have the same asymptotic properties.
On the other hand, in case that $n$ is not so large,
we have to investigate which estimator should be used 
since both estimators does not always have the same performances.
\end{en-text}

\begin{en-text}
\begin{table}[h]
\caption[Tab_TAN02]{
The mean and standard deviation (s.d.) 
of the estimators for $10000$ independent simulated sample paths
with $\theta^* = 0$. 
}
  \begin{center}
    \begin{tabular}{|c|c|c|c|c|c|c|}
      \hline
          & \multicolumn{2}{c|}{$\hat{\theta}_{n}$ with $\theta_0 =0$} 
          &   \multicolumn{2}{c|}{$\hat{\theta}_{n}$ with $\theta_0 =0.5$} 
          &   \multicolumn{2}{c|}{$\hat{\theta}_{n}$ with $\theta_0 =1$} \\
      \cline{2-7}
        $h_n$ &  mean & s.d. & mean & s.d. & mean & s.d.  \\
      \hline
  1/50  & 0.00285 & 0.29301 & 0.07702 & 0.39184 & 0.34007 & 0.57153 \\
\hline
  1/250 & 0.00095 & 0.11062 & 0.05046 & 0.23677 & 0.35170 & 0.50177 \\
\hline
  1/500 & 0.00042 & 0.07083 & 0.04230 & 0.20471 & 0.35517 & 0.49420 \\
\hline
    \end{tabular}
  \end{center}
\end{table}
\end{en-text}

\begin{table}[h]
\caption[Tab_TAN02]{
The mean and standard deviation (s.d.) 
of the estimators for $10000$ independent simulated sample paths
with $\theta^* = 0$. 
}
  \begin{center}
    \begin{tabular}{|c|c|c|c|c|c|c|}
      \hline
          & \multicolumn{2}{c|}{$\hat{\theta}_{n}$ with $\theta_0 =0.5$} 
          &   \multicolumn{2}{c|}{$\tilde{\theta}_{n}$} 
          &   \multicolumn{2}{c|}{$\hat{\theta}_{n}^{(1)}$ with $\theta_0 =\tilde{\theta}_{n}$} \\
      \cline{2-7}
        $h_n$ &  mean & s.d. & mean & s.d. & mean & s.d.  \\
      \hline
  1/50 & 0.07702 & 0.39184 & 0.00363 & 0.43532 & 0.00757 & 0.52096 \\
\hline
  1/250 & 0.05046 & 0.23677 & 0.00317 & 0.26283 & 0.00502 & 0.28079 \\
\hline
  1/500 & 0.04230 & 0.20471 & 0.00071 & 0.16614 & -0.00035 & 0.17814 \\
\hline
    \end{tabular}
  \end{center}
\end{table}


\section{A geometric criterion}\label{26721210-1}
{\cred 
Apart from analytic criteria by derivatives, 
we shall consider the following condition in the spirits of Lemma \ref{210102-2} and Remarks \ref {241205-1} and \ref{241205-2}. 
\begin{description}
\item[[A3$'$\!\!]] 
$\mbox{supp}\call\{X_\tau\}$ is compact, 
there exists a function $f:U\times\Theta\to\bbR$ for some open neighborhood $U$ of $\mbox{supp}\call\{X_\tau\}$
and the following conditions are satisfied. 
\begin{description}
\item[(i)] For some $\varrho\in(0,\infty)$ , 
$
Q(x,\theta,\theta^*)|\theta-\theta^*|^{-2}
\geq |f(x,\theta)|^\varrho
$ 
for all $(x,\theta)\in U\times(\Theta\setminus\{\theta^*\})$. 
\item[(ii)] For each $x_0\in U$, there exist a neighborhood $V$ in $U$ of $x_0$ 
and a covering $\{\Theta_k\}_{k=1,...,\bar{k}}$ of $\Theta$ 
such that for each $k=1,...,\bar{k}$, there exist 
$\xi_0\in\bbS$, 
$J\in\bbN$, some positive numbers $M,c,\ep_0$, $K_j$ ($j=1,...,J$) and 
some functions $\Psi_j:P^\perp_{\xi_0}V\times\Theta_k\to\bbR$ such that 
\begin{description}
\item[(a)] the function $P^\perp_{\xi_0}V\ni y\mapsto \Psi(y,\theta)\in\bbR$ is 
$M$-Lipschitz continuous for all $\theta\in\Theta_k$, 
\item[(b)] for $(x,\theta)\in V\times\Theta_k$, 
\beas 
|f(x,\theta)| &\geq&
c\prod_{j=1}^J\big(|\xi_0\cdot x-\Psi_j(P^\perp_{\xi_0}x,\theta)|\wedge\ep_0\big)^{K_j}.
\eeas
\end{description}
\end{description}
\end{description}

In $[A3']$, $\bar{k}$ may depend on $x_0$. 
Note that 
$
\{x\in V;\>f(x,\theta)=0\} \subset
\bigcup_{j=1}^J\{x\in V;\>\xi_0\cdot x=\Psi_j(P^\perp_{\xi_0}x,\theta)\}
$ 
under $[A3']$(ii), that is, the graph of the functions $\Psi_j$ covers locally the null set of $f$.  
\begin{theorem}
Suppose that Conditions $[A1]$, $[A2]$ and $[A3']$ are satisfied. 
Then 
the same results as Theorems \ref{230614-1} and \ref{230614-2} hold true. 
\end{theorem}
\proof 
We consider an open ball $B(x_0,\ep_{x_0})\subset V$ for each $x_0\in U$ and the covering 
$\{B(x_0,\ep_{x_0}/2)\}_{x_0\in U}$ 
of $\mbox{supp}\call\{X_\tau\}$.  
By compactness, we obtain a finite number of balls $V$, and as a result, we have an open neiborghhood 
$\calx_0$ of $\mbox{supp}\call\{X_0\}$ and we  
may assume that 
for some $\ep'>0$, every $B(x,\ep')$ ($x\in\calx_0$) can find a $V\supset B(x,\ep')$ among them. 
Call these $V$'s $\calx_\ell$ ($\ell=1,...,\bar{\ell}$). Each $\calx_\ell$ has a partition 
$\{\Theta_{\ell,k}\}_{k=1,...,\bar{k}_\ell}$. 

Let $\alpha_0>\alpha_1>\cdots>\alpha_J>0$, and $\ep>0$. 
We consider a sufficiently large $L$ and sufficiently large $n$'s. 
Fix $V=\calx_\ell$ and a $\Theta_{\ell,k}$, for which we have $\xi_0=\xi_{\ell,k}$ and $\Psi_j$ 
depending on $(\ell,k)$. 
For $x^*\in \calx_\ell$, let 
$\varsigma_i=x^*+n^{-\alpha_i}D(\xi_0,\ep)$ for $i=0,1,...,J$, and 
denote $n^{-L}$-neighborhood of $\varsigma_i$ by $\varsigma_i^n$. 
Moreover denote by $\calg_j$ 
the graph of $(x,\Psi_j(\cdot,\theta))$ in $S_{x^*}(\xi_0,\ep)$. 

We claim that there is no $\calg_j$ that intersects with two different $\varsigma_i$'s. 
Indeed, if any $\calg_j$ intersected with $\varsigma_{i_1}$ and $\varsigma_{i_2}$ for 
$i_1<i_2$, there are $x_k\in x^*+n^{-\alpha_{i_k}}D(\xi_0,\ep)$ ($k=1,2$) such that 
\bea\label{2672-1209-11}
2n^{-L}+|\Psi_j(P^\perp_{\xi_0}x_1,\theta)-\Psi_j(P^\perp_{\xi_0}x_2,\theta)|>\half n^{-\alpha_{i_2}}.
\eea
On the other hand, by the $M$-Lipschitz continuity of $\Psi_j(\cdot,\theta)$, we have 
\beas 
|\Psi_j(P^\perp_{\xi_0}x_1,\theta)- \Psi_j(P^\perp_{\xi_0}x_2,\theta)| 
&\leq& 
2M\sqrt{\ep}n^{-\alpha_{i_2}},
\eeas
which contradicts to (\ref{2672-1209-11}) if we make $\ep$ sufficiently small, and proved the claim. 
Therefore, there is at least one $\varsigma_i^n$ that does not intersect with any $\calg_j$. 
Thus
\beas 
\inf_{(x^*,\theta)\in \calx_\ell\times\Theta_{\ell,k}}
\max_{i=0,1,...,J}\min_{j=1,...,J}\inf_{\xi\in D(\xi_{\ell,k},\ep)}
\big|\xi_0\cdot (x^*+n^{-\alpha_i}\xi)-\Psi_j(P^\perp_{\xi_0}(x^*+n^{-\alpha_i}\xi),\theta)\big|
&\geq&
n^{-L}
\eeas
for large $n$ for every $\ell=1,..,\bar{\ell}$ and $k=1,..,\bar{k}_\ell$. 
Consequently, taking large $L'$, we obtain 
\bea\label{26721209-12}
\min_{\ell=1,...,\bar{\ell},\>k=1,...,\bar{k}_\ell}
\inf_{(x^*,\theta)\in \calx_\ell\times\Theta_{\ell,k}}
\max_{i=0,1,...,J}\inf_{\xi\in D(\xi_{\ell,k},\ep)}
|f(x^*+n^{-\alpha_i}\xi,\theta)|
&\geq&
cn^{-L'}
\eea
for large $n$. 
This is Lemma \ref{210102-2} with $\tilde{D}(\xi_{\ell,k},\ep_{\ell,i})$ replaced by 
$D(\xi_{\ell,k},\ep)$. Due to Remarks \ref {241205-1} and \ref{241205-2}, 
we can prove the theorem in the same way as Theorems \ref{230614-1} and \ref{230614-2}. 
\qed 
}


\section{Proof of Theorems \ref{thm1}, \ref{thm3} and \ref{thm2}}\label{230605-1}
{\colorb{
For the limit of $\bbZ_n$ given in (\ref{230808-1}), we define 
\beas 
{\mathbb Z}(u) = \exp \left( \Gamma(\theta^*)^{1/2} \zeta [u] -\frac{1}{2} \Gamma(\theta^*)[u,u] \right). 
\eeas
Then the standardized quasi Bayesian estimator  $\tilde{u}_n = \sqrt{n} ( \tilde{\theta}_n -\theta^*)$ is written by
\begin{eqnarray*} \label{Bayes2}
\tilde{u}_n = \left( \int_{{\mathbb U}_n} {\mathbb Z}_n(u) \pi(\theta^* + (1/\sqrt{n})u) du \right)^{-1}
\int_{{\mathbb U}_n} u {\mathbb Z}_n(u) \pi(\theta^* + (1/\sqrt{n})u) du. 
\end{eqnarray*}
As we will prove later, the limit of $\tilde{u}_n$ should be 
\begin{eqnarray*} \label{Bayes4}
\tilde{u} = \left( \int_{\bbR^{{\sf p}}} {\mathbb Z}(u)  du \right)^{-1}
\int_{\bbR^{{\sf p}}} u {\mathbb Z}(u) du \  = \Gamma(\theta^*)^{-1/2}\zeta .
\end{eqnarray*}
However, even existence of the integrals requires more rigorous treatment.  

In order to prove Theorems \ref{thm3} and \ref{thm2}, we will first prepare several lemmas. 
The results for the quasi maximum likelihood estimator are proved at the same time with common machinery. 
}}


\begin{lemma} \label{lem4-1} 
Let $f \in C_\uparrow^{1,1}({\bbR}^{\sf d} \times \Theta; {\bbR})$.
Assume $[H1]$. Then, for every $q > 0$,
$$
\sup_{n \in {\mathbb N}} E \left[
 \left( \sup_{\theta \in \Theta}  \sqrt{n} \left| \frac{1}{n} \sum_{k=1}^n f(\xt, \theta) 
- \frac{1}{T} \int_0^T f(X_t,\theta) dt \right| \right)^q \right] < \infty.
$$
\end{lemma}

\proof Let 
$$F_n(\theta) = \sumkn \int_{t_{k-1}}^{t_k} \left\{ f(\xt,\theta) - f(X_t,\theta) \right\} dt.
$$
An easy estimate together with $[H1]$-$(i)$ implies that for every $q >{\cred \sf p}$,
\begin{eqnarray*}
\sup_{\theta \in \Theta}  \left| \left|  \sqrt{n} F_n(\theta) \right| \right|_q^q 
&\leq& \frac{n^{3q/2}}{n^q} \sumkn \int_{t_{k-1}}^{t_k} \sup_{\theta \in \Theta} || f(\xt,\theta) - f(X_t,\theta) ||_q^q dt
< C.
\end{eqnarray*}
Therefore, 
$$\sup_{n \in {\mathbb N}} \sup_{\theta \in \Theta}  \left| \left|  \sqrt{n} F_n(\theta) \right| \right|_q < \infty.$$
Moreover, in {\cred a} similar way,
$$\sup_{n \in {\mathbb N}} \sup_{\theta \in \Theta}  \left| \left|  \sqrt{n}  \partial_\theta F_n(\theta) \right| \right|_q < \infty.$$
It follows from the Sobolev inequality that
\begin{eqnarray*}
E \left[ \sup_{\theta \in \Theta} | \sqrt{n} F_n(\theta) |^q \right] 
&\leq& E \left[ C \int_{\Theta} \left\{  | \sqrt{n} F_n(\theta) |^q +  | \sqrt{n} \partial_\theta F_n(\theta) |^q \right\} d\theta
\right] \\
&\leq& C_{\Theta} \left\{ 
\sup_{\theta \in \Theta} E \left[  | \sqrt{n} F_n(\theta) |^q \right] 
+ \sup_{\theta \in \Theta} E \left[  | \sqrt{n} \partial_\theta F_n(\theta) |^q \right] 
\right\},
\end{eqnarray*}
where $q > p$. Thus, one has that for every $ q>0$,
$$\sup_{n \in {\mathbb N}} \left| \left|  \sup_{\theta \in \Theta} \left| \sqrt{n} F_n(\theta) \right| \right| \right|_q
< \infty.$$ 
This completes the proof.
\qed

\begin{lemma} \label{lem4-2}
Let $f \in C_\uparrow^{0,1}({\bbR}^{\sf d} \times \Theta; {\bbR}^{\sf m} \otimes {\bbR}^{\sf m})$.
Assume $[H1]$. Then, for every $q > 0$,
$$
\sup_{n \in {\mathbb N}} E \left[ \left( \sup_{\theta \in \Theta}
\left| \sqrt{n}  \sum_{k=1}^n  f(\xt, \theta) \left\{ [ (\Delta_k Y)^{\otimes 2}] 
- [ ( \sigma(\xt,\theta^*) \Delta_kw )^{\otimes 2} ] \right\}  \right| \right)^q \right] < \infty. 
$$
\end{lemma}

\proof Noting that 
\begin{eqnarray*}
\Delta_k Y &=&
\int_{t_{k-1}}^{t_k} b_t dt + \int_{t_{k-1}}^{t_k} \left\{ \sigma(X_t,\theta^*) -\sigma(\xt,\theta^*) \right\}dw_t 
+ \sigma(\xt,\theta^*) \Delta_kw,
\end{eqnarray*}
\begin{eqnarray*}
\left| \left| \int_{t_{k-1}}^{t_k} b_t dt \right| \right|_q^q +
\left| \left| \int_{t_{k-1}}^{t_k} \left\{ \sigma(X_t,\theta^*) -\sigma(\xt,\theta^*) \right\}dw_t  \right| \right|_q^q &\leq& C n^{-q}, 
\end{eqnarray*}
and $\left| \left| \sigma(\xt,\theta^*) \Delta_kw  \right| \right|_q^q \leq C n^{-q/2}$ 
{\cred for every $q>{\sf p}\vee 2$,}
one has that 
\begin{en-text}
\begin{eqnarray*}
& & \sup_{\theta \in \Theta}  \left| \left| \sqrt{n}  \sum_{k=1}^n  f(\xt, \theta) \left\{ [ (\Delta_k Y)^{\otimes 2}] 
- [ ( \sigma(\xt,\theta^*) \Delta_kw )^{\otimes 2} ] \right\}    \right| \right|_q^q\\
&\leq& C\frac{n^{3q/2}}{n} \sumkn \left\{ \frac{1}{n^{2q}} + \frac{1}{n^{3q/2}} \right\} <C.
\end{eqnarray*}
Hence 
\end{en-text}
$$
\sup_{n \in {\mathbb N}} \sup_{\theta \in \Theta}  \left| \left| \sqrt{n}  \sum_{k=1}^n  f(\xt, \theta) \left\{ [ (\Delta_k Y)^{\otimes 2}] 
- [ ( \sigma(\xt,\theta^*) \Delta_kw )^{\otimes 2} ] \right\}    \right| \right|_q <\infty,
$$
{\cred with $ab-cd=(a-c)b+(b-d)c$. }
Similarly,
$$
\sup_{n \in {\mathbb N}} \sup_{\theta \in \Theta}  \left| \left| \sqrt{n}  \sum_{k=1}^n  \partial_\theta f(\xt, \theta) \left\{ [ (\Delta_k Y)^{\otimes 2}] 
- [ ( \sigma(\xt,\theta^*) \Delta_kw )^{\otimes 2} ] \right\}    \right| \right|_q <\infty.
$$
By using the Sobolev inequality, we obtain the desired inequality.
\qed

\vspace{0.5cm}
Note that for $u \in {\mathbb U}_n$,
$$
{\mathbb Z}_n(u) = \exp \left( \Delta_n[u] -\frac{1}{2} \Gamma(\theta^*)[u,u] + r_n(u) \right),  
$$
where
\begin{eqnarray*}
\Delta_n[u] &=& \frac{1}{\sqrt{n}} \partial_\theta {\mathbb H}_n(\theta^*)[u] \\
&=& -\frac{1}{2 \sqrt{n}} \sumkn \left\{ (\partial_\theta \log \det S(\xt, \theta^*) )[u] 
+ h^{-1} (\partial_\theta S^{-1})(\xt,\theta^*) [u, (\Delta_k Y)^{\otimes 2}] \right\}, 
\end{eqnarray*}
\begin{eqnarray*}
\Gamma_n(\theta)[u,u] &=& - \frac{1}{n} \partial_\theta^2 {\mathbb H}_n(\theta)[u,u] \\
&=& \frac{1}{2 n} \sumkn \left\{ (\partial_\theta^2 \log \det S(\xt, \theta) )[u^{\otimes 2}] 
+ h^{-1} (\partial_\theta^2 S^{-1})(\xt,\theta) )[u^{\otimes 2}, (\Delta_k Y)^{\otimes 2}] \right\}, \\
\Gamma(\theta^*)[u,u] &=& \frac{1}{2 T} \int_0^T \mbox{Tr} \left( (\partial_\theta S) S^{-1}
(\partial_\theta S) S^{-1}(X_t,\theta^*)[u^{\otimes 2}] \right) dt, \\
\end{eqnarray*}
\begin{eqnarray*}
r_n(u) &=& \int_0^1 (1-s) \left\{ \Gamma(\theta^*)[u,u] - \Gamma_n(\theta^* + s(1/\sqrt{n})u)[u,u] \right\}ds.
\end{eqnarray*}

\begin{en-text}
Let
\begin{eqnarray*}
{\mathbb Y}_n(\theta) &=& \frac{1}{n} \left\{ {\mathbb H}_n (\theta) - {\mathbb H}_n (\theta^*) \right\}, \\
{\mathbb Y}(\theta) &=& -\frac{1}{2 T}\int_0^T \left\{ \log \left( \frac{\det S(X_t,\theta)}{\det S(X_t,\theta^*)} \right)
+\mbox{Tr} \left( S^{-1}(X_t,\theta) S(X_t,\theta^*) -I_d \right) \right\} dt.
\end{eqnarray*}
\end{en-text}

\begin{lemma} \label{lem1}
Assume $[H1]$. Then, for every $q>0$,
\begin{description}
\item[(i)]  
$\displaystyle 
\sup_{n \in {\bf N}} E \left[ \left| \Delta_n \right|^q \right] <\infty.
$
\item[(ii)]  
$\displaystyle 
\sup_{n \in {\bf N}} E 
\left[ \left( \sup_{\theta \in \Theta} \sqrt{n} \left| {\mathbb Y}_n (\theta) - {\mathbb Y} (\theta) \right| \right)^q \right] <\infty.
$
\end{description}
\end{lemma}

\proof 
(i) Since 
$\Delta_n = M_n +R_n$,
where
\begin{eqnarray*}
M_n &=& -\frac{1}{2 \sqrt{n}} \sumkn \left\{ (\partial_\theta \log \det S(\xt, \theta^*) ) 
+ h^{-1} (\partial_\theta S^{-1}(\xt,\theta^*) )[(\sigma(\xt,\theta^*) \Delta_kw)^{\otimes 2}] \right\}, \\
R_n &=&  -\frac{\sqrt{n}}{2 T } \sumkn 
(\partial_\theta S^{-1}(\xt,\theta^*) )[ (\Delta_k Y)^{\otimes 2} - (\sigma(\xt,\theta^*) \Delta_kw)^{\otimes 2}],
\end{eqnarray*}
Lemma \ref{lem4-2} yields that for every $q>{\cred 1}$,
$
\sup_{n \in {\mathbb N}} || R_n ||_q 
< \infty.
$
Moreover, $\sqrt{n} M_n$ is the terminal value of a discrete-time martingale with respect to $({\cal F}_{t_k})_{k=0,1,\ldots,n}$
and it follows from the Burkholder inequality that
$
\sup_{n \in {\mathbb N}} || M_n ||_q < \infty.
$
Thus, one has that $\sup_{n \in {\mathbb N}} || \Delta_n ||_q < \infty$ for every $q>{\cred 1}$.

(ii)
Note that
$$
{\mathbb Y}_n(\theta) = {\mathbb Y}_n^\dagger(\theta) + M_n^\dagger(\theta) + R_n^\dagger(\theta),
$$
where
\begin{eqnarray*}
{\mathbb Y}_n^\dagger(\theta) &=& - \frac{1}{2 n} \sumkn \left\{ \log \left( \frac{\det S(\xt,\theta)}{\det S(\xt,\theta^*)} \right)
+ \mbox{Tr} \left( S^{-1}(\xt,\theta) S(\xt,\theta^*) -I_d \right) \right\}, \\
M_n^\dagger(\theta) &=& - \frac{1}{2 n} \sumkn ( S^{-1}(\xt,\theta) - S^{-1}(\xt,\theta^*) )
[h^{-1} (\sigma(\xt,\theta^*) \Delta_kw)^{\otimes 2} - S(\xt,\theta^*)], \\
R_n^\dagger(\theta) &=& - \frac{1}{2 n} \sumkn ( S^{-1}(\xt,\theta) - S^{-1}(\xt,\theta^*) )
[ h^{-1} (\Delta_k Y)^{\otimes 2} - h^{-1} (\sigma(\xt,\theta^*) \Delta_kw)^{\otimes 2}].
\end{eqnarray*}
By Lemma \ref{lem4-2}, for every $q>{\cred 1}$,
$$ \sup_{n \in {\mathbb N}}  || \sup_{\theta \in \Theta} |\sqrt{n} R_n^\dagger(\theta)| ||_q < \infty. $$
Burkholder's inequality yields that
$$
\sup_{n \in {\mathbb N}} \sup_{\theta \in \Theta} || \sqrt{n} M_n^\dagger(\theta) ||_q < \infty
\quad \mbox{and} \quad 
\sup_{n \in {\mathbb N}} \sup_{\theta \in \Theta} || \sqrt{n} \partial_\theta M_n^\dagger(\theta) ||_q < \infty
$$
for all $q>{\cred 1}$.
Hence, Sobolev's inequality implies that 
$$ \sup_{n \in {\mathbb N}}  || \sup_{\theta \in \Theta} |\sqrt{n} M_n^\dagger(\theta)| ||_q < \infty $$
{\cred for every $q>{\sf p}$.}
Furthermore, it follows from Lemma \ref{lem4-1} that
$$
\sup_{n \in {\mathbb N}}  || \sup_{\theta \in \Theta} \sqrt{n} | {\mathbb Y}_n^\dagger(\theta) - {\mathbb Y} (\theta)| ||_q < \infty
$$
for every $q>{\cred 1}$.
This completes the proof.
\qed

\vspace{0.5cm}
\begin{lemma} \label{lem2}
Assume $[H1]$. Then, for every $q>0$,
\begin{description}
\item[(i)] 
$\displaystyle 
\sup_{n \in {\bf N}} E 
\left[ \left( \sqrt{n} \left| \Gamma_n (\theta^*) - \Gamma (\theta^*) \right| \right)^q \right] <\infty.
$
\item[(ii)]  
$\displaystyle 
\sup_{n \in {\bf N}} E 
\left[ \left( \frac{1}{n} \sup_{\theta \in \Theta}  \left| \partial_\theta^3 {\mathbb H}_n (\theta) \right| \right)^q \right] <\infty.
$
\end{description}
\end{lemma}

\proof 
(i) Note that
$$
\Gamma_n(\theta^*)[u,u] = \tilde{\Gamma}_n[u,u] + \tilde{M}_n[u,u] + \tilde{R}_n[u,u],
$$
where
\begin{eqnarray*}
\tilde{\Gamma}_n[u,u] &=& \frac{1}{2 n} \sumkn \left\{ \partial_\theta^2 \log \det S(\xt,\theta^*)[u^{\otimes 2}]
+ \partial_\theta^2 S^{-1}(\xt,\theta^*) [u^{\otimes 2}, S(\xt,\theta^*) ] \right\}, \\
\tilde{M}_n[u,u] &=& \frac{1}{2 n} \sumkn  \partial_\theta^2 S^{-1}(\xt,\theta^*) 
[u^{\otimes 2}, h^{-1} (\sigma(\xt,\theta^*) \Delta_kw)^{\otimes 2} - S(\xt,\theta^*)], \\
\tilde{R}_n[u,u] &=& \frac{1}{2 n} \sumkn \partial_\theta^2 S^{-1}(\xt,\theta^*) 
[u^{\otimes 2}, h^{-1} (\Delta_k Y)^{\otimes 2} - h^{-1} (\sigma(\xt,\theta^*) \Delta_kw)^{\otimes 2}].
\end{eqnarray*}
By Lemmas \ref{lem4-1} and \ref{lem4-2}, for every $q>{\cred \sf 1}$,
\begin{eqnarray*}
\sup_{n \in {\mathbb N}} ||  \sqrt{n} ( \tilde{\Gamma}_n -\Gamma(\theta^*) ) ||_q &<& \infty, 
\quad \mbox{and} \quad 
\sup_{n \in {\mathbb N}} ||   \sqrt{n} \tilde{R}_n   ||_q < \infty.
\end{eqnarray*}
Using Burkholder's inequality and Sobolev's inequality, one has
$$
\sup_{n \in {\mathbb N}} ||   \sqrt{n} \tilde{M}_n   ||_q < \infty
$$
{\cred for $q>{\sf p}$,} and the desired inequality is obtained.

(ii) Since
\begin{eqnarray*}
\frac{1}{n} \partial_{\theta}^{\cred j} {\mathbb H}_n(\theta) &=& -\frac{1}{2 n} \sumkn \left\{ \partial_{\theta}^{\cred j}
 \log \det S(\xt, \theta) 
+ h^{-1} \partial_{\theta}^{\cred j} S^{-1}(\xt,\theta) [(\Delta_k Y)^{\otimes 2}] \right\}
\end{eqnarray*}
{\cred for j=3,4,} 
it is easy to show that
$$
\sup_{n \in {\mathbb N}} ||  \sup_{\theta \in \Theta} n^{-1} | \partial_{\theta}^3 {\mathbb H}_n(\theta) | ||_q < \infty
$$
{\cred for $q>{\sf p}$}, 
which completes the proof.
\qed

\vspace{0.5cm}

{\bf Proof of Theorem \ref{thm1}}. It is enough to check the regularity conditions 
$[A1'']$, $[A2]$, $[A3]$, $[A4']$, $[A5]$ and $[A6]$
in Theorem 2 of Yoshida (2005).
It follows from [H2] that [A2], [A3] with $\rho =2$ and [A5] are satisfied for every $L>0$.
We can take appropriate parameters satisfying $[A4']$, 
and Lemma \ref{lem2} implies [A1''] for every $L>0$.
Lemma \ref{lem1} yields [A6] for every $L>0$.
This completes the proof.
\qed

\vspace{0.5cm}
We denote $B(R) = \{ u \in \bbR^{{\sf p}} \ ; \ |u| \leq R \}$.

\begin{lemma} \label{lem3} 
Assume $[H1]$. Then,
for every $q>p$ and $R>0$, there exists  $C_0 >0$ such that
$$
\sup_{n \in {\bf N}} E \left[ \left| \log {\mathbb Z}_n (u)  \right|^q \right] \leq C_0 |u|^q
$$
for all $u \in B(R)$.
\end{lemma}

\proof 
$\log {\mathbb Z}_n(u)$ has a decomposition
\begin{eqnarray*}
\log {\mathbb Z}_n(u) &=&
\Delta_n[u] - \frac{1}{2} \Gamma_n(\theta^*)[u,u] + \tilde{r}_n(u),
\end{eqnarray*}
where
$$
\tilde{r}_n(u) = \frac{1}{2 n^{3/2}} \int_0^1 (1-s)^2  \partial_\theta^3 {\mathbb H}_n(\theta^* + s(1/\sqrt{n})u)[u^{\otimes 3}] ds.
$$
By Lemma \ref{lem1},
$\sup_{n \in {\mathbb N}} || \Delta_n[u] ||_q^q < C_1 |u|^q$. 
Lemma \ref{lem2} yields that
$\sup_{n \in {\mathbb N}} || \tilde{r}_n(u) ||_q^q < C_2 |u|^{3q}$.
Moreover, 
$\sup_{n \in {\mathbb N}} || \Gamma_n(\theta^*)[u,u] ||_q^q < C_3 |u|^{2q}$.
Thus, noting that $|u|^{2q} + |u|^{3q} \leq C_4 |u|^{q}$ for all $u \in B(R)$,
we obtain the desired inequality.
\qed

\begin{remark}\rm
By using Lemma 2 in Yoshida (2005), Lemma \ref{lem3} implies that for every $R>0$,
$$
\sup_{n \in {\bf N}} E \left[ \left( \int_{|u| \leq R}  {\mathbb Z}_n (u)  du  \right)^{-1} \right] <\infty.
$$
Consequently, one has that
$$
\sup_{n \in {\bf N}} E \left[ \left( \int_{{\mathbb U}_n}  {\mathbb Z}_n (u)  du  \right)^{-1} \right] <\infty.
$$
\end{remark}

\vspace{0.5cm}

\begin{lemma} \label{lem3-2} 
Assume ${\cred [H1^\sharp]}$. Then,
for every $R>0$, 
$ {\mathbb Z}_n (u) \rightarrow^{d_s({\cal F})} {\mathbb Z}(u)$
in $C(B(R))$ as $n \rightarrow \infty$.
\end{lemma}

\proof 
%
\begin{en-text}
By the localization method, it is enough to show the result under
the following stronger assumption $[H1'']$: 
\\
$[H1'']$ Assumption $[H1']$ holds and
there exists $K>0$ such that
$
|b_t| + |X_t| + |\tilde{b}_t| + |a_t| + |\tilde{a}_t| < K$.
\end{en-text}
Let $\sigma_t^* =\sigma(X_t,\theta^*)$, 
$\Delta_k \tilde{X} = \sigma_{t_{k-1}}^* \Delta_kw$ 
and 
\begin{eqnarray*}
H(y,u) &=& 
 -\frac{1}{2 \sqrt{n}}  \left\{ (\partial_\theta \log \det S(\xt, \theta^*) )[u] 
+ h^{-1} (\partial_\theta S^{-1})(\xt,\theta^*) [u, y^{\otimes 2}] \right\} \\
& & - \frac{1}{4 n} \left\{ (\partial_\theta^2 \log \det S(\xt, \theta^*) )[u^{\otimes 2}] 
+ h^{-1} (\partial_\theta^2 S^{-1})(\xt,\theta^*) )[u^{\otimes 2}, y^{\otimes 2}] \right\}
\end{eqnarray*}
for $y \in \bbR^{\sf d}$. Note that
\begin{eqnarray*}
\log {\mathbb Z}_n(u) &=& \sum_{k=1}^n \left\{ H(\Delta_kY,u) - H(\Delta_k \tilde{X},u) \right\} 
+ \sum_{k=1}^n H(\Delta_k \tilde{X},u) + \tilde{r}_n(u) \\
&=& \sum_{k=1}^n \left\{ J_{1,k}(u) + J_{2,k}(u) \right\} + \sum_{k=1}^n H(\Delta_k \tilde{X},u) + \tilde{r}_n(u),
\end{eqnarray*}
where
\begin{eqnarray*}
J_{1,k}(u) &=& -\frac{\sqrt{n}}{2 T} (\partial_\theta S^{-1})(\xt,\theta^*) [u, (\Delta_kY)^{\otimes 2} -(\Delta_k \tilde{X})^{\otimes 2}],  \\
J_{2,k}(u) &=& - \frac{1}{4 T} 
(\partial_\theta^2 S^{-1})(\xt,\theta^*) )[u^{\otimes 2}, (\Delta_kY)^{\otimes 2} -(\Delta_k \tilde{X})^{\otimes 2}].
\end{eqnarray*}
Since it follows from Lemma \ref{lem2} that $ ||\sqrt{n} \tilde{r}_n(u) ||_q^q < \infty$, 
one has that $\tilde{r}_n(u) \rightarrow^p 0$ as $n \rightarrow \infty$.
Lemma \ref{lem4-2} yields that $\sum_{k=1}^n J_{2,k}(u) \rightarrow^p 0$ as $n \rightarrow \infty$.

An easy calculation yields that
\begin{eqnarray*}
& & (\Delta_k Y) (\Delta_k Y)^\star 
-(\Delta_k \tilde{X}) (\Delta_k \tilde{X})^\star \\
&=& (\int_{t_{k-1}}^{t_k} b_t dt + \int_{t_{k-1}}^{t_k} (\sigma_t^* -\sigma_{t_{k-1}}^* ) dw_t)
(\int_{t_{k-1}}^{t_k} b_t dt + \int_{t_{k-1}}^{t_k} (\sigma_t^* -\sigma_{t_{k-1}}^* ) dw_t)^\star \\
& & + (\int_{t_{k-1}}^{t_k} \sigma_{t_{k-1}}^* dw_t)
(\int_{t_{k-1}}^{t_k} b_t dt + \int_{t_{k-1}}^{t_k} (\sigma_t^* -\sigma_{t_{k-1}}^*  ) dw_t)^\star \\
& & + (\int_{t_{k-1}}^{t_k} b_t dt + \int_{t_{k-1}}^{t_k} (\sigma_t^* -\sigma_{t_{k-1}}^* ) dw_t)
(\int_{t_{k-1}}^{t_k} \sigma_{t_{k-1}}^* dw_t)^\star.
\end{eqnarray*}
Let $f_t =  (\partial_\theta S^{-1})(\xt,\theta^*)[u]$ and $\bar{f}_t = f_{t_{k-1}}$ if $t \in [t_{k-1}, t_k)$.
Furthermore, setting 
\begin{eqnarray*}
H_{1,k} &=& -\frac{\sqrt{n}}{2 T}  
f_{t_{k-1}} \left[ \left( \int_{t_{k-1}}^{t_k} b_t dt + \int_{t_{k-1}}^{t_k} (\sigma_t^* -\sigma_{t_{k-1}}^* ) dw_t \right)^{\otimes 2}
\right], \\
H_{2,k} &=& -\frac{\sqrt{n}}{2 T} 
f_{t_{k-1}} \left[ (\int_{t_{k-1}}^{t_k} b_t dt ) (\int_{t_{k-1}}^{t_k}  \sigma_{t_{k-1}}^*  dw_t)^\star \right], \\
H_{3,k} &=& -\frac{\sqrt{n}}{2 T} 
f_{t_{k-1}} \left[ (\int_{t_{k-1}}^{t_k} \sigma_{t_{k-1}}^* dw_t) (\int_{t_{k-1}}^{t_k} (\sigma_t^* -\sigma_{t_{k-1}}^*  ) dw_t)^\star
\right],
\end{eqnarray*}
one has that $J_{1,k}(u) = H_{1,k} + 2 H_{2,k} + 2 H_{3,k} $. Since
{\cred 
$$
\left\| \sum_{k=1}^n H_{1,k} \right\|_1
\leq C \sqrt{n} \sum_{k=1}^n \left\| \left|
\int_{t_{k-1}}^{t_k} b_t dt + \int_{t_{k-1}}^{t_k} (\sigma_t^* -\sigma_{t_{k-1}}^* ) dw_t
\right|^2 \right\|_2
\leq C_1 h,
$$
}
we obtain $\sum_{k=1}^n H_{1,k} \rightarrow^p 0$.
\begin{en-text}
Noting that by It{\^o}'s formula,

$$
\sum_{k=1}^n H_{2,k} = - \frac{\sqrt{n}}{2T} \int_0^T \left\{ \bar{\sigma}_t^* (w_t -\bar{w}_t ) \right\}^\star \bar{f}_t b_t dt
- \frac{\sqrt{n}}{2T}  \int_0^T \left\{ \bar{f}_t (B_t - \bar{B}_t) \right\}^\star \bar{\sigma}_t^* dw_t,
$$
where $B_t = \int_0^t b_sds$,
and that by standard estimates,
\begin{eqnarray*}
& & \sqrt{n} \int_0^T \left\{ \bar{\sigma}_t^* (w_t -\bar{w}_t ) \right\}^\star \bar{f}_t b_t dt \rightarrow^p 0, \\ 
& & E \left[ \left( \sqrt{n} \int_0^T \left\{ \bar{f}_t (B_t - \bar{B}_t) \right\}^\star \bar{\sigma}_t^* dw_t \right)^2 \right]
\leq C h \rightarrow 0,
\end{eqnarray*}
one has $\sum_{k=1}^n H_{2,k} \rightarrow^p 0$.

\end{en-text}
{\cred 
In the decomposition 
\beas 
\sum_{k=1}^n H_{2,k} 
&=&
 - \frac{\sqrt{n}}{2T}\sum_{k=1}^n f_{t_{k-1}}
 \bigg[\int_{t_{k-1}}^{t_k}(b_t-b_{t_{k-1}})dt\bigg(\int_{t_{k-1}}^{t_k}\sigma_{t_{k-1}}^*dw_t\bigg)^\star\bigg]
\\&&
 - \frac{\sqrt{n}}{2T}\sum_{k=1}^n f_{t_{k-1}}
 \bigg[b_{t_{k-1}}h\bigg(\int_{t_{k-1}}^{t_k}\sigma_{t_{k-1}}^*dw_t\bigg)^\star\bigg],
 \eeas
 the first term on the right-hand side is $o_p(1)$ by standard estimates, 
\footnote{{\cred For a progressively measurable process $b_t$ satisfying $E[\int_0^T|b_t|^2dt]<\infty$, 
for any $\ep>0$, one can find a right-continuous adapted process $(\beta_t)$ such that 
$E[\int_0^T |b_t-\beta_t|^2dt]<\ep$. Then we may assume the right-continuity of $(b_t)$.}}
and the second term is $O(n^{-1/2})$ in $L^2$-norm. Thus $\sum_{k=1}^nH_{2,k}\to^p0$. 
}

In order to obtain that $\sum_{k=1}^n H_{3,k} \rightarrow^p 0$ as $n \rightarrow \infty$,
it is enough to show that as $n \rightarrow \infty$,
\begin{eqnarray*}
\sum_{k=1}^n E[  H_{3,k} | {\cal F}_{t_{k-1}} ] \rightarrow^p 0, \quad 
\sum_{k=1}^n E[  (H_{3,k})^2 | {\cal F}_{t_{k-1}} ] \rightarrow^p 0 
\end{eqnarray*}
because of Lemma 9 in Genon-Catalot and Jacod (1993).
Since
\begin{eqnarray*}
E[  (H_{3,k})^2 | {\cal F}_{t_{k-1}} ] &\leq& C_1 n 
E \left[ \left| \int_{t_{k-1}}^{t_k} \sigma_{t_{k-1}}^* dw_t \right|^4 | {\cal F}_{t_{k-1}} \right]^{1/2}
E \left[ \left| \int_{t_{k-1}}^{t_k} (\sigma_t^* -\sigma_{t_{k-1}}^*  ) dw_t \right|^4 | {\cal F}_{t_{k-1}} \right]^{1/2} \\
&\leq& C_2 h^2,
\end{eqnarray*}
one has $\sum_{k=1}^n E[  (H_{3,k})^2 | {\cal F}_{t_{k-1}} ] \rightarrow^p 0$.
Since It{\^o}'s formula implies that 
\begin{eqnarray*}
\sigma_t^* -\sigma_{t_{k-1}}^* &=& 
\int_{t_{k-1}}^t (\sum_{i=1}^d \partial_{x_i} \sigma(X_s, \theta^*) \tilde{b}_s^i 
+ \sum_{i,j=1}^d \frac{1}{2} \partial_{x_i} \partial_{x_j} \sigma(X_s, \theta^*) 
((a_s a_s^\star )^{ij} + (\tilde{a}_s \tilde{a}_s^\star )^{ij}) )ds \\
& & + \int_{t_{k-1}}^t  \sum_{i=1}^d \sum_{j=1}^r  
 \partial_{x_i} \sigma(X_s, \theta^*)  a_s^{i j} dw_s^j 
 + \int_{t_{k-1}}^t  \sum_{i=1}^d \sum_{j=1}^{\tilde{r}}  
 \partial_{x_i} \sigma(X_s, \theta^*)  \tilde{a}_s^{i j} d\tilde{w}_s^j, 
\end{eqnarray*}
one has 
\begin{eqnarray*}
\left| E [ H_{3,k} | {\cal F}_{t_{k-1}} ] \right| 
&\leq& C_1 \sqrt{n} \left|  E \left[ 
\mbox{Tr} \left( f_{t_{k-1}}  \sigma_{t_{k-1}}^*  \int_{t_{k-1}}^{t_k}(\sigma_t^* -\sigma_{t_{k-1}}^* )^\star dt \right)  
| {\cal F}_{t_{k-1}} \right] \right| \\
&\leq& C_2 h^{3/2}
\end{eqnarray*} 
and 
$\sum_{k=1}^n E [ H_{3,k} | {\cal F}_{t_{k-1}} ] \rightarrow^p 0$ as $n \rightarrow \infty$.
{\cred Therefore $\sum_{k=1}^nH_{3,k}\to^p0$, and hence $\sum_{k=1}^nJ_{1,k}(u)\to^p0$. }

Let $\xi_k := \xi_k(u) = H(\Delta_k \tilde{X},u)$ and 
\begin{eqnarray*}
\Gamma_t(\theta^*)[u,u] &=& \frac{1}{2 T}
\int_{0}^{t}  
\mbox{Tr} \left( (\partial_\theta S) S^{-1}
(\partial_\theta S) S^{-1}(X_s,\theta^*)[u^{\otimes 2}] \right) ds.
\end{eqnarray*}
It follows from Theorem 3-2 of Jacod (1997) that 
\begin{eqnarray}
\log {\mathbb Z}_n(u) &\rightarrow& \Gamma(\theta^*)^{1/2} \zeta[u] -\frac{1}{2} \Gamma(\theta^*)[u,u] 
\quad \mbox{${\cal F}$-stable} \label{LAMN}
\end{eqnarray}
if we show that
\begin{eqnarray}
& & \sup_{t \in [0,T]} \left| \sum_{k=1}^{[nt/T]} E \left[ \xi_k | {\cal F}_{t_{k-1}} \right] + \frac{1}{2}\Gamma_t(\theta^*)[u,u] \right|
\rightarrow^p 0, \label{J-1} \\
& & \sum_{k=1}^{[nt/T]} \left\{ E \left[ \xi_k^2 | {\cal F}_{t_{k-1}} \right] 
- (E \left[ \xi_k | {\cal G}_{t_{k-1}} \right])^2 \right\}
\rightarrow^p \Gamma_t(\theta^*)[u,u], \quad \mbox{for all $t \in [0,T]$}, \label{J-2}  
\\
& & \sum_{k=1}^{[nt/T]}  E \left[ \xi_k \Delta_kw | {\cal F}_{t_{k-1}} \right] 
\rightarrow^p 0, \quad \mbox{for all} \ t \in [0,T], \label{J-3} \\
& & \sum_{k=1}^{[nt/T]}  E \left[ \xi_k \Delta_kN | {\cal F}_{t_{k-1}} \right] 
\rightarrow^p 0, \quad \mbox{for all $t \in [0,T]$ and $N \in {\cal M}_b(w^\perp)$},  \label{J-3-b} \\
& & \sum_{k=1}^{n}  E \left[ \xi_k^4 | {\cal F}_{t_{k-1}} \right] 
\rightarrow^p 0, \label{J-4} 
\end{eqnarray}
where ${\cal M}_b(w^\perp)$ is the class of all bounded ${\bf F}$-martingales which is orthogonal to $w$.


{\cred The symbol $R(r_n)$ denotes a sequence of random variables for which 
$\|R(r_n)\|_q\leq C_qr_n$ for every $q>1$, where $C_q$ depend on neither $n$ nor other variables. }
Indeed, one has
\begin{eqnarray*}
E \left[ \xi_k | {\cal F}_{t_{k-1}} \right]
&=& 
 -\frac{1}{2 \sqrt{n}}  \left\{ (\partial_\theta \log \det S(\xt, \theta^*) )[u] 
+  (\partial_\theta S^{-1})(\xt,\theta^*) [u, S(\xt,\theta^*) ] \right\} \\
& & - \frac{1}{4 n} \left\{ (\partial_\theta^2 \log \det S(\xt, \theta^*) )[u^{\otimes 2}] 
+  (\partial_\theta^2 S^{-1})(\xt,\theta^*) )[u^{\otimes 2},  S(\xt,\theta^*)] \right\}, \\
&=& - \frac{1}{4 n} \mbox{Tr} \left( (\partial_\theta S) S^{-1}
(\partial_\theta S) S^{-1}(\xt ,\theta^*)[u^{\otimes 2}] \right), \\
\eeas
\beas
E \left[ \xi_k^2 | {\cal F}_{t_{k-1}} \right]
&=&  \frac{1}{4 n} \left\{ ((\partial_\theta \log \det S(\xt, \theta^*) )[u])^2 \right. \\
& & 
+ 2 (\partial_\theta \log \det S(\xt, \theta^*) )[u] (\partial_\theta S^{-1})(\xt,\theta^*) [u, S(\xt,\theta^*) ] \\
& & +  ((\partial_\theta S^{-1})(\xt,\theta^*) [u, S(\xt,\theta^*)] )^2
+ 2 \mbox{Tr} \left( (\partial_\theta S) S^{-1}
(\partial_\theta S) S^{-1}(\xt ,\theta^*)[u^{\otimes 2}] \right) 
 {\cred+  R(h^{3/2})}\\
&=& \frac{1}{2 n} \mbox{Tr} \left( (\partial_\theta S) S^{-1}
(\partial_\theta S) S^{-1}(\xt ,\theta^*)[u^{\otimes 2}] \right)
+  R(h^{3/2}).
\end{eqnarray*}
\begin{eqnarray*}
E \left[ \xi_k \Delta_kw | {\cal F}_{t_{k-1}} \right] &=& 
- \frac{\sqrt{n}}{2 T}  
E \left[ (\partial_\theta S^{-1})(\xt,\theta^*) [u, (\sigma(\xt, \theta^*) \Delta_kw)^{\otimes 2}] \Delta_kw \ | \ {\cal F}_{t_{k-1}} \right] \\
& & - \frac{1}{4 T}  
E \left[ (\partial_\theta^2 S^{-1})(\xt,\theta^*) [u^{\otimes 2}, (\sigma(\xt, \theta^*) \Delta_kw)^{\otimes 2}] 
\Delta_kw \ | \ {\cal F}_{t_{k-1}} \right] \\
&=& 0, 
\end{eqnarray*}
\begin{eqnarray*}
E \left[ \xi_k \Delta_kN | {\cal F}_{t_{k-1}} \right] &=& 
- \frac{\sqrt{n}}{2 T}  
E \left[ (\partial_\theta S^{-1})(\xt,\theta^*) [u, (\sigma(\xt, \theta^*) \Delta_kw)^{\otimes 2}] \Delta_kN \ | \ {\cal F}_{t_{k-1}} \right] \\
& & - \frac{1}{4 T}  
E \left[ (\partial_\theta^2 S^{-1})(\xt,\theta^*) [u^{\otimes 2}, (\sigma(\xt, \theta^*) \Delta_kw)^{\otimes 2}] 
\Delta_kN \ | \ {\cal F}_{t_{k-1}} \right] \\
&=& 0,
\end{eqnarray*}
where in the last estimate, we note that
$E[ (\Delta_k w^i)^2 \Delta_k N \ | \ {\cal F}_{t_{k-1}} ] 
= E[ ((\Delta_k w^i)^2 -h) \Delta_k N \ | \ {\cal F}_{t_{k-1}} ]$. 
Therefore, (\ref{J-1})-(\ref{J-3-b}) are shown.
Furthermore,
$E \left[ \xi_k^4 | {\cal F}_{t_{k-1}} \right]{\cred =R(h^2)}$, 
{\cred which proves} (\ref{LAMN}).
%
Moreover, it is easy to see that
the joint convergence for finitely many $u$'s is valid in (\ref{LAMN}). 

In order to show the tightness of the family $\{ \log {\mathbb Z}_n(u)|_{B(R)} \ ; \ n \in {\mathbb N} \}$
for every $R>0$,
it is enough to prove that
$\sup_{n \in {\mathbb N}}  E \left[ \left( \sup_{u \in {B(R)}} \left| \partial_u \log {\mathbb Z}_n(u) \right| \right)^q \right] < \infty$
for $q>1$.
For details, see Yoshida (1990).
In the same way as in the proof of Lemma \ref{lem1}, using Lemma \ref{lem4-2}, one has that 
$\sup_{n \in {\mathbb N}} \sup_{u \in {B(R)}} \left|  \left| \partial_u \log {\mathbb Z}_n(u) \right| \right|_q < \infty$
and 
$\sup_{n \in {\mathbb N}} \sup_{u \in {B(R)}} E \left|  \left|  \partial_u^2 \log {\mathbb Z}_n(u) \right| \right|_q < \infty$.
This completes the proof.
\qed

\vspace{0.5cm}

{\bf Proof of Theorem \ref{thm3}}. Note that Theorem \ref{thm1} implies Proposition 2 in Yoshida (2005).
Using Lemma \ref{lem3-2} together with this fact, one can apply Theorem 5 in Yoshida (2005),
which completes the proof.
\qed

\vspace{0.5cm}

{\bf Proof of Theorem \ref{thm2}}. By Theorem \ref{thm1} and Lemma \ref{lem3-2},  
the regularity conditions of Theorem 6 in Yoshida (2005) are satisfied
and one has that 
$$
\int f(u) {\mathbb Z}_n(u) du \rightarrow^d \int f(u) {\mathbb Z}(u) du
$$
as $n \rightarrow \infty$, for all continuous functions $f$ of at most polynomial growth.
It follows from Lemma 2 in Yoshida (2005) that
Lemma \ref{lem3} meets the condition (i) of Theorem 8 in Yoshida (2005).
Note that the condition (ii) of Theorem 8 in Yoshida (2005) is satisfied in our setting.
This completes the proof.
\qed

{\cod 
\section{{\cred Some remarks on the analytic criteria}}
In this paper, we derived Theorems \ref{230614-1} and \ref{230614-2} 
thought [N$_0$] by checking [N$_1$]. 
Because of [N$_1$], we assumed existence of supporting functions, 
this condition obviously works for the one-dimensional $X$, and 
{\cred also does in a multi-dimensional case if the null set of $f$ is locally a regular submanifold. 
However, this condition is not completely general. 
}

The supporting function supports $f$ in a neighborhood of a point $(x_0,\theta)\in\calx_0\times\Theta$. 
However, what was necessary in our logic was the existence of a function supporting $f$ 
in a particular sector depending on the location of $X_s$ 
of a good increment $X_{\tau_j}-X_s$, as was seen in the proof of 
Proposition \ref{230608-1}. 
As a matter of fact, the sectors can be chosen discontinuously 
though it was done continuously under [N$_1$]. It is because Condition [N$_1$] is a topological condition and 
this nature was used in Lemma \ref{210112-1}.

On the other hand, Condition [N$_0$] (or [N$_0^\flat$]) does not require continuity of the supporting function 
$G_{\ell,k}(x_0,x,\theta,\xi)$ in $(x_0,\theta)$ associated with a sector. 
Besides, the double sided sector is not necessary for nondegenerate diffusions; 
of course, there are cases in which it really works effectively. 
Therefore, Theorems \ref{thm3} and \ref{thm2} together with  [N$_0^\flat$] 
are much stronger than Theorems \ref{230614-1} and \ref{230614-2}. 
Though we do not go into details here, a more general condition will be as follows: 
There exist $t_0>0$, $a\in(0,1)$, $\ep>0$ and a finite subset $\Xi$ of $\bbS$ such that 
for any $(x_0,\theta)$, for some $\xi_0\in\Xi$, for all $t\in(0,t_0)$, 
\beas 
\inf_{\xi\in D(\xi_0,\ep)}|f(x_0+t\xi,\theta)|
&>&
a |f(x_0+t\xi_0,\theta)|
\eeas
and 
$\max_{j=0,...,J-1}|c_j(x_0,\theta)|>0$ for each $(x_0,\theta)$, where $c_j(x_0,\theta)$ are given by 
the derivatives of $f$ and $\xi_0$. 

}


{\cred 
If the null set of $f $ includes irregular points and if $X_\tau$ hits them, then the criteria 
like $[A3]$ 
do not work in general. 
However, even if the process $X$ starts bad points, if it moves quickly to a good area of regular points, 
it is possible to apply the idea of our criteria by some modification. 
As before, let $\{\calx_\ell\}_{\ell=1,...,\bar{\ell}}$ cover $\calx_0$ and 
$\{\Theta_{\ell,k}\}_{k=1,...,\bar{k}_\ell}$ cover $\Theta$ for each $\ell$. 
We assume that for each $(x_0,\theta)\in\calx_\ell\times\Theta_{\ell,k}$, there are a function $g$ and $\xi_{\ell,k}\in\bbS$ 
such that 
\begin{description}
\item[(i)] $\partial_x^jg$ exist and they are continuous for $j=0,...,J$, 
and $\max_{j=0,...,J-1}|c_j(x_0,\theta)|>0$ for each $(x_0,\theta)$, where $c_j(x_0,\theta)$ are given by 
the derivatives of $g$ and $\xi_0$. 
\item[(ii)] For every $(x_0,\theta)\in (\calx_\ell\cap U_n)\times\Theta_{\ell,k}$, 
$
|f(x,\theta)| \geq |g(P_{\xi_{\ell,k}}x,\theta)|
$
for all $(x,\theta)\in B(x_0,n^{-\beta_0})\times\Theta_{\ell,k}$. 
\end{description}
Moreover suppose that there is a sequence of stopping times $\tau_n$ such that 
$(1-P[\tau_n\leq T_0,\>X_{\tau_n}\in U_n])_{n\in\bbN}\in\calp$ for some $T_0\in[0,T)$. 

If $x_*$ is a singular point, we can take $U_n=\bbR^{\sf d}\setminus \overline{B(x_*,n^{-\beta_1})}$ for 
$\beta_1\in(0,\beta_0)$. 
Then it is possible to prove $[H2]$ for a nondegenerate diffusion process $X$, 
by composing the arguments in the previous sections 
with $\alpha_0>\alpha_1>\cdots>\alpha_J>\beta_0$. 
A simple example is $f(x_1,x_2,\theta)=x_1(x_1^2-\theta x_2^4)$ for $x=(x_1,x_2)$, $x_*=(0,0)$, and 
$\mbox{supp}\call\{X_0\}=\{0\}\times[0,1]$. 
}

}


\section*{References}
\begin{description}{}{
}




\item
Adams, R.\ A.\ (1975).\ 
{\it Sobolev spaces}. 
Pure and Applied Mathematics, Vol. 65. 
Academic Press,
New York-London.

\item
Adams, R.\ A.\ and Fournier, J.\ J.\ F.\ (2003).\ 
{\it Sobolev spaces}. Second edition. 
Pure and Applied Mathematics (Amsterdam), 140. Elsevier/Academic Press, Amsterdam.




\item
Dohnal, G.\
(1987).\ 
On estimating the diffusion coefficient.\ 
{\it J.\ Appl.\ Probab.\ }{\bf 24}, 105--114.




\item
Genon-Catalot, V.\ and Jacod, J.\ 
(1993).\ 
On the estimation of the diffusion coefficient 
for multidimensional diffusion processes.\
{\it Ann.\ Inst.\ Henri Poincar{\'e} Probab.\ Statist.\ }{\bf 29}, 119--151. 

\item
Genon-Catalot, V.\ and Jacod, J.\ 
(1994).\ 
Estimation of the diffusion coefficient for diffusion processes: random sampling. 
{\it Scand.\ J.\ Statist.\ }{\bf 21} (1994), 193--221.


\item
Gobet, E.\ 
(2001).\
Local asymptotic mixed normality property for elliptic diffusion: a Malliavin calculus approach.\ 
{\it Bernoulli}, {\bf 7}, 899--912.



\item
Ibragimov, I.\ A.\ and Has'minskii, R.\ Z.\ (1972).\ 
The asymptotic behavior of certain statistical estimates in the smooth case. I. 
Investigation of the likelihood ratio. (Russian) Teor. Verojatnost. i Primenen. 17, 469--486. 

\item
Ibragimov, I.\ A.\ and Has'minskii, R.\ Z.\ (1973).\ 
Asymptotic behavior of certain statistical estimates. II. 
Limit theorems for a posteriori density and for Bayesian estimates. (Russian) 
Teor. Verojatnost. i Primenen. 18, 78--93.

\item
Ibragimov, I.\ A.\ and Has'minskii, R.\ Z.\ 
(1981).\ 
{\it Statistical estimation.\ }Springer Verlag, New York. 


\item
Jacod, J.\
(1997).\
On continuous conditional Gaussian martingales and stable convergence in law.\
Seminaire de Probabilites, XXXI, 232--246, Lecture Notes in Math., 1655, Springer, Berlin.




\item
Kessler, M.\ 
(1997).\ 
Estimation of an ergodic diffusion from discrete observations.\
{\it Scand.\ J.\ Statist.\ }{\bf 24}, 211--229.







\item
Kutoyants, Yu.\ A.\ 
(1984).\ 
Parameter estimation for stochastic processes.\ 
Translated from the Russian and edited by B.\ L.\ S.\ Prakasa Rao.\ 
Research and Exposition in Mathematics, 6. Heldermann Verlag, Berlin.

\item
Kutoyants, Yu.\ A.\ 
(1994).\ 
{\it Identification of dynamical systems with small noise.\ }
Kluwer, Dordrecht.\

\item
Kutoyants, Yu.\ A.\ 
(1998).\
Statistical inference for spatial Poisson processes. 
Lecture Notes in Statistics, 134. Springer-Verlag, New York. 

\item
Kutoyants, Yu.\ A.\ 
(2004).\ 
{\it Statistical inference for ergodic diffusion processes.\ }
Springer-Verlag, London.\


\item
Le Cam, L.\
(1986).\
{\it Asymptotic methods in statistical decision theory}.\ 
Springer-Verlag, New York.

\item
Le Cam, L.\ and Yang, G.\ L.\
(1990).\
{\it Asymptotics in statistics}.\ 
Springer-Verlag, New York.


\item
Masuda, H.\
(2010).\
Approximate self-weighted LAD estimation of discretely observed ergodic Ornstein-Uhlenbeck processes. 
{\it Electronic Journal of Statistics}{\bf 4}, 525--565.





\item
Ogihara, T.\ and Yoshida, N.\
(2009).\
Quasi-likelihood analysis for the stochastic differential equation with jumps. 
To appear in {\it Statistical Inference for Stochastic Processes}.




\item
Prakasa Rao, B.\ L.\ S.\ 
(1983).\
Asymptotic theory for nonlinear least squares estimator for diffusion processes. 
{\it Math.\ Operationsforsch.\ Statist.\ Ser.\ Statist.\ }{\bf 14}, 195--209.

\item
Prakasa Rao, B.\ L.\ S.\ 
(1988).\
Statistical inference from sampled data for stochastic processes.\
{\it Contemp.\ Math.\ }{\bf 80}, 249--284.\ 
Amer.\ Math.\ Soc., Providence, RI.




\item
Shimizu, Y.\
(2006).\
Density estimation of Levy measures for discretely observed diffusion processes with jumps.\
{\it J. Japan Statist. Soc.\ }{\bf 36}, 37--62.

\item
Shimizu, Y.\ and Yoshida, N.\
(2006).\
Estimation of parameters for diffusion processes with jumps from discrete observations.\ 
{\it Statist. Infer. Stoch. Proc.\ }{\bf 9}, 227--277.



\item
S{\o}rensen, M.\ and Uchida, M.\ 
(2003).\
Small diffusion asymptotics for discretely sampled 
stochastic differential equations.\
{\it Bernoulli }{\bf 9}, 1051--1069. 



\item
Uchida, M.\
(2003).\
Estimation for dynamical systems with small noise from discrete observations.\ {\it J.\ Japan Statist.\ Soc.\ }{\bf 33}, 157--167. 

\item
Uchida, M.\
(2004).\ 
Estimation for discretely observed small diffusions 
based on approximate martingale estimating functions. 
{\it Scand. J. Statist.\ }{\bf 31}, 553-566. 







\item
Uchida, M.\ 
(2010).\ 
Contrast-based information criterion for ergodic diffusion processes from discrete observations.\ 
{\it Ann.\ Inst.\ Statist.\ Math.}, {\bf 62}, 161--187.

\item
Yoshida, N.\ 
(1990).\ 
Asymptotic behavior of M-estimator and related random field
for diffusion process.\
{\it Ann.\ Inst.\ Statist.\ Math.\ }
{\bf 42}, 221--251. 





\item
Yoshida, N.\
(2005).\
Polynomial type large deviation inequalities and quasi-likelihood analysis for stochastic differential equations.\
(to appear in {\it Ann.\ Inst.\ Statist.\ Math.})

\begin{en-text}
\item
Yoshida, N.\
(2006).\
Polynomial type large deviation inequalities and convergence of statistical random fields. 
ISM Research Memorandum 1021, Institute of Statistical Mathematics. 
\end{en-text}

\item
Yoshida, N.\
(2011).\
Polynomial type large deviation inequalities and quasi-likelihood analysis for stochastic differential equations.\
{\it Ann.\ Inst.\ Statist.\ Math.}, {\bf 63}, 431--479.


\end{description}


\end{document}